\documentclass[english,a4paper,onecolumn,12pt]{article}
\usepackage[margin=1.1in,footskip=0.25in]{geometry}
\usepackage{amssymb}
\usepackage{amsfonts}
\usepackage{babel}
\usepackage{eucal}
\usepackage{amsmath}
\usepackage{algpseudocode}
\usepackage{algorithmicx}
\usepackage{mathptmx}
\usepackage{helvet}
\usepackage{courier}
\usepackage{type1cm}         
\usepackage{makeidx}         % allows index generation
\usepackage{graphicx}        % standard LaTeX graphics tool
\usepackage{multicol}        % used for the two-column index
\usepackage[bottom]{footmisc}% places footnotes at page bottom
\newcommand{\bx}{\bar x}

\newcommand{\tx}{\tilde x}

\newcommand{\tu}{\tilde u}

\newcommand{\equal}{\buildrel \Delta \over =}
\newcommand {\bmat} {\left[\begin{array} }
	\newcommand {\emat} {\end{array}\right]}
	\newcommand {\ematrix}{\end{array}\right]}
\newcommand{\blista}{\renewcommand{\labelenumi}{(\roman{enumi})}
	\begin{enumerate}}
	\newcommand{\elista}{\end{enumerate} \renewcommand{\labelenumi}{\arabic{enumi}.}}
\newcommand {\beqn}{\begin{equation}}
	\newcommand {\eeqn}{\end{equation}}
\newcommand {\beqna}{\begin{eqnarray}}
	\newcommand {\eeqna}{\end{eqnarray}}
\newcommand {\beqnan}{\begin{eqnarray*}}
	\newcommand {\eeqnan}{\end{eqnarray*}}
\newcommand {\nn}{\nonumber}

\newcommand{\vx} {\mathbf{x}}
\newcommand{\vu} {\mathbf{u}}
\newcommand{\vz} {\mathbf{z}}

\newcommand{\vA} {\mathbf{A}}
\newcommand{\vB} {\mathbf{B}}
\newcommand{\vQ} {\mathbf{Q}}
\newcommand{\vR} {\mathbf{R}}
\newcommand{\vI} {\mathbf{I}}
\newcommand{\x} {\mathrm{x}}

\newcommand{\SetO}{ \mathcal{O}}

\newcommand{\setX}{\mathcal{X}}

\newcommand{\setY}{\mathcal{Y}}
\newcommand{\setU}{\mathcal{U}}
\newcommand{\setZ}{\mathcal{Z}}

\newcommand{\setW}{\mathcal{W}}

\def\R{ {\rm \,I\!R} }

\newcommand{\QED}{\hfill \rule{2mm}{2mm}}

\newtheorem {theorem}{Theorem}

\newtheorem {propt}{Property}

\newtheorem {lem}{Lemma}
\newtheorem {rem}{rem}
\newtheorem {exmp}{Example}

\newtheorem {assumpt}{Assumption}

\newenvironment{pf}{\noindent \textbf{Proof:}}{\QED \vspace{1ex}}
\newenvironment{pf*}[1][Proof:]{\noindent \textbf{#1} }{\QED \vspace{1ex}}

\begin{document}

\title{Internal Report\\ \textbf{Model Predictive Control for \\setpoint tracking}}

\author{by \\ Daniel Limon$^1$, Antonio Ferramosca$^2$, \\ Ignacio Alvarado$^1$, Teodoro Alamo$^1$}

\date{\small\textit{$^1$ Department of Systems Engineering and Automation, University of Seville. Camino de los Descubrimientos s/n, 41092, Seville, Spain. \{dlm,ialvarado,talamo\}@us.es\\
		$^2$Department of Management, Information and Production Enginieering, University of Bergamo. Viale Marconi 5, 24044, Dalmine (BG), Italy. antonio.ferramosca@unibg.it}}
%\author[First]{Daniel Limon} 
%\author[First]{Antonio Ferramosca}
%\author[Second]{Ignacio Alvarado}
%\author[First]{Teodoro Alamo} 

%\address[First]{Department of Management, Information and Production Engineering, University of Bergamo, 24044 Dalmine, Bergamo, Italy}
%\address[Second]{Department of Engineering, Universidad Loyola Andalucía, 41704 Dos Hermanas, Seville, Spain}

\maketitle
\section{Introduction}

The main objective of tracking control is to steer the tracking error, that is the difference between the reference and the output, to zero while the limits of operation of the plant are satisfied. This requires that some assumptions on the evolution of the future values of the reference must be taken into account. Typically a simple evolution of the reference is considered, such as step, ramp or parabolic reference signals.  It is important to notice that the tracking problem considers possible variation on the reference to be tracked, such as steps or variations of the slope of the ramps. Then the tracking control problem is inherently uncertain, since the reference may differ from expected.

In model predictive tracking control, some assumptions on the expected values of the reference must be considered in order to predict the expected evolution of the tracking error. This report will be devoted to the most extended case of tracking control: the setpoint tracking. In this case  it is assumed that the reference will remain constant along the prediction horizon.
Setpoint tracking is a relevant control problem, for instance, in position control systems or in the process industries in which the plant is typically designed to operate at an equilibrium point that maximizes the profit of the plant. Variations on the profit function or in the operation conditions of the plant may lead to changes in the setpoints of the process variables.

Tracking predictive schemes for constant references can be derived from a predictive controller for regulation, and under certain assumptions, closed-loop stability can be guaranteed if the initial state is inside the feasible region of the MPC. However, if the value of the reference is changed, then there is no guarantee that feasibility and stability properties of the resulting control law hold. Specialized predictive controllers have to be designed to deal with this problem \cite{PannocchiaAICHEJ03, MaederAUT10,BemporadTAC97,ChisciIJC03,LimonAUT08,RossiterCTA96}. This report presents the MPC for tracking approach, which ensures recursive feasibility and asymptotic stability of the setpoint when the value of the reference is changed.

\section{Problem Description}\label{prelres}

It is assumed that the system to be controlled is described by a state-space model as follows:
\beqna
x^{+}&=&Ax+Bu \label{sistema_o}\\
y&=&Cx+Du \nn
\eeqna
where $x \in \mathbb{R}^n$ is the current state
of the system, $u \in \mathbb{R}^m$ is the current input and $x^+$ is the successor
state. The solution of this system for a given sequence of control
inputs $\vu$ and initial state $x$ is denoted as
$x(j)=\phi(j;x,\vu)$, $j \in \mathbb I_{\geq 0}$, where $x=\phi(0;x,\vu)$. The objective of the tracking control is that certain output variables of the system $y \in \mathbb R^p$ track the provided setpoint. Notice that this signal may depend on the current value of both the state and the input.

The state of the system and the control input applied at sampling
time $k$ are denoted as $x(k)$ and $u(k)$ respectively. The system
is subject to hard constraints on state and control that must be fulfilled at each sampling time, throughout its evolution:
\beqnan\label{restric_oOPT} (x(k),u(k))&\in& \mathcal{Z} \eeqnan
for all $k\geq0$. $\mathcal{Z}\subset\mathbb{R}^{n+m}$ is a closed
convex polyhedron such that its interior is non-empty and the origin is inside it\footnote{This condition is not limiting and can be fulfilled by means of an appropriate change of variables of the state and input}. It is assumed that every admissible control input is bounded, that is, $Proj_u(\setZ)=\{u: (x,u) \in \setZ\}$ is compact.

% La compacidad del conjunto Z se a\~{n}ade para facilitar la demostraci\'{o}n de estabilidad. Bastar\'{\i}a con que la funci\'{o}n de coste \'{o}ptima fuese localmente acotada o bien que U sea compacto.

The following standing hypothesis is assumed:

\begin{assumpt} \label{assumption1_OPT}
	The pair (A,B) is stabilizable and the state is measured at each
	sampling time.
\end{assumpt}

If the pair $(A,B)$ is not stabilizable then the system cannot be controlled by a feedback controller and therefore the control problem of the plant should be re-studied. On the other hand, if the state of the plant is not measured, then it should be estimated by means of a suitable observer. In this case, the system should be manually steered to a suitable equilibrium point and then the observer should be triggered to start the estimation. Once the estimation error is known to be small, the automatic mode should be operated. This procedure ensures that the resulting control scheme based on the observer works appropriately, and similarly to the case of state feedback control, thanks to the separation principle.

The problem we consider is the design of an MPC controller
$$u(k)=\kappa_N(x(k), y_{sp})$$
such that the resulting controlled system
$$ x(k+1)=A x(k) + B \kappa_N(x(k),y_{sp})$$
is stable, that is, small changes in the state $x(k)$ cause small changes in the subsequent trajectory,  and, if possible, the tracking error asymptotically tends to  zero, i.e.
$$ \lim_{k \rightarrow \infty} \| y(k)-y_{sp}\| =0$$

\section{The set of reachable setpoints} \label{chaMPCT:ReachSP}

The main design requirement of a tracking controller is to asymptotically stabilize the system to an equilibrium point $(x_{sp},u_{sp})$ such that the steady output is the desired setpoint $y_{sp}$. The equilibrium point and the setpoint must satisfy the dynamic model of the plant and besides satisfy the hard constraints of the system. If both conditions are satisfied, then the setpoint is said to be reachable. But, if some of these conditions are not fulfilled, then the tracking control problem fails since it is not possible to stabilize the plant to the required output satisfying the constraints. In this case the setpoint is said unreachable. Therefore, it is important to study the set of equilibrium points and setpoints that are reachable for a tracking controller.

Firstly, the satisfaction of the dynamic model is analyzed: for a given setpoint $y_{sp}$, the associated steady state and input $(x_{sp},u_{sp})$ must satisfy (\ref{sistema_o}), that is,
$$
x_{sp}=A x_{sp} + B u_{sp}
$$
as well as the output equation
$$
y_{sp}=C x_{sp} + D u_{sp}
$$

A primary question to answer is if for any given setpoint $y_{sp}$, there exists an associated equilibrium point of the system $(x_{sp},u_{sp})$.
As proved in \cite[Lemma 1.14]{RawlingsLIB09}, this condition holds if and only if
\begin{equation}
	rank \left ( \bmat{cc} (A - I_n) & B \\ C & D \emat \right) =n+p
	\label{EqMSS}
\end{equation}
Notice that this condition depends on the matrices of dynamic model of the system $(A,B,C,D)$  and  can only  be fulfilled if the number of inputs  is greater than or equal to the number of outputs, i.e. $m\geq p$. Therefore, for every system such that condition (\ref{EqMSS}) is not satisfied, only  a certain set of setpoints can be tracked. This set can be characterized by using a geometric point of view.

The pair $(x_{sp},u_{sp})$ is an equilibrium point if and only if
\begin{equation}
	\bmat{cc} (A - I_n) & B \emat \bmat{c} x_{sp}\\u_{sp} \emat= 0
	\label{EqSS}
\end{equation}
This implies that $(x_{sp},u_{sp})$ must be contained in the null space of the matrix $$[(A-I_n) \quad B]$$ Since the system is assumed to be stabilizable, then the rank of
$[(A - I_n) \quad B]$ is equal to $n$, and then the dimension of the null space is equal to $m$.  The null space is then the subsapce spanned by the columns of a certain matrix $M_z \in \mathbb{R}^{(n+m) \times m}$ such that
\[
[(A-I_n) \quad B] M_z=0.
\]
Notice that $M_z$ is not unique but the null space is unique. The set of outputs is the subspace spanned by the columns of the matrix $$M_y=[C \quad D] M_z \in \mathbb{R}^{p \times m}$$
Thus, this set is also a linear subspace whose dimension is equal to the rank of $M_y$, and then it is lower than or equal to $m$. If the rank of matrix $M_y$ is equal to $p$, then the linear subspace spanned by the columns includes $\mathbb R^p$ and hence every setpoint $y_{sp}$ is included in it. This is the case when condition (\ref{EqMSS}) hold. In other case,  not every setpoint could be reached, but only those  contained in the linear subspace spanned by $M_y$.

According to  the number of inputs and outputs, the following cases can be found:
\begin{itemize}
	\item If the system is square, i.e. $p=m$, and condition (\ref{EqMSS}) holds, then every given setpoint $y_{sp}$ can be tracked, and for this setpoint there exists a unique equilibrium point $(x_{sp},u_{sp})$.
	\item If the system is flat, i.e. $p<m$,  and condition (\ref{EqMSS}) holds, then every given setpoint $y_{sp}$ can be tracked and for this setpoint, there exists a infinite number of equilibrium points $(x_{sp},u_{sp})$ whose output is $y_{sp}$. Therefore for a given setpoint, there exist some degrees of freedom to choose the equilibrium point of the plant, and they should be fixed considering some additional criterion.
	\item If the system is thin, i.e. $p>m$ or the condition (\ref{EqMSS}) does not hold, then only those setpoints that are in the linear subspace of matrix $M_y$ can be reached. The usual way of overcoming this problem is
	re-defining the output signals that should be tracked. Thus,  \emph{new} controlled variables, $y_c \in
	\mathbb R^{p_c}$ with $p_c \leq p$ are taken in such a way that  they are a linear combination of the
	actual outputs, i.e. $y_c=L_c y=L_c C x+L_c D u$ and the condition (\ref{EqMSS}) is satisfied for the output $y_c$.
\end{itemize}

	The set of reachable setpoints is also limited by the constraints on the system state and input. Thus, for a given setpoint $y_{sp}$ the associated equilibrium point must satisfy
	$$
	(x_{sp},u_{sp})\in \setZ
	$$
	Then, the set of setpoints that satisfy both the dynamics and the constraints are
	$$
	\setY_{sp}= \{ y_s \mid x_s=A x_s + B u_s, y_s=C x_s + D u_s,  (x_s,u_s) \in \setZ \}
	$$
	
	Notice that since $\setZ$ is a polytope and the equality constraints are linear, the set of reachable setpoints $\setY_{sp}$ is also a polytope.
	
	$\setY_{sp}$ is the  set of reachable setpoints of the constrained system and plays an important role in the tracking control problem: if the setpoint is reachable, i.e. $y_{sp} \in \setY_{sp}$, then there exists a tracking control law that can steer the system to it. But if the setpoint is not reachable, i.e.  $y_{sp} \not\in \setY_{sp}$, then it is impossible to find a tracking control law capable to steer the system to it satisfying the constraints. In this case, a suitable reachable setpoint should be calculated according to a certain condition.
	
	The presence of constraints may lead to a possible loss of controllability when some of them are active \cite{RaoAIChE99}. In effect, there may exist reachable setpoints $y_{sp}$ that can be asymptotically tracked by the controller, but once the system has converged to them, it cannot leave them due to the active constraints. This fact can be illustrated by means of a simple example.
	
	Consider the first order system given by $x^+=2 x + u$ with $y=x$ and the input constrained to $u \in [-1,1]$.
	An equilibrium point of the system is such that $x_{sp} + u_{sp}=0$. Since $u_{sp} \in [-1,1]$, then $y_{sp}=x_{sp} \in [-1,1]$.
	Then the set of reachable setpoints is $\setY_{sp}= [-1,1]$, and every setpoint $y_{sp} \in \setY_{sp}$ could be reached by a tracking control law.\\
	Consider that the initial state is $x(0)=1$, then $x(1)=2  + u(0) \in [1,3]$ since $u(0) \in [-1,1]$, and $x(1)=1$ only if $u(0)=-1$. At $k=2$ we have that  $x(2)=2 x(1) + u(1)= 4  + 2 u(0) + u(1) \in [1, 7]$, and $x(2)=1$ only if $u(0)=-1$ and $u(1)=-1$. Applying this recursion it can be proved that $x(k) \in [1,2^{k+1}-1]$, where $x(k)=1$ only if $u(j)=-1$ for $j\in [0,k-1]$.
	%Then $x(1)=2  + u(0) \in [1,3]$ and  $x(2)=2 x(1) + u(1) \geq 2 + u(1)$. By recursion it is easy to see that $x(k) \geq 2 + u(k)$. %
	Then the only sequence of control actions that makes that the system does not diverge is $u(k)=-1$. Therefore, the best that a tracking control law can achieve is to maintain the system at the initial state $x(0)=1$ from which the controlled system cannot escape.
	
	A practical method to avoid this possible loss of controlability, is to remove, from the set of reachable setpoints, those setpoints that lie at active constraints. This can be done, redefining the set of reachable setpoints\footnote{With a slight abuse of notation, this definition will be used throughout this book as the set of reachable setpoints and it will be denoted as $\setY_{sp}$}  as follows:
	
	$$
	\setY_{sp} = \{ y_s \mid x_s=A x_s + B u_s, y_s=C x_s + D u_s, (x_s,u_s) \in \lambda \setZ \}
	$$
	where $\lambda \in [0,1) $ is a given constant arbitrarily close to 1. Notice that for every setpoint contained in $\setY_{sp}$, the associated equilibrium point is contained in the interior of $\setZ$ and then the constraints are not active. Notice also that this definition does not imply a reduction of the reachable setpoints, from a practical point of view, since $\lambda$ can be chosen as close to 1 as required.
	
	Analogously,  the set of reachable steady states $\setX_{sp}$, inputs $\setU_{sp}$ and joint state-input $\setZ_{sp}$ can be defined respectively as follows
	\beqna
	\setX_{sp} &=& \{ x_s \mid x_s=A x_s + B u_s,  (x_s,u_s) \in \lambda \setZ \}\\
	\setU_{sp} &=& \{ u_s \mid x_s=A x_s + B u_s,  (x_s,u_s) \in \lambda \setZ \}\\
	\setZ_{sp} &=& \{ (x_s,u_s) \mid x_s=A x_s + B u_s,  (x_s,u_s) \in \lambda \setZ \}.
	\eeqna

	\begin{rem}
		In the design of some tracking control problems, it may result convenient to characterize the set of equilibrium points by the minimum number of variables. Following the reasoning of this section, the minimum number of variables is equal to $m$ and the characterization is given by
		\beqn (x_s,u_s)=M_\theta \theta \eeqn
		where matrix $M_\theta$ is such that $[(A\!-\!I_n) \quad B]M_\theta=0$. $\theta \in \mathbb R^m$ is the vector of parameters that univocally defines an equilibrium point.
	\end{rem}

	\section{MPC for tracking formulation}\label{newMPCform}
	
	%This section is devoted to present the so-called MPC for
	%tracking \cite{LimonAUT08,FerramoscaAUT09,AlvaradoPHD07,FerramoscaPHD11} for the case of square systems  i.e. $p=m$, that satisfies condition (\ref{EqMSS}). This control technique can be also used in the case of non-square systems, but these cases require more involved development and this will shown in the following report.
	
	This section is devoted to present the so-called MPC for
	tracking \cite{LimonAUT08,FerramoscaAUT09,AlvaradoPHD07,FerramoscaPHD11}. As it was commented before, this predictive control scheme is suitable for the setpoint tracking control problem. In the realistic tracking scenario, where the setpoint is changed, the stabilizing design of predictive controller for regulation based on a terminal constraint may lead to stability issues, due to the changes of the setpoint. The following example illustrates these stability issues and introduces the rationale behind the MPC for tracking.

	\begin{exmp}\label{chaMPCT:ejemplo_1}
		Consider the following system:
		\begin{eqnarray*} A=\left[\!\!\begin{array}{ccc}
				1 & 1
				\\ 0 & 1\end{array}\!\!\right],\!\!\quad\!\! B=\left[\!\!\begin{array}{ccc} 0.0 & 0.5 \\
				1.0 & 0.5\end{array}\!\!\right],\!\!\quad\!\!
			C=\left[\!\!\begin{array}{ccc} 1 &
				0\\ 0 & 1\end{array}\!\!\right]\nonumber\end{eqnarray*}
		The system is constrained to $\|x\|_{\infty}\leq5$ and
		$\|u\|_{\infty}\leq 0.5$. Notice that in this case the output is equal to the state. The set of reachable steady states $\setY_{sp}=\setX_{sp}$ (for a value of $\lambda=0.9999$)  is the region depicted in Figure \ref{chaMPCT:di_ec_factibilidad} in dashed line.
		
		In order to control the plant, an MPC strategy with weighting
		matrices $Q=I_2$ and $R=I_2$, and prediction horizon $N=3$ has been designed. The closed-loop stability of this controller is ensured by means of an equality terminal constraint $x(N)=x_{sp}$. This control law has been designed to track the system from the initial state  $x(0)=(0.6,2.3)$ (dot) to the setpoint $x_{sp1}=(4.9,0.245)$ (square). Figure \ref{chaMPCT:di_ec_factibilidad} shows the domain of attraction of the controller $\setX_N(x_{sp1})$ in dotted line. As $x(0)$ is contained in the feasible region, the control law will steer the system to the setpoint $x_{sp1}$ asymptotically.
		
		Assume the scenario where the setpoint changes from $x_{sp1}$ to a new setpoint  $x_{sp2}=(-4.9,0.2)$ (star) once the control law is applied. In Figure \ref{chaMPCT:di_ec_factibilidad} the dash-dotted line represents $\setX_N(x_{sp2})$,  the feasible set of the MPC control law designed to track $x_{sp2}$. As it can be seen, the designed control MPC control law can not be applied since $x(0) \not \in \setX_N(x_{sp2})$. Therefore, feasibility of the MPC control law is lost due to this change of setpoint.
		
		This feasibility loss is due to the terminal constraint $x(N)=x_{sp2}$, since the system cannot be steered from $x(0)$ to $x_{sp2}$ in only three steps, taking into account the constraints on the system. However, we have seen that the MPC control law would be feasible for the initial state if the the steady state $x_{sp1}=(4.9,0.245)$ (square) would be considered as a potential setpoint. Then the feasibility could be maintained if, instead of forcing the controller to reach the setpoint $x_{sp2}$ in $N$ steps, this constraint would be relaxed by forcing to reach some reachable setpoint $x_a$. This can be done by using the relaxed terminal constraint $x(N)=x_a$, where $x_a$ is a reachable setpoint that has to be added as new decision variable. This new reachable setpoint $x_a$ is called as artificial setpoint. Notice that the modified control law would be feasible by taking $x_a(0)=x_{sp1}$.
		
		\begin{figure}[h]
			\centering
			\includegraphics[width=0.99\textwidth]{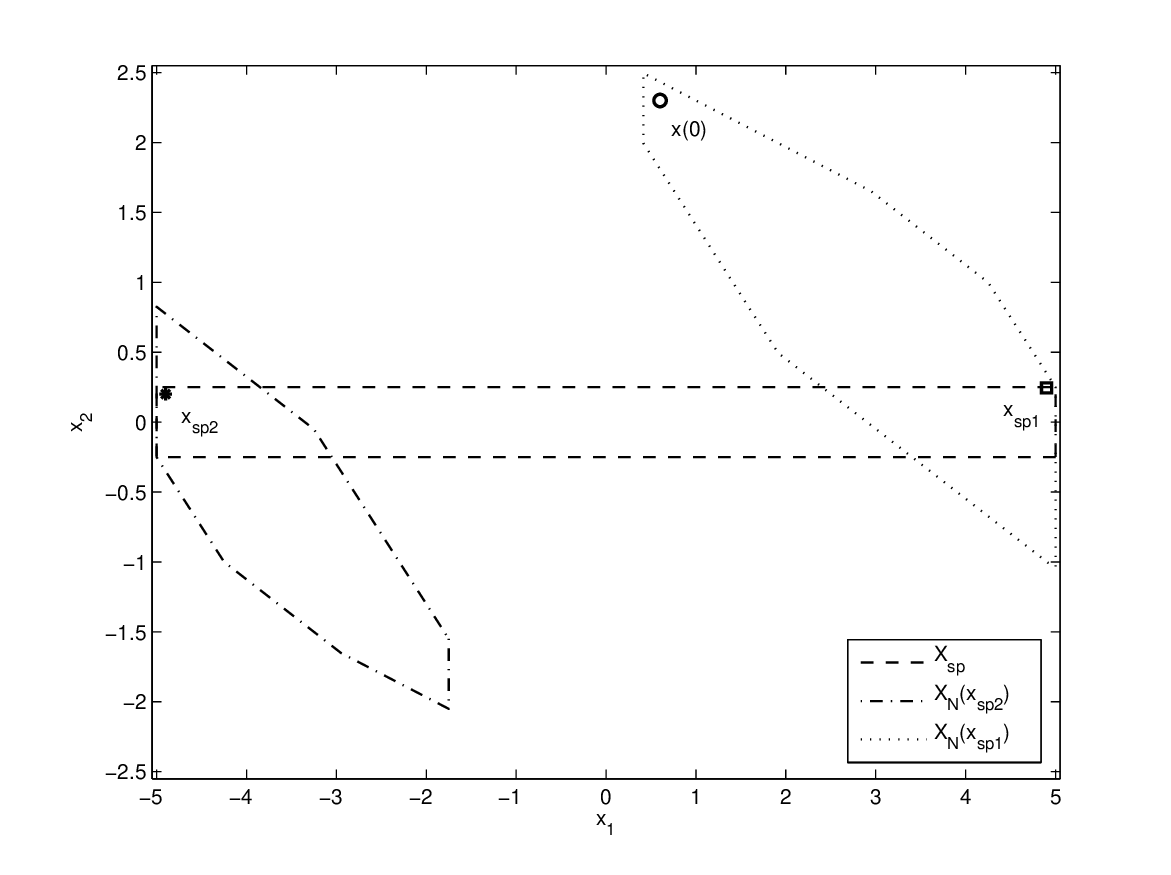}\\
			\caption{Feasibility: the initial condition $x(0)$ is clearly infeasible for an MPC for tracking $x_{sp2}$ with terminal constraint $x(N)=x_{sp2}$. However, it would be feasible for an MPC for tracking $x_a$ with terminal constraint $x(N)=x_a$.}\label{chaMPCT:di_ec_factibilidad}
		\end{figure}
		
		Using the relaxed terminal constraint, the feasibility of the predictive controller is ensured, but this does not ensure that the closed-loop system converges to the setpoint. This would be achieved if the optimal artificial setpoint $x_a$ would converge to the setpoint $x_{sp2}$ throughout the evolution of the controlled system. In order to ensure this convergence,  a term that penalizes the distance  between the artificial setpoint and the real setpoint is added  in the cost function. This term has the form of $V_O(x_a,x_{sp2})$ and it is the so-called offset cost function.
		
		This idea is illustrated in Figure \ref{chaMPCT:di_ec_factibilidad_2}. At $k=0$ the proposed predictive control scheme calculates the optimal predicted sequence of control actions together with the optimal artificial setpoint $x_a(0)$. The control input is applied and the system evolves to $x(1)=(2.65,1.55)$. At $k=1$, a new sequence of predicted control inputs and a new artificial setpoint $x_a(1)=(4.5,0.245)$ is obtained. As there is a term that penalizes the distance to the setpoint $x_{sp2}$, this new artificial setpoint $x_a(1)$ is different to $x_a(0)$ and closer to the setpoint. As it will be proved in the next section,  by doing this throughout the evolution of the plant, the artificial setpoint will converge to the setpoint and then the system state will converge to the desired setpoint.

		\begin{figure}[!h]
			\centering
			\includegraphics[width=0.99\linewidth]{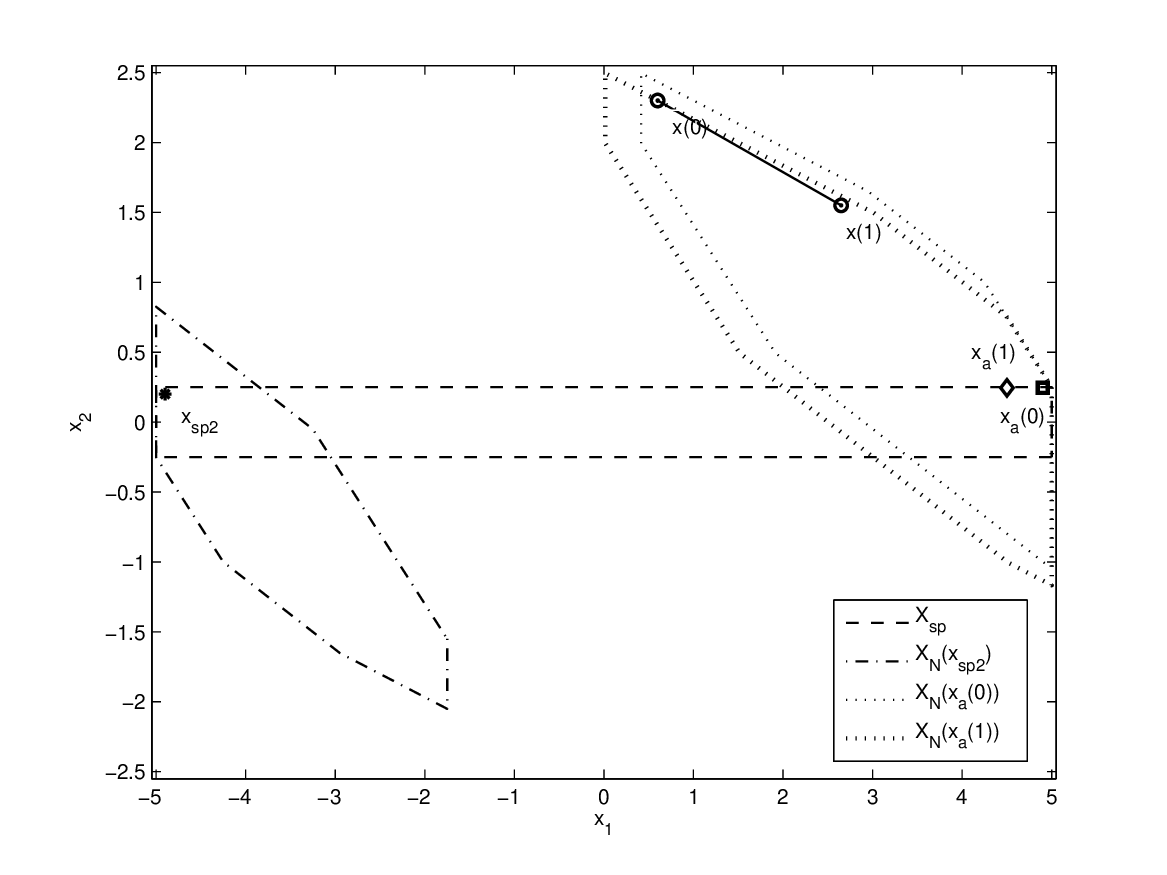}
			\caption{Feasibility: the offset cost function forces $x_a$ to move toward $x_{sp2}$ in order to ensure convergence. The MPC optimization problem remain feasible from the new initial condition $x(0)$.}\label{chaMPCT:di_ec_factibilidad_2}
		\end{figure}	
		
		%
		%	
		%Therefore, the predictive controller is not feasible for the setpoint but there exists an equilibrium point that could be tracked.	 This fact motivates the idea of relaxing the MPC terminal constraint by forcing the terminal state (i.e. the one predicted at the end of the prediction horizon) to be equal to any reachable equilibrium point $x_a$, that is $x_a \in \setX_s$. This can be done, by taking $(x_a,u_a)$ as extra decision variables of the MPC optimization problem and imposing the terminal constraint $x(N)=x_a$, and forcing that  $(x_a,u_a)$ is a reachable equilibrium point, i.e. $x_a= A x_a + B u_a$ and $(x_a,u_a) \in \lambda \setZ$.
		%******************
		%
		%designing a cost function that penalizes the distance of the predicted trajectory to this equilibrium point. We usually call $(x_a,u_a)$ as the artificial reference. Moreover, in order to force the convergence of $x_a$ to $x_{sp}$, we add to the MPC cost function, an extra cost, the offset cost function. This function penalizes the distance of $x_a$ to $x_{sp}$, and basically forces $x_a$ to move toward $x_{sp}$. This fact is depicted in Figure \ref{chaMPCT:di_ec_factibilidad_2}. At the next step, the system moves to a new point $x(1)=(2.65,1.55)$, which becomes the new initial condition. Thanks to the relaxed terminal constraint, the MPC problem remains feasible, and due to the offset cost function, the new artificial reference $x_a(1)=(4.5,0.245)$ is closer to $x_{sp}$.
		\QED
	\end{exmp}
	
	Recapping, the proposed predictive control scheme is characterized by the following features:
	
	%\begin{enumerate}
	%\item[(i)] an artificial reachable setpoint is added as decision
	%variable;
	%
	%\item[(ii) ] a stage cost that penalizes the deviation of the predicted
	%trajectory from the artificial steady conditions is considered;
	%
	%\item[(iii) ] an extra cost, the \emph{offset cost} function that penalizes the
	%deviation between the artificial setpoint and the real setpoint is added to the cost function;
	%
	%\item[(iv) ] A relaxed terminal constraint that depends on the artificial setpoint instead of on the real setpoint is considered.
	%\end{enumerate}

	\blista
	\item an artificial reachable setpoint is added as decision
	variable;
	
	\item a stage cost that penalizes the deviation of the predicted
	trajectory from the artificial steady conditions is considered;
	
	\item an extra cost, the \emph{offset cost} function, that penalizes the
	deviation between the artificial setpoint and the real setpoint is added to the cost function;
	
	\item a relaxed terminal constraint that depends on the artificial setpoint instead of on the real setpoint, is considered.
	\elista
	
	In the following sections, the MPC for tracking is precisely introduced. Besides the design procedure to ensure asymptotic stability to the setpoint is presented in both cases: based on an  terminal equality constraint and based on a terminal inequality constraint.

	\subsection{MPC for tracking with terminal equality constraint}
	
	As usual in predictive controllers, the MPC for tracking (MPCT)  is based on the solution at each sampling time of an optimization problem based on the current state and the current setpoint $(x,y_{sp})$, which are parameters of the optimization problem. The solution of this optimization will be applied using a receding horizon policy. The decision variables of this predictive controller are the sequence of control inputs $\vu$ and the steady state and input of the artificial setpoint $(x_a,u_a)$.
	
	The MPCT cost function depends on the parameters $(x,y_{sp})$ and on the decision variables $(\vu,x_a,u_a)$. This cost is composed of two terms: the first one, a dynamic term, is a quadratic cost of the expected tracking error with respect to the artificial steady state and input; the second one, a stationary term, is the offset cost function, that penalizes the deviation of the artificial setpoint $y_a$ to the setpoint $y_{sp}$ This cost function is calculated as follows:
	\beqna\label{chaMPCT:costfuc_eqconst}
	V_N(x,y_{sp};\vu,x_a,u_a)&=&\!\!\sum\limits_{j=0}^{N-1}\|x(j)\!-\!x_a\|^2_Q\!+\!\|u(j)\!-\!u_a\|^2_R+V_O(y_a,y_{sp})
	\eeqna
	where $x(j)$ denotes the prediction of the state $j$-samples ahead,
	the pair $(x_a, u_a)$ represents the artificial steady state
	and input, and $y_a=Cx_a+Du_a$ the artificial output or artificial setpoint; $y_{sp}$ is the desired setpoint.
	
	The function $V_O(y_a,y_{sp})$  is the so-called offset cost function and measures the distance between the artificial setpoint $y_a$ and the real setpoint $y_{sp}$. It is  assumed that this function  satisfies the following condition
	\begin{assumpt}\label{chaMPCT:def_optimo}
		Let the offset cost function $V_O:\mathbb{R}^{p} \rightarrow
		\mathbb{R}$ be a convex, positive definite and subdifferentiable
		function\footnote{Notice that a subdifferentiable function \cite{BoydLIB06} is a function that admits subgradients. Given a function $f$, $g$ is a subgradient of $f$ at $x$ if $$f(y)\geq f(x)+g'(y-x) \quad \forall y$$ Notice also that, the term  subdifferential defines the set of all subgradients of $f$ at $x$ and is noted as $\partial f(x)$. This set is a nonempty closed convex set.
		}, with $V_O(0,0)=0$, such that the minimizer of
		$$\min\limits _{y\in \mathcal Y_{sp}} V_O(y,y_{sp})$$ is unique.
	\end{assumpt}
	%

	% From the previous definition, it is clear that the constraints on the choice of $V_O(\cdot,\cdot)$ are not that restrictive. As it will be better explained in the following reports of this book, this function may be chosen as a norm of the distance to the desired steady state \cite{LimonAUT08,FerramoscaAUT09}, i.e. $\|y-y_{sp}\|_\infty$, as a distance to a target set \cite{FerramoscaJPC10}, or as a generic economic cost function as well, providing significant properties to the closed-loop system.\\
	% \hfill

	In the case of terminal equality constraint, the controller is derived from the solution of
	the following optimization problem:
	\begin{subequations}\label{chaMPCT:optprob_eqconst}
		\beqna
		V_N^{0}(x,y_{sp})&=&\min\limits_{\vu,x_a,u_a} V_N(x,y_{sp};\vu,x_a,u_a)\\
		&s.t.& x(0)=x, \label{chaMPCT:MPCT1_Eq1}\\
		&& x(j+1)=A x(j)+B u(j), \label{chaMPCT:MPCT1_Eq2}\\
		&& (x(j),u(j)) \in \setZ, \quad  \!j\!=\!0,\cdots, N\!-\!1 \label{chaMPCT:MPCT1_Eq3} \\
		&& x_a=Ax_a+Bu_a, \label{chaMPCT:MPCT1_Eq4}\\
		&& y_a=Cx_a+Du_a, \label{chaMPCT:MPCT1_Eq5}\\
		&& (x_a,u_a) \in \lambda \setZ \label{chaMPCT:MPCT1_Eq6}\\
		&& x(N)=x_a \label{chaMPCT:MPCT1_Eq7}
		\eeqna
	\end{subequations}

	Constraints (\ref{chaMPCT:MPCT1_Eq1})-(\ref{chaMPCT:MPCT1_Eq3}) force to the predicted trajectory to be consistent with the dynamic model equations while  the constraints are fulfilled.  Constraints (\ref{chaMPCT:MPCT1_Eq4})-(\ref{chaMPCT:MPCT1_Eq5}) make the artificial state and input $(x_a, u_a)$ to be a steady state and input for the prediction model. Constraint (\ref{chaMPCT:MPCT1_Eq6}) ensures that the artificial equilibrium point $(x_a,u_a)$ is admissible (and the constraints are not active). The last equation ( \ref{chaMPCT:MPCT1_Eq7}) is the relaxed terminal equality constraint that forces to the terminal state, that is the predicted state at the end of the prediction horizon, to be equal to the artificial state.

	The optimization problem can be posed as a quadratic programming problem and can be solved using specialized and extraordinarily efficient algorithms \cite{Nocedal99numerical}. In section \ref{chaMPCT:implementation_QP} it will be detailed how to formulate this problem as a quadratic programming problem.  The optimal solution of this optimization problem is denoted as $(\vu^0,x_a^0, u_a^0)$ and  depends on the parameters of the optimization problem $(x,y_{sp})$.  Considering the receding horizon policy, the control law is given by
	$$\kappa_N(x,y_{sp})=u^0(0;x,y_{sp})$$

	The feasibility region of this optimization problem is the set of parameters $(x,y_{sp})$ where the optimization problem has a solution, that is, it is feasible. Since the constraints (\ref{chaMPCT:MPCT1_Eq1})-(\ref{chaMPCT:MPCT1_Eq7}) do not depend on the setpoint $y_{sp}$, the feasibility of this optimization problem does not depend on the setpoint $y_{sp}$, but only on the current state  $x$. This means that the feasibility of the optimization problem cannot be lost due to a change of the setpoint.
	
	Furthermore, as the constraints of the optimization problem are linear, then the feasible region results to be a polyhedron, that is, the intersection of a set of linear inequalities.
	
	For a given reachable setpoint $x_a$,  $\setX_N(x_a)$ denotes the set of states that can be steered to $x_a$ in $N$ steps satisfying the constraints on the state and input throughout its evolution. This is the domain of attraction of a regulation MPC designed to regulate the system to $x_a$. In MPCT, the artificial setpoint is a decision variable, which means that its feasible region  $\setX_N$ is the  set of states $x$ that
	can steered to \textit{any} reachable steady state in $N$ steps, satisfying the constraints. This can be read as follows
	
	$$
	\setX_N= \bigcup\limits_{x_a \in \setX_{sp}} \setX_N(x_a)
	$$

	\begin{exmp}\label{chaMPCT:ejemplo_2}
		Consider the system described in Example \ref{chaMPCT:ejemplo_1}, where an MPC with $N=3$ must be designed to track the system from the initial state  $x(0)=(0.6,2.3)$ to the setpoint  $y_{sp}=x_{sp}=(-4.9,0.2)$. The MPC for tracking with weighting matrices $Q=I_2$ and $R=I_2$ and $V_O(y_a,y_{sp})=10\|y_a-y_{sp}\|_\infty$ as offset cost function has been designed.

		\begin{figure}[!h]
			\centering
			\includegraphics[width=0.99\textwidth]{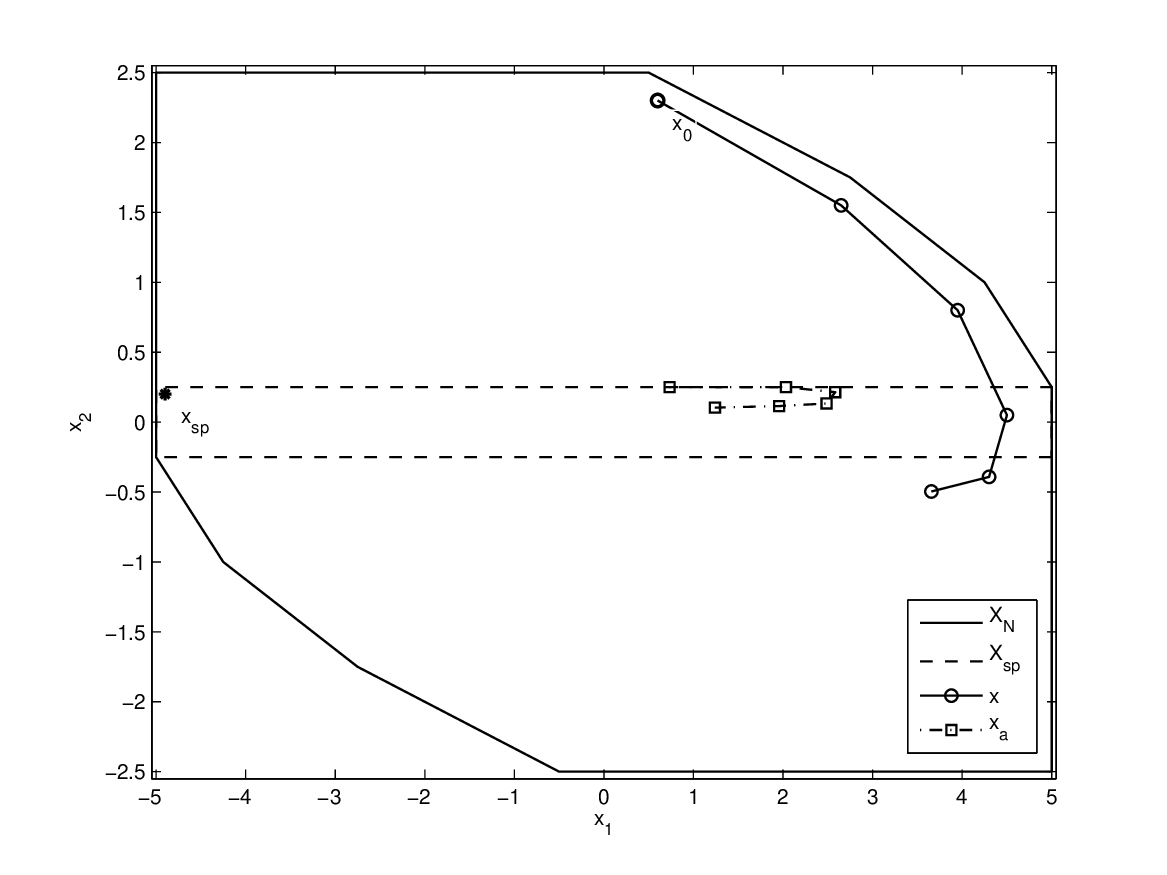}
			\caption{State space evolution of the first 5 steps of simulation: the artificial reference $x_a$ (dahsed-dotted line) never leaves $\setX_s$ (dashed line) and moves toward $x_{sp}$ (star).}\label{chaMPCT:di_ec_state_space_f2}
		\end{figure}
		
		In Figure \ref{chaMPCT:di_ec_state_space_f2} the state space evolution of the first 5 steps of the simulation is represented. The desired setpoint is depicted as a star, the evolution of the artificial reference $x_a$ in dash-dotted line, and the evolution of the closed-loop system in solid line. The feasible set of the MPCT, $\setX_N$, is plotted in solid line. Notice that the $x(0)$ is now inside $\setX_N$, and hence the optimization problem is feasible, even for $N=3$. Notice also how the evolution of $x_a$ changes along the time to ensure the feasibility of the optimization problem at each sampling time and its evolution converges toward $x_{sp}$. Notice also that $x_a$ never leaves sets $\setX_{sp}$ (in dashed line), which means that $x_a$ is always an admissible steady state.
		
		Therefore, the trajectory of the closed-loop system tracks the artificial reference, which evolves to maintain the recursive feasibility but converges to the setpoint thanks to the offset cost function. This is particularly clear in Figures \ref{chaMPCT:di_ec_state_space} and \ref{chaMPCT:di_ec_y_evol}, which represent the state space evolution and the time evolution of the complete simulation. The artificial reference $x_a$ never leaves set $\setX_s$ and eventually converges to $x_{sp}$. The closed-loop system tracks the artificial reference and is driven to the desired setpoint $x_{sp}$.
		
		\begin{figure}[!h]
			\centering
			\includegraphics[width=0.99\textwidth]{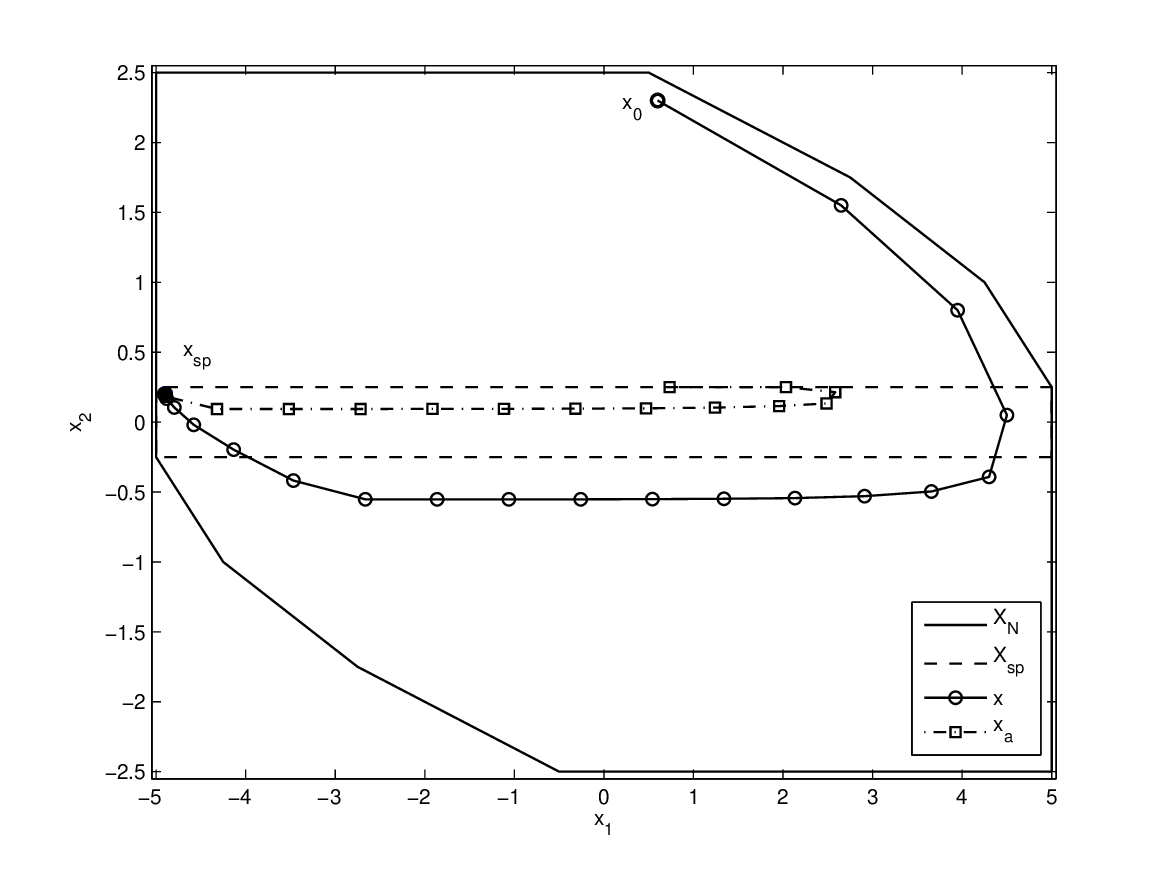}
			\caption{State space evolution of the complete simulation: both artificial reference (dash-dotted line) and closed-loop system (solid line) converge to the desired setpoint (star).}\label{chaMPCT:di_ec_state_space}
		\end{figure}
		
		\begin{figure}[!h]
			\centering
			\includegraphics[width=0.99\textwidth]{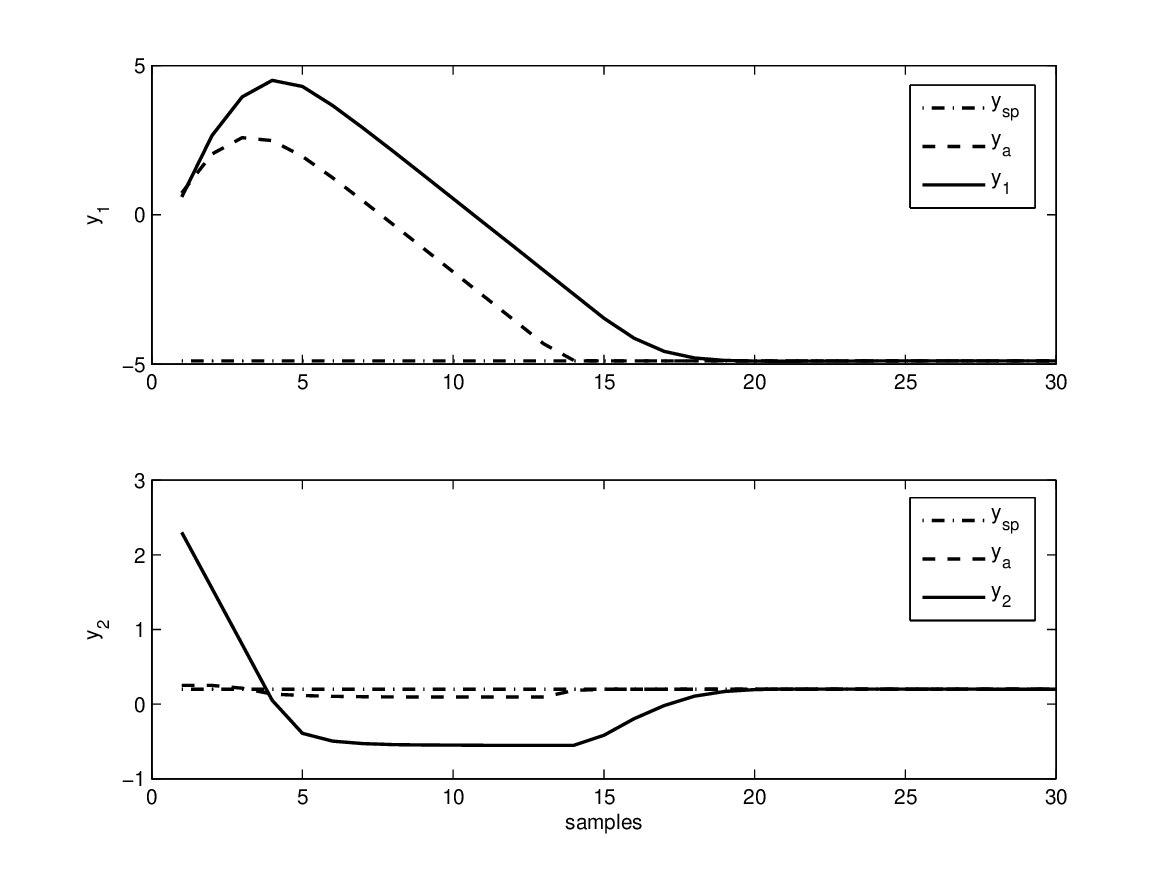}
			\caption{Time evolution of the complete simulation: the artificial reference (dashed line) converges to the desired setpoint (dash-dotted line). The closed-loop system (solid line) follows the artificial reference and eventually converges to $y_{sp}$.}\label{chaMPCT:di_ec_y_evol}
		\end{figure}
		
	\end{exmp}
	
	% % % % % % % % % % % % % % % % % % % % % % % % % % % % % % % % % % % % % % % % % % % % % % % % % % % %
	\subsubsection{Stabilizing design}
	
	The controller proposed in this section can be designed to ensure asymptotic stability of any reachable setpoint in the Lyapunov sense. The tuning parameters of the proposed controller are the weighting matrices of the stage cost $Q$ and $R$, the prediction horizon $N$ and the offset cost function $V_O(\cdot)$.
	The stability will be stated under the following assumptions on the system and on the controller parameters:
	\begin{assumpt}\label{chaMPCT:assumptionQR_eqconst}
		\begin{enumerate}
			\item[]
			%\item The system is square, i.e. $p=m$, and condition (\ref{EqMSS}) holds.
			
			\item The pair $(A,B)$ is controllable, and  $n_c\geq 1$ is the minimum integer such that the matrix
			$$
			\bmat{cccc} A^{n_c-1}B, & A^{n_c-2}B, & \cdots & B \emat
			$$
			is full rank.
			
			\item Matrix  $R \in \mathbb{R}^{m \times m}$ is a positive definite matrix and $Q \in \mathbb{R}^{n \times n}$ is a positive semi-definite matrix such that the pair $(Q^{1/2},A)$ is observable.
			%\footnote{A pair $(A,C)$ is observable if the observability matrix $[C^T,(CA)^T, (CA^2)^T, \cdots, (CA^{n-1})^T]^T$ is full rank.}.
			
		\end{enumerate}
	\end{assumpt}
	%
	
	%The set of admissible steady outputs is given by:
	%$$\setY_s=\{y=Cx_a+Du_a \mid  (x_a,u_a) \in \setZ_s\}$$
	%
	%This set is equal to  the set of all admissible outputs for system
	%\eqref{sistema_o} subject to \eqref{restric_oOPT}, that is,
	%$\setY_s$.
	
	Taking into account the proposed conditions on the controller
	parameters, in the following theorem asymptotic stability and
	constraints satisfaction of the controlled system are proved.
	\\
	
	% Para garantizar que existe la cota superior del coste \'{o}ptimo es necesario demostrar que el coste
	% \'{o}ptimo es locally bounded. Esto est\'{a} demostrado solo si los ingredientes son continuos y U es compacto.
	% Me surge la duda si la funci\'{o}n de coste \'{o}ptima, al ser convexa,  es LB
	% Entonces por sencillez del teorema  se asume que el conjunto Z est\'{a} acotado en la definici\'{o}n
	
	\begin{theorem}[Asymptotic Stability]\label{chaMPCT:TeorEst_eqconst}
		
		Consider that Assumptions \ref{assumption1_OPT}, \ref{chaMPCT:def_optimo} and
		\ref{chaMPCT:assumptionQR_eqconst} hold and the prediction horizon is such that $N\geq n_c$. Then for a given setpoint $y_{sp}$ and for any feasible initial state $x_0 \in \setX_N$,
		the system controlled by the MPC controller $\kappa_N(x,y_{sp})$
		is stable, fulfills the constraints throughout the time and, besides
		
		\begin{enumerate}
			\item[(i)] If the setpoint is reachable, i.e. $y_{sp}\in \setY_{sp}$,
			then the closed-loop system asymptotically converges to the steady
			state, input and output $(x_{sp},u_{sp},y_{sp})$.
			\item [(ii) ] If the setpoint is not reachable, i.e. $y_{sp} \not\in \setY_{sp}$, then the closed-loop system asymptotically converges to a reachable steady
			state, input and output  $(x_t,u_t,y_t)$ where the offset cost function is minimal, that is,
			\[
			y_t=\arg \min_{y\in \setY_{sp}}V_O(y,y_{sp})
			\]
		\end{enumerate}
		
	\end{theorem}
	
	\begin{pf}
		To derive the proof, only the second statement must be proved, since if the setpoint is feasible, then
		\[
		y_{sp}=\arg \min_{y\in \setY_{sp}}V_O(y,y_{sp})
		\]

		Consider that $x \in \setX_N$ at time $k$, then the optimal cost
		function is given by $V_N^0(x,y_{sp})=V_N(x,y_{sp};\vu^0(x),x_a^0(x),u_a^0(x))$, where
		$(\vu^0(x),x_a^0(x),u_a^0(x))$ defines the optimal solution to Problem \eqref{chaMPCT:optprob_eqconst}
		and $\vu^0(x)=\{u^0(0;x),u^0(1;x),...,u^0(N-1;x)\}$\footnote{The dependence from $y_{sp}$ will be omitted for the sake of clarity.}. The resultant
		optimal sequence of predicted states associated to $\vu^0(x)$  is given by
		$\vx^0(x)=\{x^0(0;x),x^0(1;x),...,x^0(N-1;x),x^0(N;x)\}$, where
		$x^0(j;x)=\phi(j;x,\vu^0(x))$ and $x^0(N;x)=x_a^0(x)$.
		
		The proof will be carried out in two stages: first the recursive feasilibity will be shown, and next, it will be demonstrated that there exists a Lyapunov function based on the optimal cost function.
		
		As standard in MPC \cite[Chapter 2]{RawlingsLIB09},
		define the successor state at time $k+1$, $x^+=Ax+Bu^0(0;x)$ and
		define also the following sequences:
		\begin{eqnarray*}
			\tilde \vu(x)&\!\equal\!&\{u^0(1;x),\cdots, u^0(N\!-\!1;x), u_a^0(x)\}\\
			\tilde x_a(x)&\!\equal\!& x_a^0(x)\\
			\tilde u_a(x)&\!\equal\!& u_a^0(x)
		\end{eqnarray*}
		Since $x^+=x(1;x)$, then the sequence of predicted states starting from $x^+$, when  the feasible solution $(\tilde \vu(x),\tx_a(x),\tu_a(x))$ is applied, is given by
		$$\tilde \vx=\{x^0(1;x),x^0(2;x),...,x^0(N;x),x^0(N+1;x)\}$$ where
		$x^0(N+1;x)=Ax^0(N;x)+Bu_a^0(x)=x_a^0(x)$. As the optimal solution is feasible and fulfills the constraints, then the trajectories $\tilde \vx(x) $ and $\tilde \vu(x)$ will be also feasible.  Therefore, the optimization problem is recursively feasible.
		
		In order to derive the Lyapunov function, we recur to the following result.
		Comparing the optimal cost $V_N^0(x,y_{sp})$, with the cost
		given by $(\tilde \vu(x),\tx_a(x),\tu_a(x))$, at time $k+1$, $\tilde V_N(x^+,y_{sp};\tilde
		\vu(x),\tx_a(x),\tu_a(x))$, we have
		\beqnan\tilde V_N\!(x^+\!\!,\!y_{sp}; \tilde
		\vu(x),\tx_a(x),\tu_a(x)\!)\!-\!V_N^{0}\!(x,y_{sp}\!)\!\!\!\!\!\!&=&\!\!\!\!\!\!-\|x\! -\!
		x_a^0(x)\|^2_Q\!-\!\|u^0(0;x)\! -\! u_a^0(x)\|^2_R\!\\
		\!\!\!\!\!\!&&\!\!\!\!\!\! -\!\! \!\sum\limits_{j=1}^{N-1}\!\!\!\|x^0\!(j;x)\! -\!
		x_a^0\!(x)\!\|^2_Q\!-\!\|u^0\!(j;x)\! -\! u_a^0\!(x)\!\|^2_R\\
		\!\!\!\!\!\!&&\!\!\!\!\!\! -V_O(y_s,y_t)\! +\!V_O(y_s,y_t)\\
		\!\!\!\!\!\!&&\!\!\!\!\!\!+\!\!\!\sum\limits_{j=1}^{N-1}\!\!\!\|x^0\!(j;x)\! -\!
		x_a^0\!(x)\!\|^2_Q\!-\!\|u^0\!(j;x)\! -\! u_a^0\!(x)\!\|^2_R \\
		\!\!\!\!\!\!&&\!\!\!\!\!\!+\|x^0\!(N;x)\! -\! x_a^0(x)\|^2_Q\!+\!\|u_a^0(x)\! -\!
		u_a^0(x)\|^2_R\! \\
		\!\!\!\!\!\!&=&\!\!\!\!\!\! -\|x\! -\!x_a^0(x)\|^2_Q\!-\!\|u^0(0;x)\! -\!
		u_a^0(x)\|^2_R \eeqnan
		By optimality, we have that $V_N^{0}(x^+,y_{sp})\leq \tilde V_N(x^+,y_{sp}; \tilde
		\vu,\tx_a(x),\tu_a(x))$ and then:
		\beqnan V_N^{0}(x^+,y_{sp})-V_N^{0}(x,y_{sp})&\leq& -\|x -
		x_a^0(x)\|^2_Q-\|u^0(0;x) - u_a^0(x)\|^2_R  \eeqnan
		Taking into account that the cost function is positive definite, the previous inequality implies that there exists a $\mathcal{K}$-function $\alpha$ such that:
		\begin{equation}
			V_N^{0}(x^+,y_{sp})-V_N^{0}(x,y_{sp})\leq-\alpha(\| x-x_a^0(x)\|) \label{chaMPCT:EQ_costo_decre}
		\end{equation}
		Now we are ready to prove that the function
		$$J(x)=V^{0}_N(x,y_{sp})-V_O(y_t,y_{sp})$$
		is a Lyapunov function for the closed-loop system in the domain of attraction $\setX_N$. That is, that there exist three $\mathcal K_\infty$-function,$\alpha_1$, $\alpha_2$ and $\alpha_3$, such that for all $x \in \setX_N$ the following conditions hold
		\begin{eqnarray}
			J(x) &\geq & \alpha_1(\|x-x_{t}\|)\\
			J(x) &\leq & \alpha_2(\|x-x_{t}\|)\\
			J(Ax+B \kappa(x,y_{sp})) - J(x) &\leq & -\alpha_3(\|x-x_{t}\|).
		\end{eqnarray}

		This function is well defined in $\setX_{N}$. Define also $e(x)=x-x_a^0$. Notice that, since $Q$ and $R$ are positive definite, $J(x)\geq\alpha(\|e(x)\|)$, for all $x \in \setX_{N}$; due to \eqref{chaMPCT:EQ_costo_decre}, we have that $J(x^+)-J(x)\leq -\alpha(\|e(x)\|)$, for all $x \in \setX_{N}$.\\
		From Lemma \ref{chaMPCT:Lem_Func_DefPos}, it follows that $$\alpha(\|e(x)\|)\geq \alpha(\alpha_e(\|x-x_{sp}\|))=\alpha_J(\|x-x_{sp}\|)$$ where $\alpha_e$ and $\alpha_J$ are $\mathcal{K}$-functions. Then, we can conclude that:
		\blista
		
		\item $J(x)\geq\alpha_1(\|x-x_{sp}\|)$, for all $x \in \setX_{N}$.
		
		From the definition of the function $J(x)$ we have that
		$$ J(x)=V^{0}_N(x,y_{sp})-V_O(y_t,y_{sp})\geq \|x - x^0_a\|_Q^2+ V_O(y^0_a,y_{sp}) -V_O(y_t,y_{sp})$$
		
		From the optimality of $y_t$, it is derived that $V_O(y^0_a,y_{sp}) -V_O(y_t,y_{sp}) \geq 0$. Defining $\alpha_{11}(\|x - x^0_a\|) =\|x - x^0_a\|_Q^2$, we have that
		$$ J(x)=V^{0}_N(x,y_{sp})-V_O(y_t,y_{sp})\geq \alpha_{11}(\|x - x^0_a\|) $$
		From Lemma \ref{chaMPCT:Lem_Func_DefPos}, it follows that there exists a $\mathcal K_\infty$ function $\alpha_{12}$ such that  $\|x - x^0_a\| \geq \alpha_{12}(\|x - x_t\|)$, we have
		$$ J(x)=V^{0}_N(x,y_{sp})-V_O(y_t,y_{sp})\geq \alpha_1(\|x - x^0_a\|) $$
		for some  $\mathcal K_{\infty}$ function.

		%From the definition of the function $J(x)$ we have that
		%$$ J(x)=V^{0}_N(x,y_{sp})-V_O(y_t,y_{sp})\geq \|x - x^0_a\|_Q^2+ V_O(y^0_a,y_{sp}) -V_O(y_t,y_{sp})$$
		%From the optimality of $y_t$ and the strict convexity of the offset cost function, there exists a $\mathcal K_\infty$ function such that
		%$$ V_O(y^0_a,y_{sp}) -V_O(y_t,y_{sp}) \geq \alpha_{12}(\|y^0_a - y_t\|)$$
		%As it was discussed in section \ref{chaMPCT:ReachSP}, since the system is square and condition (\ref{EqMSS}) holds, then for a given setpoint $y_{s}$  there exists a unique equilibrium point $(x_{s},u_{s})$. Therefore, there exists constant $c>0$ such that $ \|y^0_a - y_t\| \geq c \|x^0_a - x_t\|$.
		%Then there exists a $\mathcal K_\infty$ function such that
		%$$ V_O(y^0_a,y_{sp}) -V_O(y_t,y_{sp}) \geq \alpha_{13}(\|x^0_a - x_t\|)$$
		%
		%Defining $\alpha_{11}(\|x - x^0_a\|) =\|x - x^0_a\|_Q^2$, we have that
		%
		%$$ J(x)\geq \alpha_{11}(\|x - x^0_a\|) + \alpha_{13}(\|x^0_a - x_t\|)$$
		%and then there exists a $\mathcal K_\infty$ function such that
		%$$ J(x) \geq \alpha_{1}(\|x - x^0_a\| + \|x^0_a - x_t\|) \geq \alpha_{1}(\|x  - x_t\|)$$

		\item $J(x)\leq\alpha_2(\|x-x_t\|)$, for all $x \in \setX_{N}$.
		
		Since the stage cost function is quadratic and the model is linear, the optimal cost function $J(x)=V_N^0(x,y_{sp}) -V_O(y_t,y_{sp}) $ is a locally bounded continuous function and $J(x_t)=0$, then there exists a $\mathcal K_\infty$ function $\alpha_2(\cdot)$ such that $J(x)\leq\alpha_2(\|x-x_t\|)$, for all $x \in \setX_{N}$ (see Propositions 1 and 2 of the postface to the book \cite{RawlingsLIB09}).
		
		% aqu\'{\i} se debe citar mejor el paper de estabilidad con Picasso

		\item $J(Ax+B \kappa(x,y_{sp})) - J(x) \leq  \alpha_3(\|x-x_{t}\|)$ for all $x \in \setX_N$.
		
		From equation (\ref{chaMPCT:EQ_costo_decre}), we have that
		
		$$J(Ax+B \kappa(x,y_{sp})) - J(x)\leq-\alpha(\| x-x_a^0(x)\|) $$
		
		On the other hand, we have previously proved that there exists a $\mathcal K_\infty$ function $\alpha_{12}$ such that  $\|x - x^0_a\| \geq \alpha_{12}(\|x - x_t\|)$. Then there exists  a $\mathcal K_\infty$ function $\alpha_3(\cdot)$ such that $\alpha(\|x - x^0_a\|) \geq \alpha_3(\|x - x_t\|)$. Thus, it is easily derived that
		$$J(Ax+B \kappa(x,y_{sp})) - J(x)\leq-\alpha_3(\| x-x_t\|) $$
		
		\elista
		
		We have proved that $J(x)$ is a Lyapunov function for the controlled system, and then  $(x_t,u_t)$ is an asymptotically stable equilibrium point and its domain of attraction is $\setX_N$.

		%Hence $J(x)$ is a Lyapunov function and $x_{sp}$ is an asymptotically stable equilibrium point for the closed-loop system, that is, there exists a $\mathcal{KL}$-function $\vartheta$ such that $$\|x(k)-x_{sp}\|\leq \vartheta(\|x(0)-x_{sp}\|,k)$$ for all $x(0) \in \setX_{N}$.

	\end{pf}
	
	From this theorem we have that the resulting controller steers the system from a feasible initial state  to any admissible setpoint satisfying the constraints on the input and state throughout its trajectory. Besides, as it is customary in stabilizing MPC, the optimal cost function is a Lyapunov function of the controlled system. Further properties of this controller will be shown later on.

	\subsection{MPCT with terminal inequality constraint}
	
	 MPC schemes can be stabilized by adding a terminal cost function and an inequality constraint on the terminal state. The terminal constraint forces the terminal state to be in a region that is a domain of attraction of the system stabilized by a local controller, called terminal control law. The terminal cost function is chosen to be a Lyapunov function of the system controlled by the terminal control law. These assumptions ensure the existence of a feasible solution based on optimal solution at last sampling time and ensure that the optimal cost function is a Lyapunov function for the controlled system.
	The main advantages of this methodology is that the resulting controller has a larger domain of attraction and a better closed-loop performance than the controller that uses an equality constraint.
	
	In this section it will be shown how the MPC for tracking can be designed using both, a terminal cost function and a terminal inequality constraint. As in the regulation case, a larger domain of attraction and a better closed-loop performance will be achieved.
	
	The underlying idea is similar to the one of the regulation case: a (linear) terminal control law must be designed to stabilize the system to any admissible equilibrium point $(x_a,u_a)$. Then a Lyapunov function for the controlled system, $V_f(x-x_a)=\|x-x_a\|_P^2$,  is used as the terminal cost function and a domain of attraction of the controlled system $(x,x_a,u_a) \in \Omega_t^a$ is used as terminal constraint. Notice that the terminal region depends on the admissible equilibrium point.
	
	%
	%With the aim of enlarging the domain of attraction, a formulation with a terminal inequality constraint can also be implemented. In this case, the relaxed terminal constraint can be formulated by means of a terminal invariant set, called the \emph{invariant set for tracking}.
	%In this formulation, we include in the cost function, a cost-to-go from $N$ to infinity, as well as a polytopic terminal inequality constraint in the optimization problem. The MPCT cost function is given by:
	
	Then, the cost function of the optimization problem of the MPC for tracking is  given by

	\beqna\label{chaMPCT:costfuc_ineqconst}
	V_N(x,y_{sp};\vu,x_a,u_a)\!\!&=&\!\!\sum\limits_{j=0}^{N-1}\!\!\|x(j)\!-\!x_a\|^2_Q\!+\!\|u(j)\!-\!u_a\|^2_R \nonumber \\ && +\|x(N)\!-\!x_a\|^2_P\!+V_O(y_a,y_{sp})
	\eeqna
	The controller is derived from the solution of
	the following optimization problem:
	\begin{subequations}\label{chaMPCT:optprob_ineqconst}
		\beqna
		V_N^{0}(x,y_{sp})&=&\min \limits_{\vu,x_a,u_a} V_N(x,y_{sp};\vu,x_a,u_a)\\
		&s.t.& x(0)=x, \\
		&& x(j+1)=A x(j)+B u(j), \\
		&& (x(j),u(j)) \in \setZ, \quad  \!j\!=\!0,\cdots, N\!-\!1 \\
		% && x_a=A x_a + B u_a\\
		&& y_a=Cx_a+Du_a,\\
		&& (x(N),x_a,u_a) \in \Omega_{t}^a
		\eeqna
	\end{subequations}
	Considering the receding horizon policy, the control law is given by
	$$\kappa_N(x,y_{sp})=u^0(0;x,y_{sp})$$

	As in the case of equality terminal constraint, the set of constraints of Problem \eqref{chaMPCT:optprob_ineqconst} does not depend on the setpoint $y_{sp}$, and then the feasible region does not depend on $y_{sp}$. Let define $\Omega_t$ as the projection of $\Omega_t^a$ onto $x$, then the feasible region is the set of states that can be admissible steered to the set $\Omega_t$ in $N$ steps. This set will be denoted as $\setX_N(\Omega_t)$.
	
	Next, a design procedure to provide closed loop stability is shown.

	\subsubsection{Stabilizing design}
	
	The design parameters of this controller are the weighting matrices $Q$ and $R$,  the prediction horizon $N$, the offset cost function $V_O(\cdot,\cdot)$, the weighting matrix of the quadratic terminal cost function $P$ and the extended terminal constraint set $\Omega_t^a$. These parameters must fulfill the following assumption

	\begin{assumpt}\label{chaMPCT:assumption_STAB_ineqcons}
		\begin{enumerate}
			\hfill
			\item Let  $R \in \mathbb{R}^{m \times
				m}$ be a positive definite matrix and $Q \in \mathbb{R}^{n \times
				n}$ a positive semi-definite matrix such that the pair $(Q^{1/2},A)$
			is observable.
			\item Let $K \in \mathbb{R}^{m \times n}$ be a stabilizing control gain such that $(A+BK)$ has the eigenvalues inside the unit circle.
			\item Let $P \in \mathbb{R}^{n \times n}$ be a positive definite matrix such that:
			\begin{eqnarray*}
				\!\!(A\!+\!BK)'\!P\!(A\!+\!BK)\!-\!P\!\!\leq \!\!-\!(Q\!+\!K'RK)
			\end{eqnarray*}
			\item Let $\Omega_{t}^a\subseteq\mathbb{R}^{n+m}$ be an admissible
			polyhedral invariant set for tracking for system \eqref{sistema_o}
			subject to \eqref{restric_oOPT}, for a given gain $K$. That is, for all
			$(x,x_a,u_a) \in \Omega_{t}^a$, the following conditions must hold
			\beqnan
			(x,K(x-x_a)+u_a) &\in& \setZ\\
			(x_a,u_a) &\in& \setZ_{sp}\\
			(A x + B (K (x-x_a)+u_a),x_a,u_a) &\in& \Omega_{t}^a
			\eeqnan
		\end{enumerate}
	\end{assumpt}

	Notice that the assumption on matrix $P$ is a Lyapunov condition and ensures that $V_f(x-x_a)= \|x-x_a\|^2_P$ is a Lyapunov-type function for the system controlled by the terminal control law $u=K (x-x_a) + u_a$, that is $x^+=A x + B (K (x-x_a) + u_a)$. The existence of this matrix is ensured  thanks to the stabilizing control gain $K$.
	
	On the other hand, the invariant set for tracking condition of $\Omega_t^a$ ensures that for all $(x, x_a, u_a) \in \Omega_t^a$, $(x_a,u_a)$ is an admissible equilibrium point, the control input $u=K(x-x_a) + u_a$ is admissible, i.e. $(x,u) \in \setZ$, and the successor state $x^+$  remains in $\Omega_t^a$ for the same equilibrium point $(x_a,u_a)$, i.e. $(x^+,x_a,u_a) \in \Omega_t^a$.
	Therefore, for all initial state $x(0)$ and admissible equilibrium point $(x_a,u_a)\in \setZ_{sp}$, such that $(x(0), x_a, u_a) \in \Omega_t^a$, the terminal control law $u=K(x-x_a) + u_a$ steers the controlled system to $(x_a,u_a)$ , the trajectory is such that $(x(k),u(k)) \in \setZ$ and $(x(k),x_a,u_a) \in \Omega_t^a$ for all $k\geq 0$ and $V_f(x(k)-x_a)= \|x(k)-x_a\|^2_P$ is strictly decreasing.
	
	It is interesting to notice that  $(x_a,x_a,u_a) \in \Omega_t^a$, which means that any admissible steady state $x_a$ is contained in $\Omega_t$.
	%It can be seen that $\Omega_{t}^a$ contains the set of
	%equilibrium points $ \setZ_s$.
	%
	%The set of admissible steady outputs consistent with the invariant
	%set for tracking $\Omega_{t}^a$ is given by:
	%$$\{y_a=Cx_a+Du_a \mid  x_a=Ax_a+Bu_a, \mbox{and } (x,x_a,u_a) \in \Omega_{t}^a \}$$
	%
	%This set is equal to  the set of all admissible outputs for system
	%\eqref{sistema_o} subject to \eqref{restric_oOPT}, that is,
	%$\setY_s$.
	
	Now, asymptotic stability of the proposed controller is stated.
	
	\begin{theorem}[Stability]\label{chaMPCT:TeorEst_ineqcons}
		
		Consider that Assumptions \ref{assumption1_OPT}, \ref{chaMPCT:def_optimo} and
		\ref{chaMPCT:assumption_STAB_ineqcons}  hold and consider a given setpoint $y_{sp}$. Let $\kappa_N(x,y_{sp})$ be the MPC control law resulting from the optimization problem \eqref{chaMPCT:optprob_ineqconst}. Then for any feasible initial state $x_0 \in \setX_N(\Omega_t)$, the closed-loop system $x^+=A x + B \kappa_N(x,y_{sp})$	is stable, fulfills the constraints throughout the time and, besides
		
		\begin{enumerate}
			\item[(i)] If the setpoint is reachable, i.e. $y_{sp}\in \setY_{sp}$,
			then the closed-loop system asymptotically converges to the steady
			state, input and output $(x_{sp},u_{sp},y_{sp})$.
			\item [(ii) ] If the setpoint is not reachable, i.e. $y_{sp} \not\in \setY_{sp}$, then the closed-loop system asymptotically converges to a reachable steady
			state, input and output  $(x_t,u_t,y_t)$ where the offset cost function is minimal, that is,
			\[
			y_t=\arg \min_{y\in \setY_{sp}}V_O(y,y_{sp})
			\]
		\end{enumerate}

		%Consider that Assumptions \ref{assumption1_OPT} and
		%\ref{chaMPCT:assumption_STAB_ineqcons} hold and consider a given setpoint point $y_t$. Then for any feasible initial state $x_0 \in \setX_N$,
		%the system controlled by the proposed MPC controller $\kappa_N(x,y_{sp})$
		%is stable, fulfils the constraints throughout the time and, besides
		%
		%\begin{itemize}
		%\item[(i)] If $y_{sp}\in \setY_{sp}$
		%then the closed-loop system asymptotically converges to a steady
		%state and input $(x_{sp},u_{sp})$ such that $y_{sp}=C x_{sp} + D u_{sp}$.
		%\item [(ii)] If $y_{sp} \not \in \setY_{sp}$, the closed-loop system asymptotically converges to a steady
		%state and input $(x_t,u_t)$ and $y_t=C x_t + D u_t$ where
		%\[
		%y_t=\arg \min_{y\in \setY_{sp}}V_O(y,y_{sp})
		%\]
		%\end{itemize}
		
	\end{theorem}

	\begin{pf}
		The proof to this theorem follows same arguments as the proof to Theorem \ref{chaMPCT:TeorEst_eqconst}.
		
		Consider that $x \in \setX_N(\Omega_t)$ at time $k$, then the optimal cost
		function is given by $V_N^0(x,y_{sp})=V_N(x,y_{sp};\vu^0(x),x_a^0(x),u_a^0(x))$, where
		$(\vu^0(x),x_a^0(x),u_a^0(x))$ defines the optimal solution to Problem \eqref{chaMPCT:optprob_ineqconst}
		and $\vu^0(x)=\{u^0(0;x),u^0(1;x),...,u^0(N-1;x)\}$\footnote{The dependence from $y_{sp}$ will be omitted for the sake of clarity.}. The resultant
		optimal state sequence associated to $\vu^0(x)$  is given by
		$\vx^0(x)=\{x^0(0;x),x^0(1;x),...,x^0(N-1;x),x^0(N;x)\}$, where
		$x^0(j;x)=\phi(j;x,\vu^0(x))$ and the last predicted state is such that $x^0(N;x) \in \Omega_t$.
		
		As standard in MPC \cite[Chapter 2]{RawlingsLIB09},
		define the successor state at time $k+1$, $x^+=Ax+Bu^0(0;x)$ and
		define also the following sequences:
		\begin{eqnarray*}
			\tilde \vu(x)&\!\equal\!&\{u^0(1;x),\cdots, u^0(N\!-\!1;x), K
			(x^0(N;x)-x_a^0(x))+u_a^0(x)\}\\
			\tilde x_a (x)&\!\equal\!& x_a^0(x)\\
			\tilde u_a (x)&\!\equal\!& u_a^0(x)
		\end{eqnarray*}
		Since $x^+=x^0(1;x)$, the state sequence associated to $(\tilde \vu(x), \tx_a (x),\tu_a (x))$ is
		$$\tilde \vx=\{x^0(1;x),x^0(2;x),...,x^0(N;x),x^0(N+1;x)\}$$
		where
		$x^0(N+1;x)=A x^0(N;x)+B (K(x^0(N;x)- x_a^0(x)) +  u_a^0(x))$.   Given that
		$$(x^0(N;x),x_a^0(x), u_a^0(x)) \in \Omega_t^a$$
		the control action
		$\tilde u(N-1;x)=K (x^0(N;x)-x_a^0(x))+u_a^0(x)$ is admissible, which means that $(x^0(N;x),K (x^0(N;x)-x_a^0(x))+u_a^0(x)) \in \setZ$. Besides the terminal state $x^0(N+1;x)$ is also
		feasible thanks to the properties of the invariant set for tracking, that is, $(x^0(N+1;x),\tilde x_a (x),\tilde u_a (x)) \in \Omega_t^a$.
		
		Then, we find that $(\tilde \vu(x), \tx_a (x),\tu_a (x))$ is a feasible solution to Problem \eqref{chaMPCT:optprob_ineqconst} at time $k+1$.

		Compare now the optimal cost, $V_N^0(x,y_{sp})$, with $\tilde V_N(x^+,y_{sp};\tilde \vu(x), \tx_a (x),\tu_a (x))$, that is the cost
		given by $(\tilde \vu(x), \tx_a (x),\tu_a (x))$. Taking into account the properties of the
		feasible nominal trajectories for $x^+$, the condition \emph{(4)} of
		Assumption \ref{chaMPCT:assumption_STAB_ineqcons} and using standard procedures in
		MPC \cite[Chapter 2]{RawlingsLIB09} it is possible to
		obtain:
		\beqnan\tilde V_N\!(x^+\!,y_{sp}; \tilde
		\vu,\tx_a (x),\tu_a (x)\!)\!-\!V_N^{0}\!(x,y_{sp})\!\!\!\!\!\!&=&\!\!\!\!\!\!-\|x\! -\!
		x_a^0(x)\|^2_Q\!-\!\|u^0(0;x)\! -\! u_a^0(x)\|^2_R\\
		\!\!\!\!\!\!&&\!\!\!\!\!\!-\!\|x^0(N;x)\!
		-\!x_a^0(x)\|^2_P\!-\!V_O(y_a,y_t)\\
		\!\!\!\!\!\!&&\!\!\!\!\!\!+\|x^0\!(N;x)\! -\! x_a^0\!(x)\!\|^2_Q\!+\!\|K\!(x^0\!(N;x)\! -\!
		x_a^0(x)\!)\!\|^2_R\!\\
		\!\!\!\!\!\!&&\!\!\!\!\!\!+\|x^0(N\!+\!1;x) \!-\! x_a^0(x)\|^2_P\! +\!V_O(y_a,y_t)\\
		\!\!\!\!\!\!&\leq &\!\!\!\!\!\! -\|x\! -\!x_a^0(x)\|^2_Q\!-\!\|u^0(0;x)\! -\!
		u_a^0(x)\|^2_R \eeqnan
		By optimality, we have that $V_N^{0}(x^+,y_{sp})\leq \tilde V_N(x^+,y_{sp}; \tilde
		\vu,\tx_a (x),\tu_a (x))$ and then:
		\beqnan V_N^{0}(x^+,y_{sp})-V_N^{0}(x,y_{sp})&\leq& -\|x -
		x_a^0(x)\|^2_Q-\|u^0(0;x) - u_a^0(x)\|^2_R  \eeqnan
		Taking into account that the cost function is positive definite, the previous inequality implies that there exists a $\mathcal{K}$-function $\alpha$ such that:
		\begin{equation}
			V_N^{0}(x^+,y_{sp})-V_N^{0}(x,y_{sp})\leq-\alpha(\| x-x_a^0(x)\|) \label{chaMPCT:INEQ_costo_decre}
		\end{equation}

		Define now the function $J(x)=V^{0}_N(x,y_{sp})-V_O(y_{t},y_{sp})$. Following the same arguments as in the proof to Theorem \ref{chaMPCT:TeorEst_eqconst}, it can be verified that $J(x)$ is a Lyapunov function and $(x_{t},u_t,y_t)$ is an asymptotically stable equilibrium point for the closed-loop system.
	\end{pf}

	In order to illustrate the proposed stabilizing design, the constrained double integrator example used in the Example \ref{chaMPCT:ejemplo_2} is controlled by using the inequality terminal constraint.
	
	\begin{exmp}\label{chaMPCT:ejemplo_3}
		Consider the tracking control problem presented in Example \ref{chaMPCT:ejemplo_2}. An  MPC for tracking has been designed using an inequality terminal constraint. In this case the terminal control law is a Linear Quadratic Regulator and the associated Lyapunov matrix is used as matrix $P$. This pair satisfies the conditions 2 and 3 of Assumption \ref{chaMPCT:assumption_STAB_ineqcons}. A polyhedral terminal region $\Omega_t^a$ has been calculated satisfying condition 4 (see the next section for the calculation of this set). Thus, the resulting controller asymptotically stabilize the system.
		
		In Figure \ref{chaMPCT:di_ic_state_space} the state space evolution of the closed-loop system is represented. The desired setpoint is depicted as a star, the evolution of the artificial reference $x_a$ in dash-dotted line, and the evolution of the closed-loop system in solid line. The projection of the invariant set for tracking, $\Omega_{t}$ is depicted in dash-dotted line. The time evolution is plotted in Figure \ref{chaMPCT:di_ic_y_evol}. Notice how the evolution of $x_a$ moves toward $x_{sp}$, and how the closed-loop system follows $x_a$. Notice also that the domain of attraction is larger than the domain in Example \ref{chaMPCT:ejemplo_2}.

		\begin{figure}[!h]
			\centering
			\includegraphics[width=0.99\textwidth]{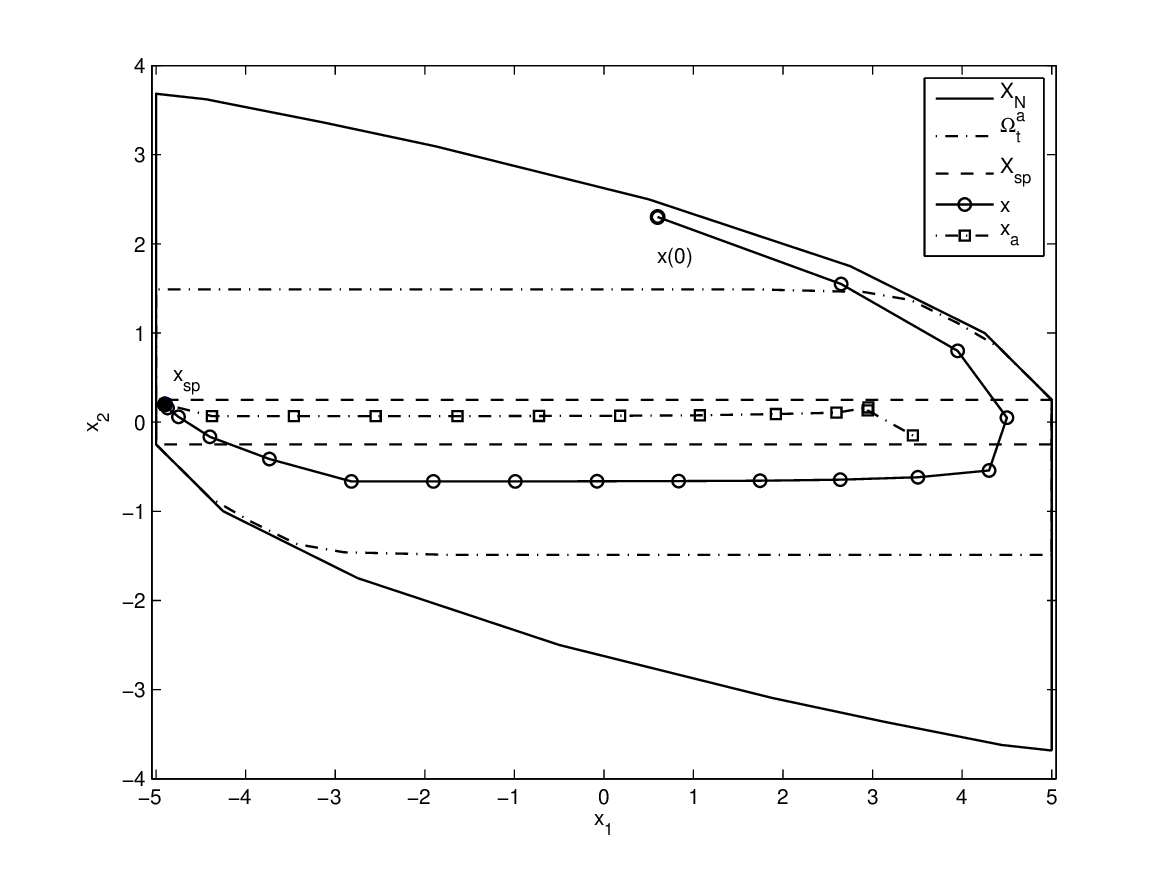}
			\caption{State space evolution of the complete simulation with terminal inequality constraint.}\label{chaMPCT:di_ic_state_space}
		\end{figure}
		
		\begin{figure}[!h]
			\centering
			\includegraphics[width=0.99\textwidth]{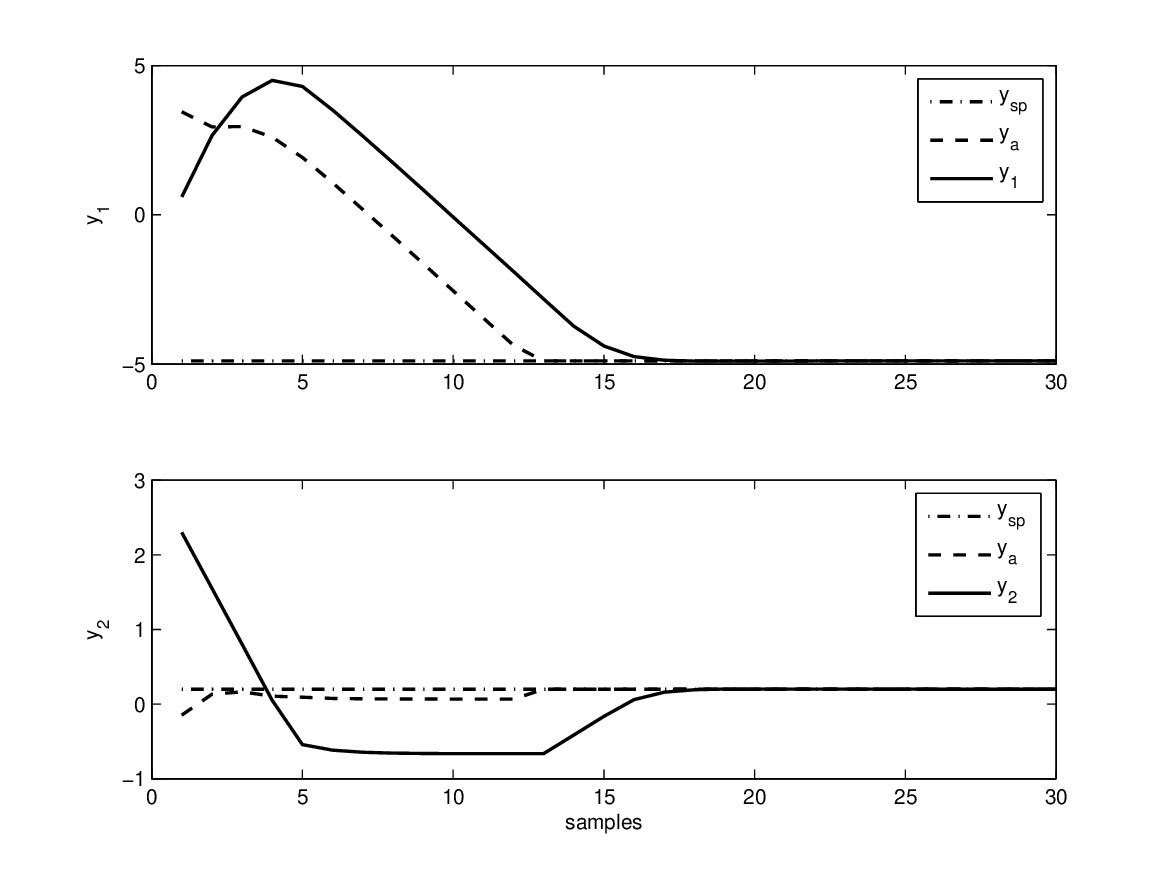}
			\caption{Time evolution of the complete simulation with terminal inequality constraint.}\label{chaMPCT:di_ic_y_evol}
		\end{figure}
		
		In Figure  \ref{chaMPCT:di_XN_comparison}, the feasible set of the MPCT with terminal inequality constraint $\setX^{ic}_N=\setX_N(\Omega_t)$ (solid line), compared to the feasible set of the MPC with terminal equality constraint $\setX^{ec}_N=\setX_N$ (dash-dotted line).
		
		Set $\setX^{ic}_N$ is clearly larger than set $\setX^{ec}_N$, since $\setX_s \subseteq \Omega_t$.
		
		\begin{figure}[!h]
			\centering
			\includegraphics[width=0.99\textwidth]{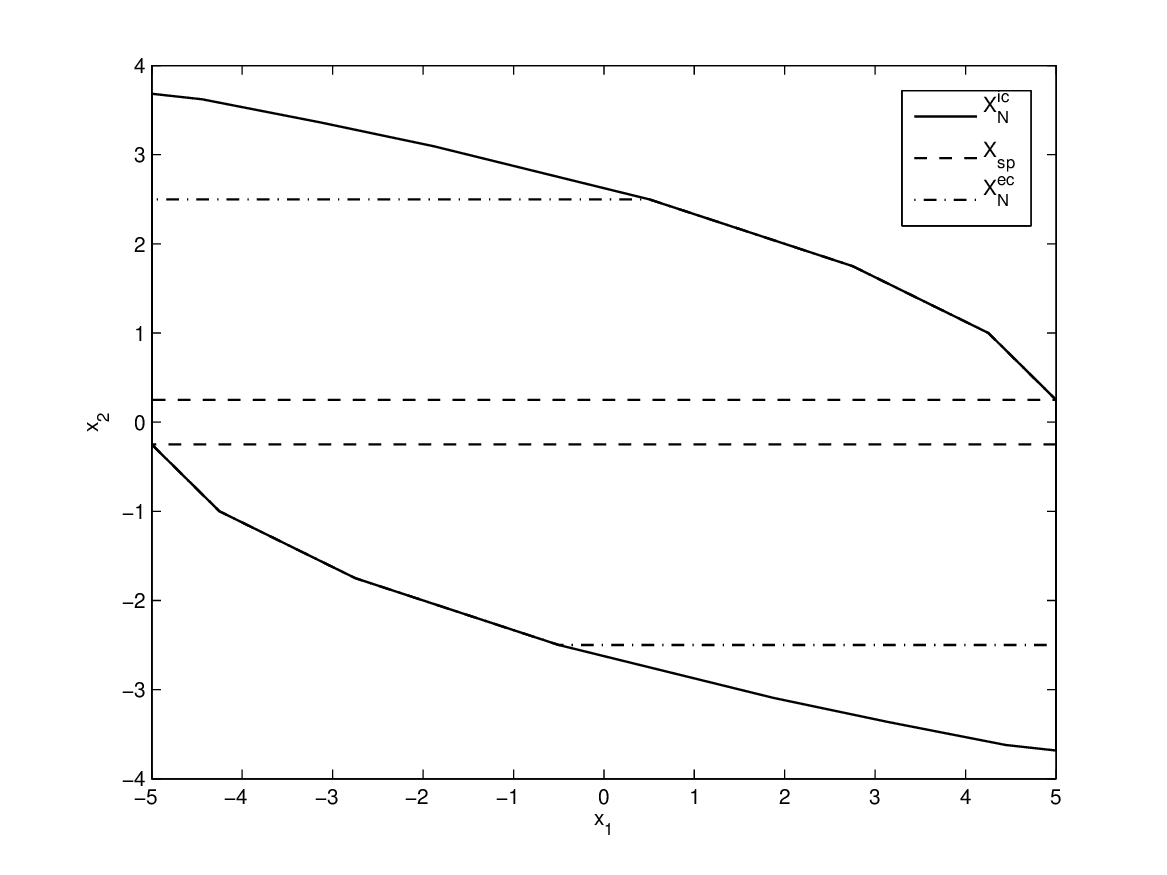}
			\caption{Feasible set of the MPCT with terminal inequality constraint in solid line, compared to the feasible set of the MPC with terminal equality constraint, in dash-dotted line.}\label{chaMPCT:di_XN_comparison}
		\end{figure}

	\end{exmp}

	\subsubsection{Calculation of the invariant set for tracking}\label{chaMPCT:invariante}
	
	Consider the following controller for a given admissible equilibrium point $(x_a,u_a) \in \setZ_{sp}$
	\beqn u=K(x-x_a)+u_a,
	\label{EqCtrlK_lin}
	\eeqn
	It is well known that if  $A+BK$ has all its eigenvalues
	inside the unit circle then the system is steered to the equilibrium point $(x_a,u_a)$. Since the system is constrained, this controller leads
	to an admissible evolution of the system only in a neighborhood of
	the steady state.
	
	Based on this control gain, a procedure to calculate the invariant set for tracking is derived. Let define
	the augmented state $\bx=(x,x_a,u_a)$. Assuming that the equilibrium point remains constant, the system controlled by the linear controller can be defined by the following equation
	\beqn \bmat{c} x \\
	x_a \\ u_a \emat^+ = \bmat{ccc} A+B K & -BK & B  \\ 0 & I_n & 0 \\ 0 & 0 & I_{m} \emat
	\bmat{c} x
	\\ x_a\\ u_a \emat \label{EqSysXsr} \eeqn
	that is, $\bx^+=\bar A \bx$. It may be seen that for a given initial augmented state $\bx(0)=(x(0),x_a,u_a)$, the  trajectory of the augmented autonomous system $\bx(k)$ is given by  $\bx(k)=(x(k),x_a,u_a)$, where $x(k)$ is the trajectory of the system controlled by (\ref{EqCtrlK_lin}) for the initial state $x(0)$.
	
	The evolution of the system must satisfy the constraints $(x_a,u_a) \in \setZ_{sp}$ and $(x,K(x-x_a) + u_a) \in \setZ$. This can be written in terms of the augmented state by the following set:
	\[
	\bar \setX=\{\bx=(x,x_a,u_a) \mid (x,K(x-x_a)+u_a) \in \setZ, \,
	x_a=A x_a + B u_a, (x_a,u_a) \in \lambda \setZ \}
	\]
	where $\lambda \in(0,1) $ is a constant arbitrarily close to 1 (see section \ref{chaMPCT:ReachSP}).
	
	Define, now the set
	\[
	\SetO_j=\{\bx \mid \bar A^i \bar x \in \bar\setX, \forall i \in [0,j] \}
	\]
	This is the set of initial augmented states such that the trajectory of the augmented system will be admissible along the first $j$ time instants. Notice that this set satisfies the following recursion
	$$\SetO_{j+1}= \SetO_j \cap \{\bx \mid \bar A^{j+1}\bx \in \bar \setX\} $$
	and then $\SetO_{j+1}\subseteq \SetO_j$. Besides $\SetO_0=\bar \setX$, which implies that $\SetO_j \subseteq \bar \setX$.
	
	The set $\SetO_\infty$ is called the maximal admissible invariant set and it satisfies that
	for all $\bx \in \SetO_\infty$, then $\bar A \bx \in \SetO_\infty$.  This is equivalent to say that for all $(x,x_a,u_a) \in \SetO_\infty \subset \bar \setX$, we have that $(x^+,x_a,u_a) \in \SetO_\infty$, where $x^+=A x+ B(K (x-x_a) + u_a)$. Therefore,  $\SetO_\infty$ is an invariant set for tracking.

	It has been proved in \cite{GilbertTAC91} that there exists a finite number $M$ such that the set $\SetO_\infty=\SetO_M$, and then this can be calculated in a finite number of steps of the recursion
	$$\SetO_{j+1}= \SetO_j \cap \{\bx \mid \bar A^{j+1}\bx \in \bar \setX\} $$
	with $\SetO_{0}=\bar \setX$. This recursion will stop when  the condition $\SetO_{j+1}= \SetO_j$ holds.
	Notice that the resulting set is a polyhedron.
	\section{Implementation: how to formulate the QP problem}\label{chaMPCT:implementation_QP}
	
	The optimization problems \eqref{chaMPCT:optprob_eqconst} and \eqref{chaMPCT:optprob_ineqconst}
	are convex mathematical programming problems, which can be
	efficiently solved using appropriate algorithms \cite{BoydLIB06}. For some realizations of  the offset cost function, problems \eqref{chaMPCT:optprob_eqconst} and \eqref{chaMPCT:optprob_ineqconst} can be posed as quadratic programming (QP) problems.
	In this section, we will show how to practically formulate problems \eqref{chaMPCT:optprob_eqconst} and \eqref{chaMPCT:optprob_ineqconst} as QPs.  To this aim, consider the canonical form of a QP problem \cite{BoydLIB06}, which can be written down as
	\beqn
	\begin{array}{lll}\label{chaMPCT:QP_canonico}
		&\displaystyle\min_{\vu_e}& \displaystyle\frac{1}{2}\vu_e'H\vu_e+f'\vu_e+r\\
		&s.t.&G\vu_e\leq W\\
		&&F\vu_e=S	
	\end{array}	
	\eeqn
	The objective of this section is to show how to obtain the ingredients necessary to pose Problems \eqref{chaMPCT:optprob_eqconst} and \eqref{chaMPCT:optprob_ineqconst} in the same form as Problem \eqref{chaMPCT:QP_canonico}.

	\subsection{MPCT with terminal inequality constraint}
	\textit{Cost function}
	
	Let us start by manipulating cost function \eqref{chaMPCT:costfuc_ineqconst}  and considering that the offset cost function is a quadratic function given by $V_O(y_a,y_{sp})=\|y_a-y_{sp}\|^2_T$, where $T$ is a suitable positive definite matrix. Then, function \eqref{chaMPCT:costfuc_ineqconst} can be written as:

	\beqnan\label{chaMPCT:cost_fucnt_ineqconst_QP}
	V_N(x,y_{sp};\vu,x_a,u_a)\!\!\!\!\!\!&=&\!\!\!\!\!\!\sum\limits_{j=0}^{N-1}\|x(j)\!-\!x_a\|^2_Q\!+\!\|u(j)\!-\!u_a\|^2_R+\|x(N)\!-\!x_a\|^2_P\!+\|y_a-y_{sp}\|^2_T
	\eeqnan

	%At time $k$, take $x(0)=x(k)$ and $u(0)=u(k)$.
	Define the sequence of predicted states and inputs as:
	\beqn
	\vx=\bmat{c}x(0)\\x(1)\\ x(2)\\ \vdots \\x(N-1)\\x(N) \emat,  \quad \vu=\bmat{c} u(0) \\ u(1) \\ u(2)\\ \vdots \\ u(N-2) \\u(N-1) \emat \nn
	\eeqn
	where $\vx \in \R^{(N+1)n}$, $\vu \in \R^{Nm}$. Then, taking into account the prediction model, the sequence of predicted states are given by
	\beqn
	\vx=\bold Ax(0)+\bold B \vu
	\eeqn
	
	where  $\vA \in \R^{(N+1)n \times n}$, and $\vB \in \R^{(N+1)n \times Nm}$ are given by
	
	\beqn
	\bold A=\bmat{c} I_n\\ A\\ A^2\\ \vdots \\ A^{N-1} \\A^{N}\emat, \quad
	\bold B=\bmat{llllll} 0 & 0 & 0 & \cdots & 0 & 0\\
	B & 0 & 0 & \cdots & 0 & 0\\
	AB & B & 0 & \cdots & 0 & 0\\
	\vdots & \vdots & \vdots & \cdots & \vdots & \vdots\\
	A^{N-2}B & A^{N-3}B & A^{N-4}B & \cdots & B & 0\\
	A^{N-1}B & A^{N-2}B & A^{N-3}B & \cdots & AB & B\emat \nn
	\eeqn
	
	Define also the following diagonal matrices $ \bold{Q} \in \R^{(N+1)n \times (N+1)n }$ and $ \bold{R} \in \R^{Nm \times Nm }$:
	
	\beqn
	\bold Q= \bmat{ccccc} Q & 0 & \cdots & 0 & 0 \\ 0 & Q & \cdots & 0 & 0 \\
	\vdots & \vdots & \ddots & \vdots & \vdots \\
	0 & 0 & \cdots & Q & 0 \\
	0 & 0 & \cdots & 0 & P \emat, \quad
	\bold R= \bmat{ccccc} R & 0 & \cdots & 0 & 0 \\ 0 & R & \cdots & 0 & 0 \\
	\vdots & \vdots & \ddots & \vdots & \vdots \\
	0 & 0 & \cdots & R & 0 \\
	0 & 0 & \cdots & 0 & R \emat \nn
	\eeqn
	
	Notice that, in Problem \eqref{chaMPCT:optprob_ineqconst}, the decision variables are $\vu_e=(\vu,x_a,u_a)$. Therefore, define:
	
	\beqn
	%\vu_e=\bmat{c} \vu \\ x_a \\ u_a \emat, \quad
	\bold B_e=\bmat{ccc} \bold B, & -I_{x_a}, & \mathbf{0}_{(N+1)n\times m}\emat, \quad \bold I_e=\bmat{ccc} I_{Nm},& \mathbf{0}_{Nm\times n} & -I_{u_a}\emat \nn \eeqn
	where $\mathbf{0}_{r\times l} \in \R^{r\times l}$ is a matrix with all zero elements, and  $ I_{x_a} \in \R^{(N+1)n\times n}$ and $ I_{u_a} \in \R^{Nm \times m}$, are given by
	\beqn
	\begin{array}{cc}
		I_{x_a}=\bmat{c} I_n \\ I_n \\ \vdots \\ I_n\emat, & I_{u_a}=\bmat{c} I_m \\I_m\\ \vdots \\I_m \emat
	\end{array}\nn	
	\eeqn
	Taking into account that $y_a=Cx_a+Du_a$, define also $F_e=\bmat{ccc}\mathbf{0}_{p \times Nm}, & C, & D\emat$.
	
	Given all these ingredients, we can now rewrite function $V_N(x,y_{sp};\vu, x_a, u_a)$ as:
	
	\beqnan
	V_N(x,y_{sp};\vu_e)\!\!\!\!&=&\!\!\!\!(\vA x+\vB_e\vu_e)'\vQ(\vA x+\vB_e\vu_e)+\vu_e'\vI_e'\vR\vI_e\vu_e\\
	\!\!\!\!&&\!\!\!\!+(F_e\vu_e-y_{sp})'T(F_e\vu_e-y_{sp})\\
	\!\!\!\!&=&\!\!\!\!x'\vA' \vQ \vA x + \vu_e'\vB_e' \vQ \vB_e\vu_e + 2 x'\vA' \vQ \vB_e \vu_e+\vu_e'\vI_e'\vR\vI_e\vu_e\\
	\!\!\!\!&&\!\!\!\!+\vu_e'F_e'TF_e\vu_e+ y_{sp}'Ty_{sp}-2y_{sp}'TF_e\vu_e
	\eeqnan
	
	Rearranging the equality above, we can rewrite function \eqref{chaMPCT:costfuc_eqconst} as
	
	\beqn\label{chaMPCT:costfunc_QP}
	V_N(x,y_{sp};\vu_e)=\frac{1}{2}\vu_e'H\vu_e+f'\vu_e+r
	\eeqn
	where
	\beqnan
	H&=& 2(\vB_e' \vQ \vB_e+\vI_e'\vR\vI_e+F_e'TF_e)\\
	f&=& 2(x'\vA' \vQ \vB_e-y_{sp}'TF_e)\\
	r&=& x'\vA' \vQ \vA x + y_{sp}'Ty_{sp}
	\eeqnan
	
	\textit{Constraints}
	
	Let us now take into account the constraints to problem \eqref{chaMPCT:optprob_ineqconst}. First of all, notice that (i) $x(0)=x$, (ii) $x(j+1)=A x(j)+B u(j)$, and (iii) $y_a=Cx_a+Du_a$, are actually taken into account in the ingredients defined above. So there is no need to reconsider them again. Let us now focus on the inequality constraints.
	
	Let us consider the inequality constraint $(x(j),u(j)) \in \setZ, \quad  \!j\!=\!0,\cdots, N\!-\!1$. $\setZ$ is a set of linear constraints. It usually represents upper and lower bound on $x$ and $u$, but it can also represent bounds on some linear combination of them. This set can be written in the form of a linear inequality as
	\beqn
	\bmat{cc} \tilde G_{\setZ,x}, & \tilde G_{\setZ,u} \emat \bmat{c} x\\ u\emat \leq W_{\setZ}\nn
	\eeqn
	Considering the entire sequences of future states and future inputs, we have:
	\beqn
	\bmat{cc} G_{\setZ,x}, & G_{\setZ,u} \emat \bmat{c} \vx\\\vu\emat \leq W_{\setZ}\nn
	\eeqn
	where
	\beqn
	G_{\setZ,x}= \bmat{cccccc} \tilde G_{\setZ,x} & 0 & \cdots & 0 & 0 & 0 \\ 0 & \tilde G_{\setZ,x} & \cdots & 0 & 0 & 0\\
	\vdots & \vdots & \ddots & \vdots & \vdots & \vdots \\
	0 & 0 & \cdots & \tilde G_{\setZ,x} & 0 & 0 \\
	0 & 0 & \cdots & 0 & \tilde G_{\setZ,x} & 0 \emat \nn,\eeqn
		\beqn G_{\setZ,u}= \bmat{ccccc} \tilde G_{\setZ,u} & 0 & \cdots & 0 & 0 \\ 0 & \tilde G_{\setZ,u} & \cdots & 0 & 0 \\
	\vdots & \vdots & \ddots & \vdots & \vdots \\
	0 & 0 & \cdots & \tilde G_{\setZ,u} & 0 \\
	0 & 0 & \cdots & 0 & \tilde G_{\setZ,u} \emat \nn
	\eeqn

	Recalling that $\vx=\bold Ax(0)+\bold B \vu$, the previous inequality can be rewritten as:
	\beqn\label{chaMPCT:QP_ineqconst_1}
	(G_{\setZ,x}\mathbf{B}+G_{\setZ,u})\vu \leq W_{\setZ} - G_{\setZ,x}\mathbf{A}x(0)
	\eeqn

	Let us consider now, the terminal inequality constraint $(x(N),x_a,u_a) \in \Omega^a_{t}$. This constraint is also a set of linear inequalities of the form
	\beqn
	\bmat{ccc} G_{\Omega^a_{t},x}, & G_{\Omega^a_{t},x_a} & G_{\Omega^a_{t},u_a} \emat \bmat{c} x(N)\\ x_a \\ u_a\emat \leq W_{\Omega^a_{t}}\nn
	\eeqn
	Recalling that $x(N)= A^Nx(0)+\bold B_N \vu$, with
	\[\bold B_N=\bmat{llllll} A^{N-1}B & A^{N-2}B & A^{N-3}B & \cdots & AB & B\emat \]
	the previous inequality can be rewritten as:
	\beqn\label{chaMPCT:QP_ineqconst_2}
	\bmat{ccc} G_{\Omega^a_{t},x}\mathbf{B}_N, & G_{\Omega^a_{t},x_a} & G_{\Omega^a_{t},u_a} \emat \bmat{c} \vu\\ x_a \\ u_a\emat \leq W_{\Omega^a_{t}}-G_{\Omega^a_{t},x}A^Nx(0)
	\eeqn
	Combining \eqref{chaMPCT:QP_ineqconst_1} and \eqref{chaMPCT:QP_ineqconst_2} and taking $\vu_e=(\vu,x_a,u_a)$, we can pose the inequality constraints of problem \eqref{chaMPCT:optprob_ineqconst} in the form $$G\vu_e\leq W$$ where
	\beqn
	\begin{array}{cc}
		G=\bmat{ccc} G_{\setZ,x}\mathbf{B}+G_{\setZ,u} & \mathbf{0}_{n_z \times n} & \mathbf{0}_{n_z \times m} 	\\ G_{\Omega^a_{t},x}\mathbf{B}_N, & G_{\Omega^a_{t},x_a} & G_{\Omega^a_{t},u_a}\emat, &
		\mbox{ and } W=\bmat{c} W_{\setZ} - G_{\setZ,x}\mathbf{A}x(0) \\ W_{\Omega^a_{t}}-G_{\Omega^a_{t},x}A^Nx(0) \emat	 
	\end{array}	
	\eeqn
	where $n_z$ is the number of linear inequalities that define $\setZ$.
	
	Notice that, in the MPCT with terminal inequality constraint, we do not have equality constraints of the form $F\vu_e=S$.

	\begin{rem}
		If the offset cost function is a $\infty$-norm (or a 1-norm), i.e. $V_O(y_a,y_{sp})=\|y_a-y_{sp}\|_\infty$ (or $V_O(y_a,y_{sp})=\|y_a-y_{sp}\|_1$), the optimization problem can still be formulated as a quadratic programming, by posing the $\infty$-norm (or the 1-norm) in an epigraph form: that is, the offset cost function is actually taken as $V_O(y_a,y_{sp})=\lambda$, where $\lambda \in \R$ is an auxiliary optimization variable, and two constraints are added to the optimization problem: (i) $\|y_a-y_{sp}\|_\infty\leq \lambda$ (or $\|y_a-y_{sp}\|_1\leq \lambda$), and (ii) $\lambda \geq 0$.
		
		Function $V_N(x,y_{sp};\vu_e)$ can still be written as in Equation \eqref{chaMPCT:costfunc_QP}, but in this case:
		
		\beqn
		\begin{array}{cccc}
			\vu_e=\bmat{c}\! \vu \\ x_a \\ u_a \\ \lambda \!\emat, & H=\bmat{cc}\! 2(\vB_e' \vQ \vB_e+\vI_e'\vR\vI_e) & 0\\ 0 & 0\!\emat, &
			f= 2x'\vA' \vQ \vB_e+I_\lambda, & r= x'\vA' \vQ \vA x
		\end{array}\nn
		\eeqn
		where $I_\lambda=\bmat{cc} \mathbf{0}_{1 \times Nm+n+m}, & 1\emat$.
		
		The inequality constraint $G\vu_e\leq W$ (for instance in case of a $\infty$-norm) may be modified as follow:
		\beqn
		\begin{array}{cc}
			G=\bmat{cccc} G_{\setZ,x}\mathbf{B}+G_{\setZ,u} & \mathbf{0}_{n_z \times n} & \mathbf{0}_{n_z \times m} & \mathbf{0}_{n_z \times 1}	\\
			G_{\Omega^a_{t},x}\mathbf{B}_N, & G_{\Omega^a_{t},x_a} & G_{\Omega^a_{t},u_a} & \mathbf{0}_{n_\Omega \times 1}	\\
			\mathbf{0}_{p \times Nm} & C & D & -\mathbf{1}_{p \times 1}	\\
			\mathbf{0}_{p \times Nm} & -C & -D & -\mathbf{1}_{p \times 1} \\
			\mathbf{0}_{1 \times Nm} & \mathbf{0}_{1 \times n} & \mathbf{0}_{1 \times m} & -1 \emat, &
			\mbox{ and } W=\bmat{c} W_{\setZ} - G_{\setZ,x}\mathbf{A}x(0) \\ W_{\Omega^a_{t}}-G_{\Omega^a_{t},x}A^Nx(0) \\ y_{sp} \\ -y_{sp} \\ 0\emat \nn	
		\end{array}	
		\eeqn
		where $n_z$ is the number of linear inequalities that define $\setZ$, $n_\Omega$ is the number of linear inequalities that define $\Omega^a_{t}$, and $\mathbf{1}_{r \times 1} \in \R^r$ is an array with all its elements equal to 1.
		
		This solution can actually be adopted for any offset cost function such that the region $\{y \mid \, V_O(y,y_{sp})\leq \lambda\}$ is  polyhedral for any $\lambda>0$. \QED
	\end{rem}
	
	\begin{rem}
		From a practical point of view, in order to reduce the number of optimization variables, the steady state and
		input $(x_a,u_a)$ can be parameterized as a linear
		combination of a vector $\theta \in \mathbb R^m$, that is
		\beqn \label{EqSolSs} (x_a,u_a)=M_\theta \theta \eeqn
		where matrix $M_\theta$ is such that $$\bmat{cc}(A\!-\!I_n) & B\emat M_\theta=0$$ The
		steady outputs are given by
		\beqn \label{EqSolSs_y} y_a=N_\theta \theta \eeqn
		where $N_\theta=\bmat{cc}C & D\emat M_\theta$ (see section \ref{chaMPCT:ReachSP}).
		
		Then, in Problem \eqref{chaMPCT:optprob_ineqconst}, the decision variables $(\vu,x_a,u_a)$, become $(\vu,\theta)$. Therefore:
		
		\beqn
		\vu_e=\bmat{c} \vu \\ \theta \emat, \quad \bold B_e=\bmat{ccc} \bold B, & -I_{x_a}M_{\theta,x}\emat, \quad \bold I_e=\bmat{ccc} I_{Nm},&  -I_{u_a}M_{\theta,u}\emat \nn \eeqn
		where $M_{\theta,x}$ and $M_{\theta,u}$ are such that
		\beqn
		\bmat{c} x_a \\ u_a  \emat=\bmat{c} M_{\theta,x} \\ M_{\theta,u}  \emat \theta
		\nn	
		\eeqn

		As for the invariant set for tracking, notice that the terminal control law $u=K(x-x_a)+u_a$ can be rewritten as
		
		\beqnan u&=& K x + \bmat{cc}-K & I_m\emat \bmat{c} x_a \\ u_a \emat \\
		&=&K x + \bmat{cc}-K & I_m\emat M_\theta \theta\\
		&=& K x+L \theta \eeqnan
		where $L=\bmat{cc}-K & I_m\emat M_\theta \in \mathbb{R}^{m \times m}$. Consider
		the augmented state $\bx=(x,\theta)$, then the closed-loop augmented
		system can be defined by the following equation
		\beqn \bmat{c} x \\
		\theta \emat^+ = \bmat{cc} A+B K & B L  \\ 0 & I_{m} \emat
		\bmat{c} x
		\\ \theta \emat \label{EqSysXsr} \eeqn
		that is, $\bx+=\bar A \bx$.  The set of constraints is defined as $\bar \setX=\{(x,\theta) \mid (x,Kx + L \theta) \in \setZ, M_\theta \theta \in \lambda \setZ\}$. Then, the invariant set for tracking is calculated as proposed in section \ref{chaMPCT:invariante}.
		The terminal inequality constraint becomes $(x(N),\theta) \in \Omega^a_{t}$, and can be posed as a set of linear inequalities of the form
		\beqn
		\bmat{cc} G_{\Omega^a_{t},x}, & G_{\Omega^a_{t},\theta} \emat \bmat{c} x(N)\\ \theta\emat \leq W_{\Omega^a_{t}}\nn
		\eeqn
		Recalling that $x(N)= A^Nx(0)+\bold B_N \vu$, with
		\[\bold B_N=\bmat{llllll} A^{N-1}B & A^{N-2}B & A^{N-3}B & \cdots & AB & B\emat \]
		the previous inequality can be rewritten as:
		\beqn\label{chaMPCT:QP_ineqconst_theta}
		\bmat{cc} G_{\Omega^a_{t},x}\mathbf{B}_N, & G_{\Omega^a_{t},\theta} \emat \bmat{c} \vu\\ \theta \emat \leq W_{\Omega^a_{t}}-G_{\Omega^a_{t},x}A^Nx(0)
		\eeqn

		The dimension of $\theta$ is $m$, which is the dimension of the
		subspace of steady states and inputs that can be parameterized by a
		minimum number of variables. Hence, equations \eqref{EqSolSs} and \eqref{EqSolSs_y}
		represent a mapping of $(x_a,u_a)$ and $y_{a}$ onto the subspace of
		$\theta$. The set of setpoints $y_{sp}$ that can be admissibly reached
		is the subspace spanned by the columns of $N_\theta$.  \QED
	\end{rem}

	\subsection{MPCT with terminal equality constraint}
	The formulation with terminal equality constraint is a particular case of the previous one. Let us start with the cost function. In this case, with the offset cost function given by $V_O(y_a,y_{sp})=\|y_a-y_{sp}\|^2_T$, where $T$ is a suitable matrix, function \eqref{chaMPCT:costfuc_eqconst} can be written as:
	\beqna\label{chaMPCT:cost_fucnt_eqconst_QP}
	V_N(x,y_{sp};\vu,x_a,u_a)&=&\!\!\sum\limits_{j=0}^{N-1}\|x(j)\!-\!x_a\|^2_Q\!+\!\|u(j)\!-\!u_a\|^2_R+\|y_a-y_{sp}\|^2_T
	\eeqna
	
	All the ingredients introduced in the previous section are the same, but $ \bold{Q} \in \R^{(N+1)n \times (N+1)n }$, which is now given by
	\beqn
	\bold Q= \bmat{ccccc} Q & 0 & \cdots & 0 & 0 \\ 0 & Q & \cdots & 0 & 0 \\
	\vdots & \vdots & \ddots & \vdots & \vdots \\
	0 & 0 & \cdots & Q & 0 \\
	0 & 0 & \cdots & 0 & 0 \emat \nn
	\eeqn
	since in this formulation there is no cost-to-go from $N$ to $\infty$.
	
	Taking into account this fact, we can once again rewrite function \eqref{chaMPCT:cost_fucnt_eqconst_QP} as:
	\beqn
	V_N(x,y_{sp};\vu_e)=\frac{1}{2}\vu_e'H\vu_e+f'\vu_e+r \nn
	\eeqn
	with
	\beqnan
	H&=& 2(\vB_e' \vQ \vB_e+\vI_e'\vR\vI_e+F_e'TF_e)\\
	f&=& 2(x'\vA' \vQ \vB_e-y_{sp}'TF_e)\\
	r&=& x'\vA' \vQ \vA x + y_{sp}'Ty_{sp}
	\eeqnan

	As for the constraints, as in the previous case, (i) $x(0)=x$, (ii) $x(j+1)=A x(j)+B u(j)$, and (iii) $y_a=Cx_a+Du_a$, are taken into account in the ingredients that define the cost function. Since there is no terminal inequality constraint, the inequality constraints of problem \eqref{chaMPCT:optprob_eqconst} can be posed in the form $$G\vu_e\leq W$$ with
	\beqn
	\begin{array}{cc}
		G=\bmat{ccc} G_{\setZ,x}\mathbf{B}+G_{\setZ,u} & \mathbf{0}_{n_z \times n} & \mathbf{0}_{n_z \times m} \emat, &
		\mbox{ and } W=W_{\setZ} - G_{\setZ,x}\mathbf{A}x(0)
	\end{array}	
	\eeqn
	where $n_z$ is the number of linear inequalities that define $\setZ$.
	
	Finally, it remains to formulate the equality constraints:
	$$x_a=Ax_a+Bu_a, \quad x(N)=x_a$$
	in the form $F\vu_e=S$.
	
	Recalling that $x(N)= A^Nx(0)+\bold B_N \vu$, this can easily be done by taking:
	\beqn
	F=\bmat{lll} \mathbf{0}_{n \times Nm}, & I_n-A, & -B \\ \mathbf{B}_N, & -I_n, &  \mathbf{0}_{n \times m}\emat, \mbox{ and } S=\bmat{c} \mathbf{0}_{n \times 1} \\ -A^Nx(0)\emat \nn
	\eeqn

	\begin{rem}
		Also in this case, if the offset cost function is a $\infty$-norm (or a 1-norm), i.e. $V_T(y_a,y_{sp})=\|y_a-y_{sp}\|_\infty$ (or $V_T(y_a,y_{sp})=\|y_a-y_{sp}\|_1$), the optimization problem can be formulated as a quadratic programming, by posing the $\infty$-norm (or the 1-norm) in an epigraph form. \QED
	\end{rem}

	\subsection{Off-line implementation as an explicit controller.}\label{chaMPCT:prop_explicit}
	The control law derived by the solution of problems \eqref{chaMPCT:optprob_eqconst} and \eqref{chaMPCT:optprob_ineqconst}, $\kappa_N(x,y_{sp})$, is a function of the parameters $(x,y_{sp})$, which appear in the cost function (both), and in the constraints (only $x$).
	
	Moreover, as shown in this Section, these problems can be posed as a Quadratic Programming problem like \eqref{chaMPCT:QP_canonico}, whose ingredients can be easily calculated.
	
	Taking into account these two facts, we can express the cost function $V_N(x,y_{sp};\vu_e)$ and the optimal solution $\kappa_N(x,y_{sp})$ as en explicit function of the parameter $(x,y_{sp})$.
	
	To this aim, let the offset cost function be given by $V_O(y_a,y_{sp})=\|y_a-y_{sp}\|^2_T$, where $T$ is a suitable matrix. Let us define
	$$\vz=\vu_e+L_x x+ L_{sp} y_{sp}$$
	where $L_x=H^{-1}f_x'$, $L_{sp}=H^{-1}f_{sp}'$, and $f_x=2\vA'\vQ'\vB_e$ and $f_{sp}=-2TF_e$. Then, problem
	\eqref{chaMPCT:QP_canonico} can be rewritten as
	\beqn
	\begin{array}{lll}\label{chaMPCT:QP_mp}
		&\displaystyle\min_{\vz}& \displaystyle\frac{1}{2}\vz'H\vz\\
		&s.t.&G\vz\leq \bar W+W_x x +W_{sp} y_{sp} \\
		&&F\vz=\bar S + S_x x + S_{sp} y_{sp} 	
	\end{array}	
	\eeqn
	where
	\begin{itemize}
		\item in case of terminal inequality constraint we have $W_x=\tilde W + GL_x$, and $W_{sp}=GL_{sp}$, with
		\beqn
		\bar W = \bmat{c} W_{\setZ} \\ W_{\Omega^a_{t}}\emat, \quad \tilde W=\bmat{c} - G_{\setZ,x}\mathbf{A} \\ -G_{\Omega^a_{t},x}A^N \emat\nn
		\eeqn
		and no equality constraint;
		
		\item in case of terminal equality constraint we have $W_x=\tilde W + GL_x$, and $W_{sp}=GL_{sp}$, with
		\beqn
		\bar W =  W_{\setZ}, \quad \tilde W=- G_{\setZ,x}\mathbf{A} 	\nn
		\eeqn
		and $S_x=\tilde S + FL_x$, and $S_{sp}=FL_{sp}$, with
		\beqn
		\bar S =  \mathbf{0}_{2n\times 1}, \quad \tilde S= S 	\nn
		\eeqn
	\end{itemize}	
	
	Problem \eqref{chaMPCT:QP_mp} is a multiparametric Quadratic Programming (mp-QP).
	The structure of this problem shows that the
	control law $\kappa_N(x,y_{sp})$ is a piecewise affine function of the parameters
	$(x,y_{sp})$, and it can by explicitly calculated off-line, by means of the
	existing multiparametric programming tools.  Thus,  the feasibility region $\setX_N$ can be divided in a collection of disjoints polyhedrons $\Gamma_j$ such that
	$$
	\setX_N=\bigcup_{j} \Gamma_j
	$$
	The control law can be posed as a piecewise affine function defined at each of these regions, that is, for all $x \in \Gamma_j$, the control law is given by
	$$\kappa_N(x,y_{sp})=K^x_j x + K_j^y y_{sp} + c_j$$
	
	The partition and  the matrices of the control law can be calculated by means of specialized algorithms \cite{BemporadAUT02}.
	
	Notice that, thanks to the calculation of the explicit control law,  the range of applicability of the MPC for tracking can also include those situations where the on-line computation of the MPC control law may be prohibitive, such as those arising in the automotive and aerospace industries \cite{BemporadAUT02}.

	% % % % % % % % % % % % % % % % % % % % % % % % % % % % % % % % % % % % % % % % % % % % % % % % % % % %
	\section{Properties of the proposed controller}\label{chaMPCT:properties}

	As it has been proved, the proposed controller guarantees stability and convergence to the setpoint when this is reachable. Besides this controller has a number of interesting properties that are highlighted hereafter.
	
	\begin{enumerate}
		
		\item \textbf{Stability under changing setpoints.}
		
		The control law $u=\kappa_N(x,y_{sp})$ is derived from the solution of the optimization problem $P_N(x,y_{sp})$. Since the set of constraints of this optimization problem does not depend on the setpoint $y_{sp}$, feasibility cannot be lost due to  changes in the setpoint, even when the change is significant. Therefore, if the setpoint is changed to a new reachable setpoint, the controller will be well posed and  will steer the plant to the new setpoint.\\
		
		\item \textbf{Unreachable setpoints} \label{chaMPCT:prop_SSTO}
		
		It is not unusual that the given setpoint $y_{sp}$  is not reachable, that is, $y_{sp}$ is not contained in $\setY_{sp}$. For instance, when the setpoint is provided by the Real Time Optimizer (RTO),  this may be not reachable  due to the difference between the nonlinear model used in the RTO and the linear prediction model used in the MPC.
		%To deal with this situation, the standard solution in the hierarchical control structure is to add a Steady State Target Optimizer (SSTO) between the RTO and the MPC, to calculate the best admissible target for the controller \cite{YingAICHE1999performance}.
		
		The proposed controller will be feasible when the provided setpoint is not reachable and besides it will steer the system to a reachable equilibrium point $(x_t,u_t,y_t)$ such that
		$$ y_t=\arg \min_{y \in \setY_{sp}} V_O(y,y_{sp})$$
		
		Then, in case of unreachable setpoints, the controlled system will be steered to the reachable equilibrium point that minimizes the offset cost function. This property constitutes a criterion for the selection of the offset cost function, as long as this function is convex, positive definite and subdifferentiable.\\

		%The special formulation of the MPC for tracking  drives the system to the optimal operating
		%point according to the offset cost function $V_O(\cdot,\cdot)$, that is, steers the system to the admissible steady output $y_t$ such that
		%$$ y_t=\arg \min_{y \in \setY_{sp}} V_O(y,y_{sp})$$
		%
		%Then it can be
		%considered that the proposed controller integrates the SSTO in its own formulation and $V_O(\cdot,\cdot)$ defines the function to optimize.

		%\item \textbf{Offset cost function and stability.}\label{chaMPCT:prop_eco}
		%The stability results in Theorems \ref{chaMPCT:TeorEst_eqconst} and \ref{chaMPCT:TeorEst_ineqcons}, always hold for any offset cost function satisfying Definition \ref{chaMPCT:def_optimo}. Therefore, if $V_O(\cdot,\cdot)$ varies with
		%the time, the results of these Theorems still hold.\\
		%This property allows one to tune the cost function in such a way that it provides certain (economic) properties to the MPC for tracking. This property will be clearer in the following chapters.
		
		\item \textbf{Larger domain of attraction.}\label{chaMPCT:prop_domain}
		
		The domain of attraction of the MPC for a given prediction horizon is the set of states that can be admissible steered terminal set in $N$ steps. Since the terminal set designed for the MPC for tracking $\Omega_t$ is potentially much larger than the terminal set designed for the MPC for regulation, the domain of attraction of MPC for tracking is potentially much larger.
		
		Consider also, the case of using a terminal equality constraint. The terminal set for the MPC for regulation is the steady state for the setpoint $x_{sp}$ while the terminal set of the MPC for tracking is the set of all reachable steady states $\setX_{sp}$, which is much larger.
		
		Another remarkable property is that for any prediction horizon, every admissible steady states is contained in the domain of attraction of the MPC for tracking. That is
		
		$$ \setX_{sp} \subseteq \setX_N, \qquad \forall N\geq 1$$
		
		This means that if the initial state is an admissible equilibrium point, the proposed MPC for tracking will steer the system to any admissible setpoint irrespective of the prediction horizon.
		
		Notice that in practice, the plant to be controlled is typically  operated manually to an admissible equilibrium point and then the controller operates the plant in closed loop. Therefore the initial state for the controller is an admissible equilibrium point, which guarantees that the control low is feasible and well-posed.

		This property was illustrated  in Example \ref{chaMPCT:ejemplo_1}, where it was shown that the proposed controller, for a given horizon, provides a larger feasible set then standard MPC.
		
		\begin{figure}[!h]
			\centering
			\includegraphics[width=0.99\textwidth]{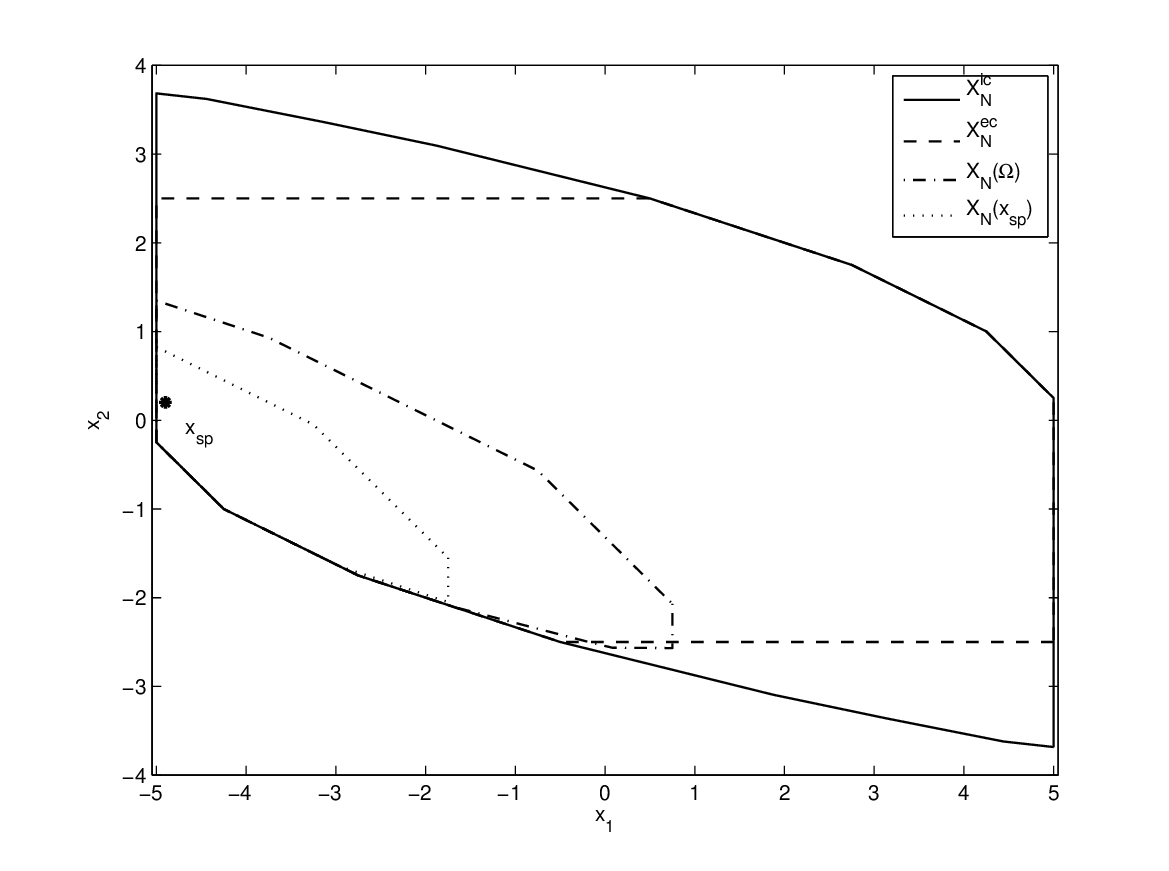}
			\caption{Comparison of the feasible sets of MPC with terminal equality constraint $x(N)=x_{sp}$ ($\setX_N(x_{sp})$, dotted line), MPC with terminal inequality constraint $x(N) \in \Omega$ ($\setX_N(\Omega)$, dash-dotted line), MPC with relaxed terminal equality constrain $x(N)=x_a$ ($\setX_N^{ec}$, dashed line), MPC with relaxed terminal inequality constrain $x(N) \in \Omega_t$ ($\setX_N^{ic}$, solid line).}\label{chaMPCT:di_feasible_set_comp}
		\end{figure}
		
		Figure \ref{chaMPCT:di_feasible_set_comp} shows a comparison of feasible sets for the case study in Example \ref{chaMPCT:ejemplo_1}. The feasible set for and MPC with terminal equality constraint $x(N)=x_{sp}$ is drawn in dotted line ($\setX_N(x_{sp})$), and the one for an MPC with terminal inequality constraint $x(N) \in \Omega$ is plotted in dash-dotted line ($\setX_N(\Omega)$). At the same time, the feasible set for an MPC with relaxed terminal equality constrain $x(N)=x_a$ is drawn in dashed line ($\setX_N^{ec}$), and the one for an MPC with relaxed terminal inequality constrain $(x(N),x_a,u_a) \in \Omega_t^a$ is plotted in solid line ($\setX_N^{ic}$). It is evident how for a same prediction horizon, the relaxed terminal constraint provides a feasible region larger than standard MPC formulation.

		This property makes the
		proposed controller interesting even for regulation objectives.
		
		\item \textbf{Robustness and output feedback}\label{chaMPCT:prop_robustness}
		
		It has been shown in \cite{GrimmAUT04} that asymptotically stabilizing predictive
		control laws may exhibit zero-robustness, that is, any disturbance
		may lead the MPC controller to lose feasibility or asymptotic stability. In the case of the MPC for tracking presented in this report, taking into
		account that the control law is derived from a multiparametric
		convex problem, the closed-loop system is input-to-state stable under
		sufficiently small
		uncertainties \cite{LimonISS09}.\\
		This property is very interesting
		for an output feedback formulation \cite{MessinaAUT05}, since it
		allows to ensure asymptotic stability of the closed-loop, under a control law based on
		an estimated state, using an asymptotically stable observer.
	\end{enumerate}
	
	%%%%%%%%%%%%%%%%%%%%%%%%%%%%%%%%%%%%%%%%%%%%%%%%%%%%%%%%%%%%%%%%%%%%%%%%%%%%%%%%
	\section{Local Optimality}\label{chaMPCT:local_optimality}
	%The aim of this section is to present the property of local optimality of an MPC control law, and how this property can be enjoyed by the MPC for tracking.

	Consider that system \eqref{sistema_o} is controlled by the control
	law $u=\kappa(x,y_{sp})$, in order to steer the system to the setpoint $y_{sp}\in
	\setY_{sp}$. According to the considered quadratic stage cost function, the performance of the controlled system can be measured by means of the following cost-to-go function: 
	
	\beqn  
	V_\infty(x,y_{sp}; \kappa(\cdot,y_{sp}))\!=\!\sum\limits_{j=0}^{\infty}\|x(j)\!-\!x_{sp}\|^2_Q + \|\kappa(x(j),y_{sp})\!-\!
	u_{sp}\|^2_R \label{cha1:Vinfty}
	\eeqn
	where $x(j)=\phi(j;x,\kappa(\cdot,y_{sp}))$ is calculated from the
	recursion $x(i+1)=Ax(i)+B\kappa(x(i),y_{sp})$ for $i=0,\cdots, j-1$
	with $x(0)=x$. A control law $\kappa_\infty(x,y_{sp})$ is said to be
	optimal if it is admissible (namely, the constraints are fulfilled
	throughout the closed-loop evolution) and it is the one that minimizes
	the cost $V_\infty(x,y_{sp};\kappa(\cdot,y_{sp}))$ for all admissible $x$. Let us denote the optimal cost function as
	$$V_\infty^0(x,y_{sp})=V_\infty(x,y_{sp};\kappa_\infty(\cdot,y_{sp}))$$
	From a practical point of view it is very interesting to find the optimal control law since 
	this is the one that provide the best possible closed-loop performance measured by the function (\ref{cha1:Vinfty}), $V_\infty^0(x,y_{sp})$.

	It is well known that the constrained  Linear
	Quadratic Regulator is the optimal control law to be designed
	according to the given quadratic performance index.  While the optimal control law  for an
	unconstrained system is easily obtained by solving the Riccati's
	equation, its
	calculation in case of a constrained system may be
	computationally not affordable. Model predictive control can be considered as suboptimal
	since the cost function is only minimized over a finite
	prediction horizon, but its optimality is enhanced as long as the prediction horizon is enlarged.
	
	In order to study the optimality of the proposed controller, it is convenient to write the optimization problem of the MPC control law to regulate the
	system to the target $y_{sp}$, $\kappa^r_N(x,y_{sp})$, as follows
	\begin{subequations}\label{chaMPCT:prob_regulation}
		\beqna
		V_N^{r,0}(x,y_{sp})\!\!&\!\!=\!\!&\!\!\min\limits_{\vu,x_a,u_a}\sum\limits_{j=0}^{N-1}\|x(j)\!-\!x_a\|^2_Q\!+\!\|u(j)\!-\!u_a\|^2_R+\|x(N)\!-\!x_a\|^2_P \\
		&s.t.& x(0)=x,\\
		&& x(j+1)=A x(j)+B u(j), \\
		&& (x(j),u(j)) \in \setZ,  \!\quad\!  j=0,\cdots, N\!\!-\!\!1 \\
		&& y_a=Cx_a+Du_a,\\
		&& (x(N),x_a,u_a) \in \Omega^a_t\\
		&& \|y_a-y_{sp}\|_{\infty}=0\label{chaMPCT:prob_regulation_eq_constr}
		\eeqna
	\end{subequations}
	Notice that this is similar to the optimization problem of the MPC for tracking, but with the additional constraint (\ref{chaMPCT:prob_regulation_eq_constr}) that forces the artificial equilibrium point to be equal to the setpoint. This optimization problem is feasible for any $x$ in a polyhedral
	region denoted as $\setX_N^r(y_{sp})$.

	Under certain assumptions \cite{MayneAUT00}, for any setpoint $y_{sp}$ and any feasible initial state $x \in
	\setX_N^r(y_{sp})$, the control law $\kappa^r_N(x,y_{sp})$ steers the
	system to the desired setpoint fulfilling the constraints. However, this
	control law is suboptimal in the sense that it does not minimizes
	$V_\infty(x,y_{sp};\kappa^r_N(\cdot,y_{sp}))$. Fortunately, as stated in
	the following lemma, if the terminal cost function is the optimal
	cost of the unconstrained LQR, then the resulting finite horizon MPC
	is equal to the constrained LQR in a neighborhood of the terminal
	region \cite{HuTAC02}.\\
	
	\begin{lem}[Local Optimality]\label{chaMPCT:locopt_LemHu}
		Consider that assumptions \ref{assumption1_OPT} and
		\ref{chaMPCT:assumption_STAB_ineqcons} hold. Consider also that the terminal control gain
		$K$ is the one of the unconstrained Linear Quadratic Regulator.
		Define the set $\Upsilon_N(y_{sp}) \subset \mathbb{R}^n$ as
		\[
		\Upsilon_N(y_{sp})=\{x \in \mathbb{R}^n \mid (\phi(N;x,\kappa_\infty(\cdot,y_{sp})),x_{sp},u_{sp}) \in \Omega^a_t \}
		\]
		where $(x_{sp},u_{sp})$ is the equilibrium point associated to the setpoint $y_{sp}$.
		Then for all $x \in \Upsilon_N(y_{sp})$,
		$V_N^{r,0}(x,y_{sp})=V^0_\infty(x,y_{sp})$ and
		$\kappa_N^r(x,y_{sp})=\kappa_\infty(x,y_{sp})$.\\
	\end{lem}
	This lemma directly stems from \cite[Thm. 2]{HuTAC02} and it states that the MPC for regulation designed using the LQR  a terminal control law ensures local optimality in a neighborhood of the setpoint and besides this region is enlarged as long as the prediction horizon is enlarged.
	
	The MPC for tracking presented here, might not ensure this local optimality
	property under the assumptions of Lemma \ref{chaMPCT:locopt_LemHu} due to the
	artificial steady state and input, and the cost function to
	minimize. However, as it is proved in the following lemma,
	under Assumption \ref{chaMPCT:assum_locopt} on the offset cost function $V_O(\cdot,\cdot)$, this
	property holds.\\
	
	\begin{assumpt}\label{chaMPCT:assum_locopt}
		Let the offset cost function $V_O(\cdot,\cdot)$  defined by Assumption \ref{chaMPCT:def_optimo}, be such that
		\[
		V_O(y,y_{sp}) \geq \gamma \|y-y_{sp}\| , \quad \forall y \in
		\setY_{sp}
		\]
		where $\gamma$ is a positive real constant.
	\end{assumpt}
	
	\begin{propt}\label{chaMPCT:OPT:ProptOpt}
		\hfill\\
		Consider that assumptions \ref{assumption1_OPT},
		\ref{chaMPCT:assumption_STAB_ineqcons} and \ref{chaMPCT:assum_locopt} hold. Then there exists
		a $\gamma^*>0$ such that for all $\gamma\geq \gamma^*$:
		\begin{itemize}
			\item The proposed MPC for tracking
			provides the same control law as the MPC for regulation, that is
			$\kappa_N(x,y_{sp})=\kappa^r_N(x,y_{sp})$ and
			$V_N^{0}(x,y_{sp})=V_N^{r,0}(x,y_{sp})$ for all $x \in \setX_N^r(y_{sp})$.
			
			\item If the terminal control gain
			$K$ is the one of the unconstrained Linear Quadratic Regulator, then
			the MPC for tracking control law $\kappa_N(x,y_{sp})$ is equal to the
			optimal control law
			$\kappa_\infty(x,y_{sp})$ for all $x \in \Upsilon(y_{sp})$.\\
		\end{itemize}
	\end{propt}

	\begin{pf}
		First, define the following optimization problem:
		\begin{subequations}\label{chaMPCT:prob_locopt}
			\beqna
			V_{N,\gamma}^{0}(x,y_{sp},\gamma)\!\!&\!\!=\!\!&\!\!\min\limits_{\vu,x_a,u_a}\sum\limits_{j=0}^{N-1}\!\|x(j)\!-\!x_a\|^2_Q\!+\!\|u(j)\!-\!u_a\|^2_R \nonumber\\
			&& +\|x(N)\!-\!x_a\|^2_P+V_O(y_a,y_{sp})\\
			&s.t.& x(0)=x, \\
			&& x(j+1)=A x(j)+B u(j), \\
			&& (x(j),u(j)) \in \setZ, \!\quad\!  j=0,\cdots, N\!\!-\!\!1 \\
			&& y_a=Cx_a+Du_a,\\
			&& (x(N),x_a,u_a) \in \Omega^a_t
			\eeqna
		\end{subequations}
		Considering $V_O(y_a,y_{sp})=\gamma\|y_a-y_{sp}\|_1$, with $\|.\|_1$ dual of $\|.\|_\infty$\footnote{The dual $\|.\|_p$ of a given norm $\|.\|_q$ is defined
			as $\|u\|_p\triangleq \displaystyle\max_{\|v\|_q\leq 1}u'v$. For
			instance, $p=1$ if $q=\infty$ and vice versa, or $p=2$ if $q=2$
			\cite{LuenbergerLIB84}.}, the optimization problem \ref{chaMPCT:prob_locopt} results
		from optimization problem \ref{chaMPCT:prob_regulation} with the last
		constraint \eqref{chaMPCT:prob_regulation_eq_constr} posed as an exact penalty function
		\cite{LuenbergerLIB84}.
		
		Therefore, there exists a finite constant
		$\gamma^*>0$ such that for all $\gamma \geq \gamma^*$,
		$V_{N,\gamma}^{0}(x,y_{sp})=V_N^{r,0}(x,y_{sp})$ for all $x \in
		\setX_N^r(y_{sp})$ \cite{LuenbergerLIB84}. Hence $V_N^{0}(x,y_{sp})=V_N^{r,0}(x,y_{sp})$.
		
		The second claim is derived from Lemma \ref{chaMPCT:locopt_LemHu} observing that
		$\Upsilon_N(y_{sp}) \subseteq \setX_N^r(y_{sp})$.
	\end{pf}

	It is worth to notice that, in virtue of the well-known result on the exact penalty functions
	\cite{LuenbergerLIB84}, the constant $\gamma$ can be chosen such
	that $\|\nu(x,y_{sp})\|_1\leq\gamma$, where $\nu(x,y_{sp})$ is the
	Lagrange multiplier of the equality constraint
	\eqref{chaMPCT:prob_regulation_eq_constr} of the optimization problem \ref{chaMPCT:prob_regulation}.
	Since the optimization problem depends on the parameters $(x,y_{sp})$,
	the value of this Lagrange multiplier also
	depends on $(x,y_{sp})$.

	Moreover, notice that the local optimality property can be ensured using any norm (not just 1-norms and $\infty$-norms), thanks
	to the properties of the duality of the norms and of equivalence of the norms, that is $\exists c>0$
	such that $\|x\|_q\geq c\|x\|_1$. On the other hand, the square of a norm
	cannot be used. With the $\|.\|_q^2$ norm, in fact, there will be
	always a local optimality gap for a finite value of $\gamma$ since
	$\|.\|_q^2$ is a (not exact) penalty function,
	\cite{LuenbergerLIB84}. That gap can be reduced by means of a
	suitable penalization of the offset cost function,
	\cite{AlvaradoPHD07}.\\
	
	%
	%\begin{rem}\label{remark1}
	%Assumption \ref{assumption3} can be easily satisfied for any
	%function $\hat V_O(.)$ considering as offset cost function $V_O(y)=
	%\max (\hat V_O(y),\alpha \|y\|)$ which is a convex function. If
	%$\setY_s$ is bounded, the upper bound condition is directly
	%fulfilled.
	%\end{rem}

	\begin{exmp}\label{chaMPCT:ejemplo_4}
		
		\begin{figure}[!h]
			\centering
			\includegraphics[width=0.99\textwidth]{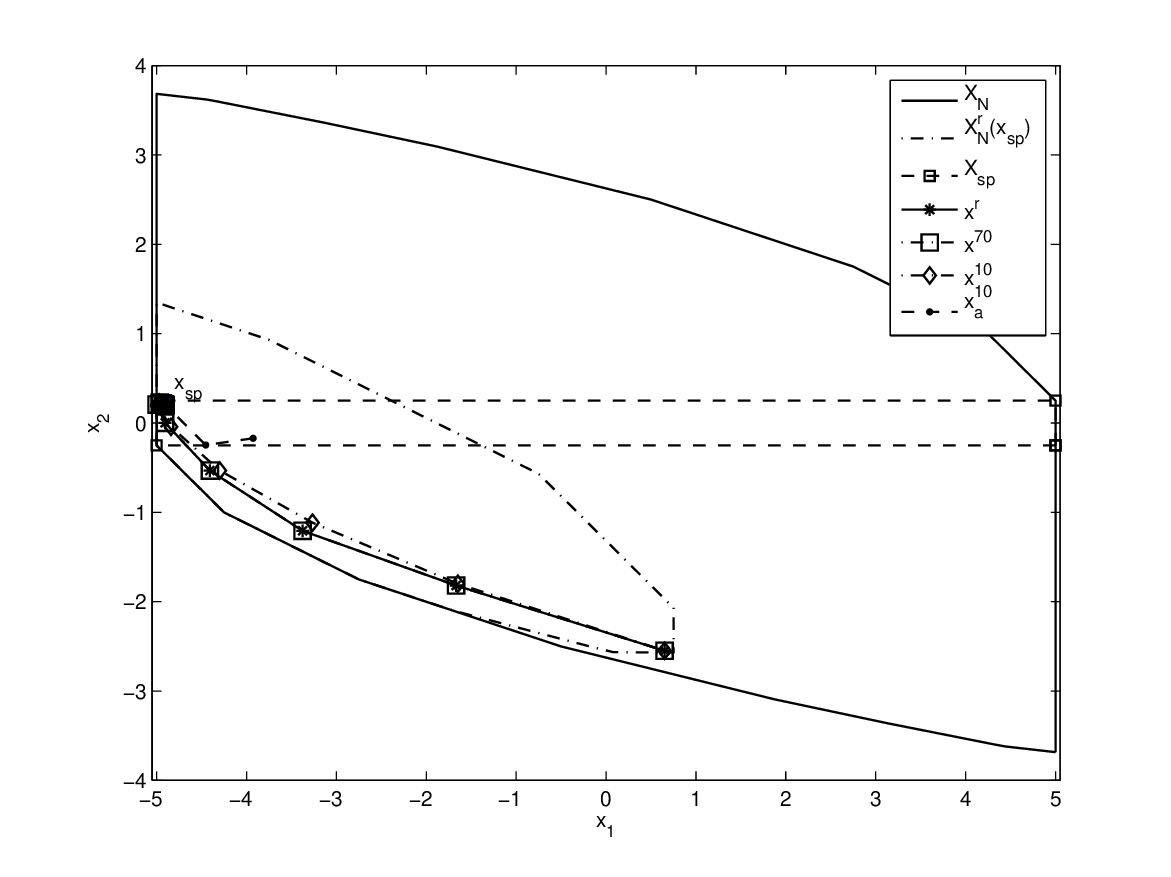}
			\caption{State space evolution of MPCT with $\gamma=10$ (red dashed line), $\gamma=70$ (red solid line) and MPC for regulation (black dotted line).}\label{chaMPCT:di_locopt_track}
		\end{figure}
		
		\begin{figure}[!h]
			\centering
			\includegraphics[width=0.99\textwidth]{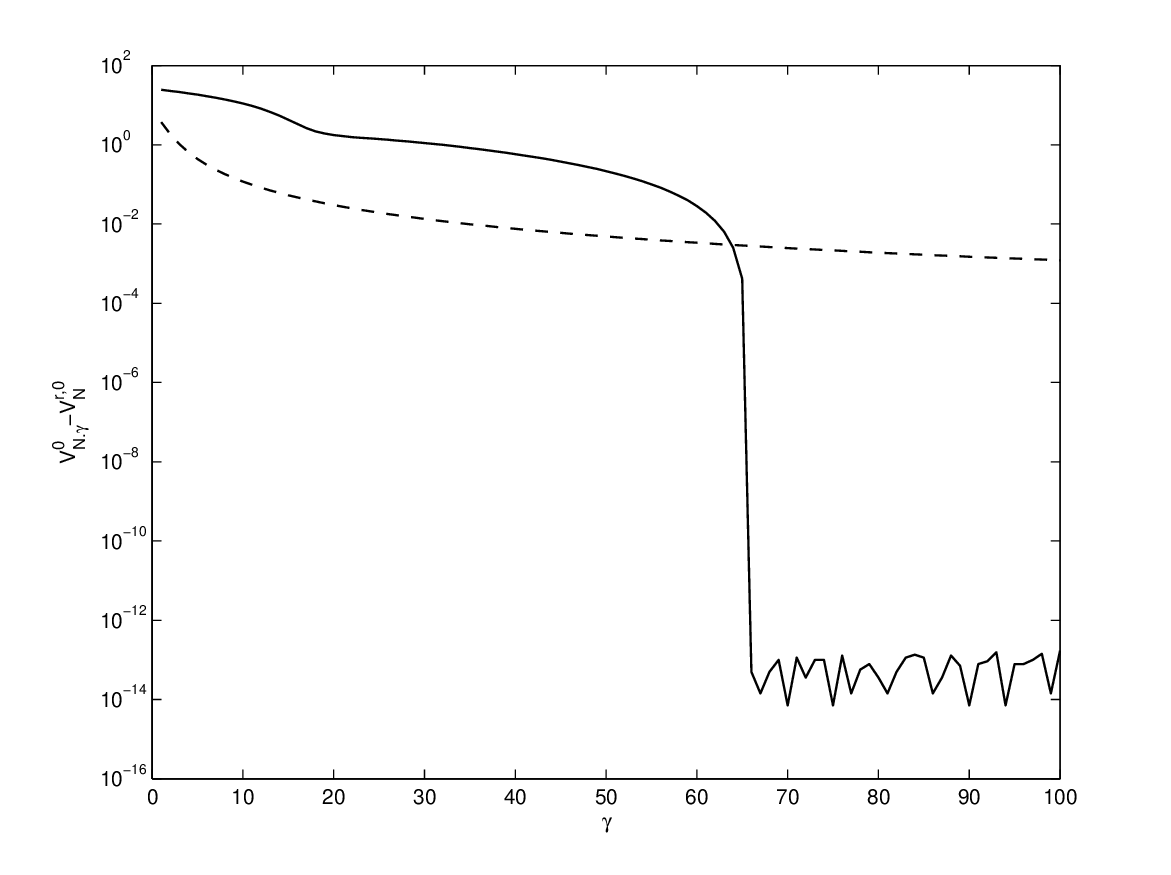}
			\caption{Difference between the optimal regulation cost and the optimal tracking
				cost versus $\gamma$.}\label{chaMPCT:di_locopt_coste}
		\end{figure}
		
		Let us study the Local Optimality property in the case study of Example \ref{chaMPCT:ejemplo_2}. In this simulations we take as initial condition a point $x(0)$ inside $\setX_N^r(x_{sp})$, which is $x(0)=(0.65,-2.55)$, with $x_{sp}=(-4.9,0.2)$ as in the previous examples.
		
		First of all, the optimal trajectory of the MPCT with offset cost function $V_O=\gamma\|y_a-y_{sp}\|_1$ is compared with the optimal trajectory of the MPC for regulation, for different values of $\gamma$. In Figure \ref{chaMPCT:di_locopt_track} the state space evolution of an MPCT with $\gamma=10$ (dash-dotted line with diamond markers), $\gamma=70$ (solid line with square markers) and MPC for regulation (solid line with star markers) are presented. Notice how the evolution of the MPCT with $\gamma=70$ is exactly the same as the one of the MPC for regulation, which is the optimal one. This is because the value of $\gamma$ is greater than the value
		of the norm of the Lagrange multiplier of the equality constraint of the
		regulation Problem \eqref{chaMPCT:prob_regulation}, which is $\|\nu(x,y_{sp})\|_1=65.69$. The evolution of the MPCT with $\gamma=10$ has a different, and suboptimal, trajectory. Notice also that, when $\gamma=70$ the value of the artificial reference is exactly $x_a=x_{sp}$, while in case of $\gamma=10$, the artificial reference takes different values, before converging to $x_{sp}$ (dashed line with dot markers).

		In Figure \ref{chaMPCT:di_locopt_coste}, the difference of the optimal cost value of the MPCT, $V_{N,\gamma}^0$ to the one of the MPC for regulation,
		$V_0^{r,0}$ has been compared for a varying values of $\gamma$. In particular, the black solid line refers to an MPCT with $V_O=\gamma\|y_a-y_{sp}\|_1$, while the black dashed line refers to an MPC with $V_O=\gamma^2\|y_s-y_t\|^2_{T}$, with $T=100I_2$.
		As it can be
		seen, in case of a quadratic cost function, the difference between the optimal costs tends to zero asymptotically, while in case of the $1$-norm, this difference drops to (practically) zero when the value of $\gamma$ becomes grater than $\|\nu(x,y_{sp})\|_1=65.69$. This result shows that the optimality gap can be made
		arbitrarily small by means of a suitable penalization of the square
		of the 2 norm, and this value asymptotically converge to zero
		\cite{AlvaradoPHD07}, while in the case of the $1$-norm, the
		difference between the optimal value of the MPC for tracking cost
		function and the standard MPC for regulation cost function becomes
		zero. \QED
	\end{exmp}

	%%%%%%%%%%%%%%%%%%%%%%%%%%%%%%%%%%%%%%%%%%%%%%%%%%%%%%%%%%%%%%%%%%%

	\subsection{Characterization of the region of local optimality}
	
	Some questions arise from the above result, such as (i) how  a suitable value of
	the parameter $\gamma$ can be determined for all possible set of
	parameters; (ii) if there exists a region where local
	optimality property holds for a given value of $\gamma$. These issues are analyzed in \cite{FerramoscaIJSS11}.
	
	In such a work it is shown that, in order to characterize the region where the Local Optimality Property holds, it is interesting to study the region where the norm of
	the Lagrange multiplier $\nu(x,y_{sp})$ is lower than or equal to
	$\gamma$. The characterization of this
	region can be done resorting to well known results on the multiparametric quadratic programming problems
	\cite{BemporadAUT02,jones2006primal,morari2008multiparametric}.
	
	To this aim, notice that  Problem \eqref{chaMPCT:prob_regulation}
	is a multiparametric problem, and it can be posed in the form of Problem \eqref{chaMPCT:QP_mp}.
	
	It is shown in \cite{FerramoscaIJSS11} that, by solving the the Karush-Kuhn-Tucker (KKT) optimality conditions \cite{BoydLIB06} of Problem \eqref{chaMPCT:prob_regulation}, the maximum  and the minimum value of
	$\|\nu(x,y_{sp})\|_1$ for all possible values of $(x,y_{sp})$ can be
	computed, that is, the values of $\gamma_{min}$ and $\gamma_{max}$
	such that for all $(x,y_{sp})\in \Gamma$,
	$\gamma_{min}\leq\|\nu(x,y_{sp})\|_1\leq\gamma_{max}$. Notice that, since Problem \eqref{chaMPCT:prob_regulation}
	is such that the solution of its KKT conditions is unique, then the
	value of $\gamma_{max}$ is finite.
	
	The set of parameters $(x,y_{sp})$,
	$\Gamma(\gamma)$, such that the norm of $\nu(x,y_{sp})$ is bounded by $\gamma$, is:
	\[
	\Gamma(\gamma)=\{(x,y_{sp}) \mid \exists (\lambda,\nu) \mbox{ and } \|\nu\|_1\leq\gamma \}
	\]
	where $\lambda$ and $\nu$ are respectively the Lagrange multipliers of the inequality and equality constraint of Problem \eqref{chaMPCT:prob_regulation}. Given this set, we can now characterize the region where the property of Local Optimality holds.
	\begin{lem}\label{chaMPCT:Lem2Opt}
		\hfill\\
		Consider that Lemma \ref{chaMPCT:OPT:ProptOpt} holds. Then:
		\begin{itemize}
			\item For all $\gamma>\gamma_{min}$, there exists a polygon $\Gamma(\gamma)$
			such that if $(x,y_{sp})\in\Gamma(\gamma)$, then
			$V_N^r(x,y_{sp})=V_N(x,y_{sp})$.
			\item For all $\gamma_{min}<\gamma_{a}\leq\gamma_b$, $\Gamma(\gamma_a)
			\subseteq\Gamma(\gamma_b)$. That is, $\Gamma(\gamma)$ grows
			monotonically with $\gamma$.
			\item For all $\gamma\geq\gamma_{max}$, $\Gamma(\gamma)=\Gamma=\{(x,y_{sp}) \mid x \in \setX_N^r(y_{sp}) \}$.
		\end{itemize}
	\end{lem}

	\begin{theorem}[Region of local optimality]\label{chaMPCT:LocalOpt}
		\hfill\\
		Consider that Lemma \ref{chaMPCT:OPT:ProptOpt} and Lemma \ref{chaMPCT:Lem2Opt}
		hold. Define the following region
		$$\setW(\gamma,y_{sp})=\{x\in\Upsilon_N(y_{sp}) \mid
		(\phi(i;x,\kappa_N(\cdot,y_{sp})),y_{sp})\in\Gamma(\gamma), \forall i\geq
		0\}$$ and let the terminal control gain $K$ be the one of the
		unconstrained LQR. Then:
		\begin{enumerate}
			\item For all $\gamma>\gamma_{min}$, $\setW(\gamma,y_{sp})$ is a non-empty polygon and it is
			a positively invariant set of the controlled system.
			\item If $\gamma_{min}<\gamma_a\leq\gamma_b$, then
			$\setW(\gamma_a,y_{sp})\subseteq\setW(\gamma_b,y_{sp})$.
			\item If $\gamma>\gamma_{min}$,  $x(0)$ and $y_{sp}$ are such that
			$x(0)\in\setX_N^r(y_{sp})$, then
			\begin{enumerate}
				\item There exists an instant $\bar k$ such that  $x(\bar
				k)\in\setW(\gamma,y_{sp})$ and
				$\kappa_N(x(k),y_{sp})=\kappa_{\infty}(x(k),y_{sp})$, for all $k\geq\bar
				k$.
				\item If $\gamma\geq\gamma_{max}$ then $\kappa_N(x(k),y_{sp})=\kappa_N^r(x(k),y_{sp})$ for all $k\geq
				0$ and there exists an instant $\bar k$ such that $x(\bar k) \in
				\Upsilon_N(y_{sp})$ and $\kappa_N(x(k),y_{sp})=\kappa_{\infty}(x(k),y_{sp})$
				for all $k\geq \bar k$.
			\end{enumerate}
		\end{enumerate}
	\end{theorem}

	What the last theorem infers is that, for every
	$\gamma\geq\gamma_{min}$, the MPC for tracking is locally optimal in
	a certain region. In particular the value of $\gamma_{min}$ is
	interesting from a theoretical point of view, because it is the
	critical value that ensures the existence of a region of Local Optimality. Moreover, Theorem \eqref{chaMPCT:LocalOpt} also shows that the region of Local Optimality grows monotonically with $\gamma$. The maximal region of Local Optimality is given for any $\gamma\geq\gamma_{max}$, and it is equal to the feasible set of Problem \eqref{chaMPCT:prob_regulation}.

	% % % % % % % % % % % % % % % % % % % % % % % % % % % % % % % % % % % % % % % % % % % % % % % % % % % %
	
	\section{Application to the a four tanks plant}
	In this section, the properties of the controller presented in this
	report, are proved in an application to the four
	tanks plant, inspired by the educational quadruple-tank process proposed in \cite{JohanssonTCST00}.
	
	\subsection{Tracking reachable and unreachable setpoints}
	The aim of the first test is to show the property of offset
	minimization of the controller. The offset cost function has been
	chosen as $V_O=\alpha\|y_a-y_{sp}\|_{\infty}$. In the test, five
	references have been considered:$y_{sp,1}=(0.3,0.3)$,
	$y_{sp,2}=(1.25,1.25)$, $y_{sp,3}=(0.35,0.8)$, $y_{sp,4}=(1,0.8)$ and
	$y_{sp,5}=(h^0_{1},h^0_{2})$. Notice that $y_{sp,3}$ is not reachable. The initial state is
	$x_0=(0.65,0.65,0.6658,0.6242)$. An MPC with terminal inequality constraint and $N=3$ has been
	considered. The weighting matrices have been chosen as $Q=I_4$ and
	$R=0.01\times I_2$. Matrix $P$ is the solution of the Riccati
	equation and $\alpha=50$.
	
	\begin{figure}[h]
		\centering
		\includegraphics[width=0.98\textwidth]{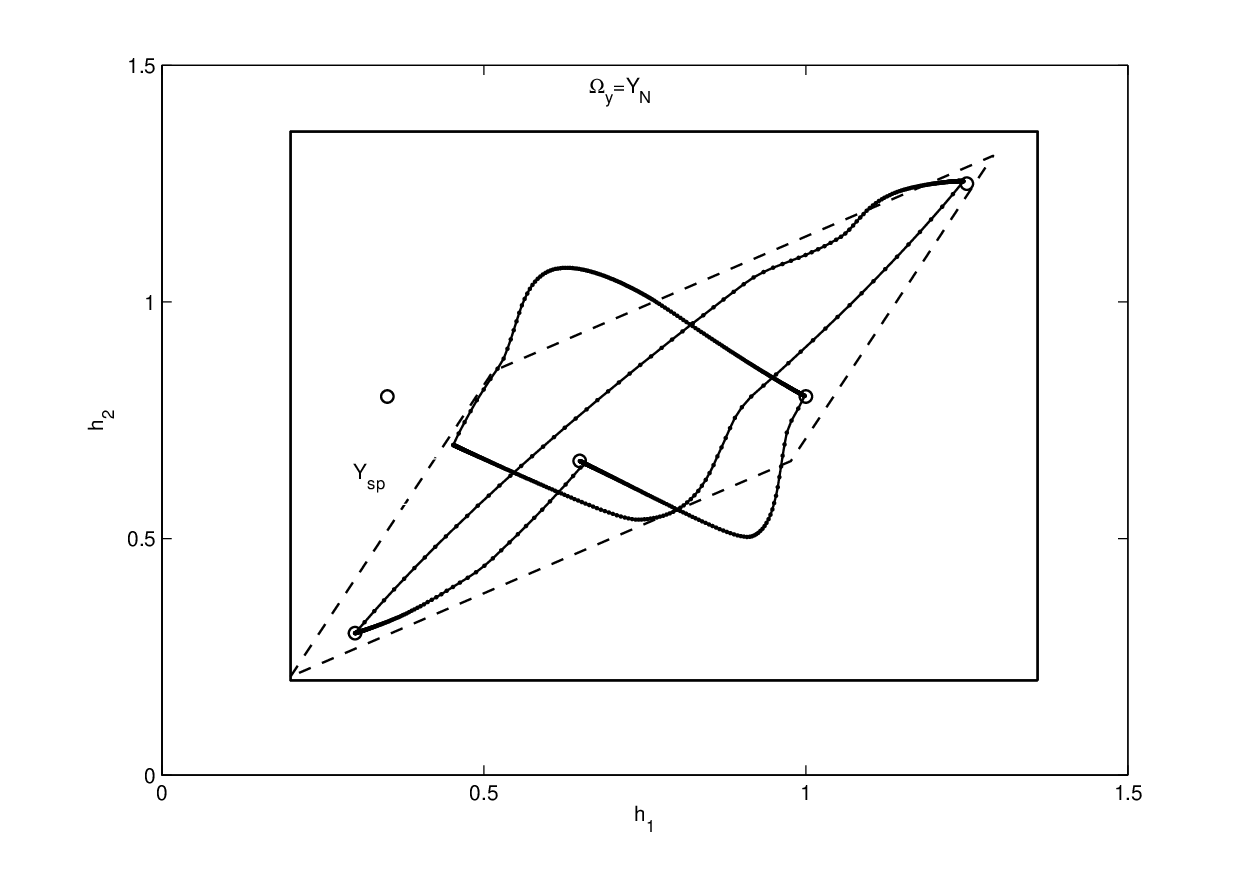}
		\caption{Application to the four tanks plant. Output-space evolution of the closed-loop system.}\label{cha1:4T_sets}
	\end{figure}
	
	\begin{figure}[h]
		\centering
		\includegraphics[width=0.98\textwidth]{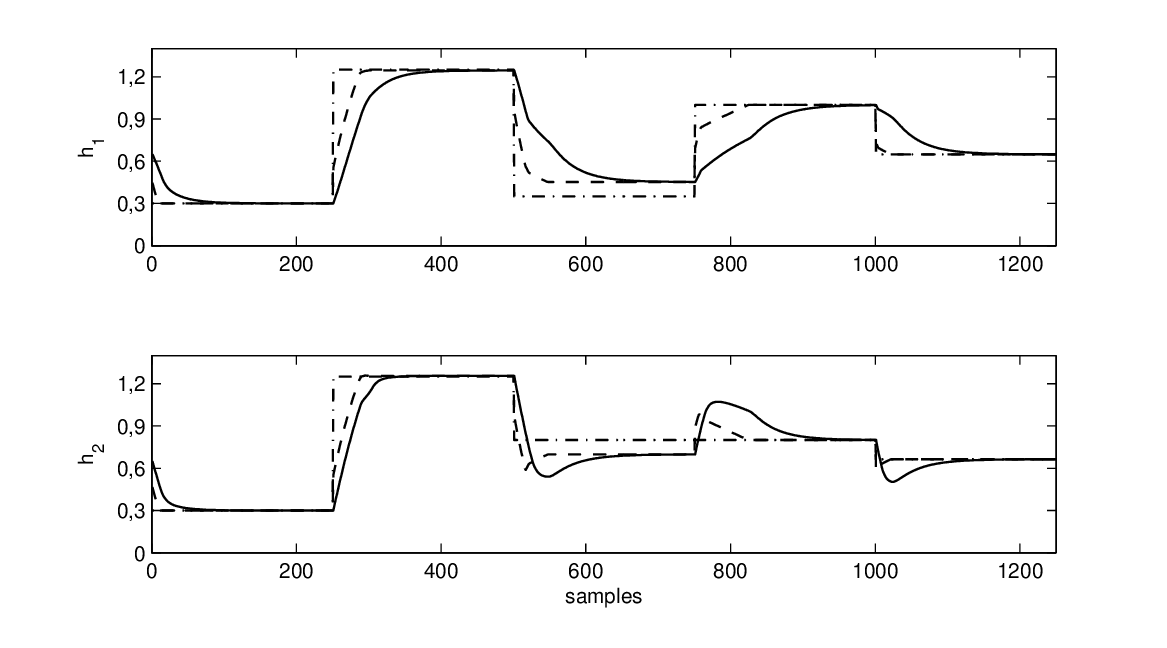}
		\caption{Application to the four tanks plant. Time evolution of the levels $h_1$ and
			$h_2$: reference in dashed-dotted line, artificial reference in dashed line, real evolution of the system in solid line.}\label{cha1:4T_h1h2}
	\end{figure}
	
	\begin{figure}
		\centering
		\includegraphics[width=0.98\textwidth]{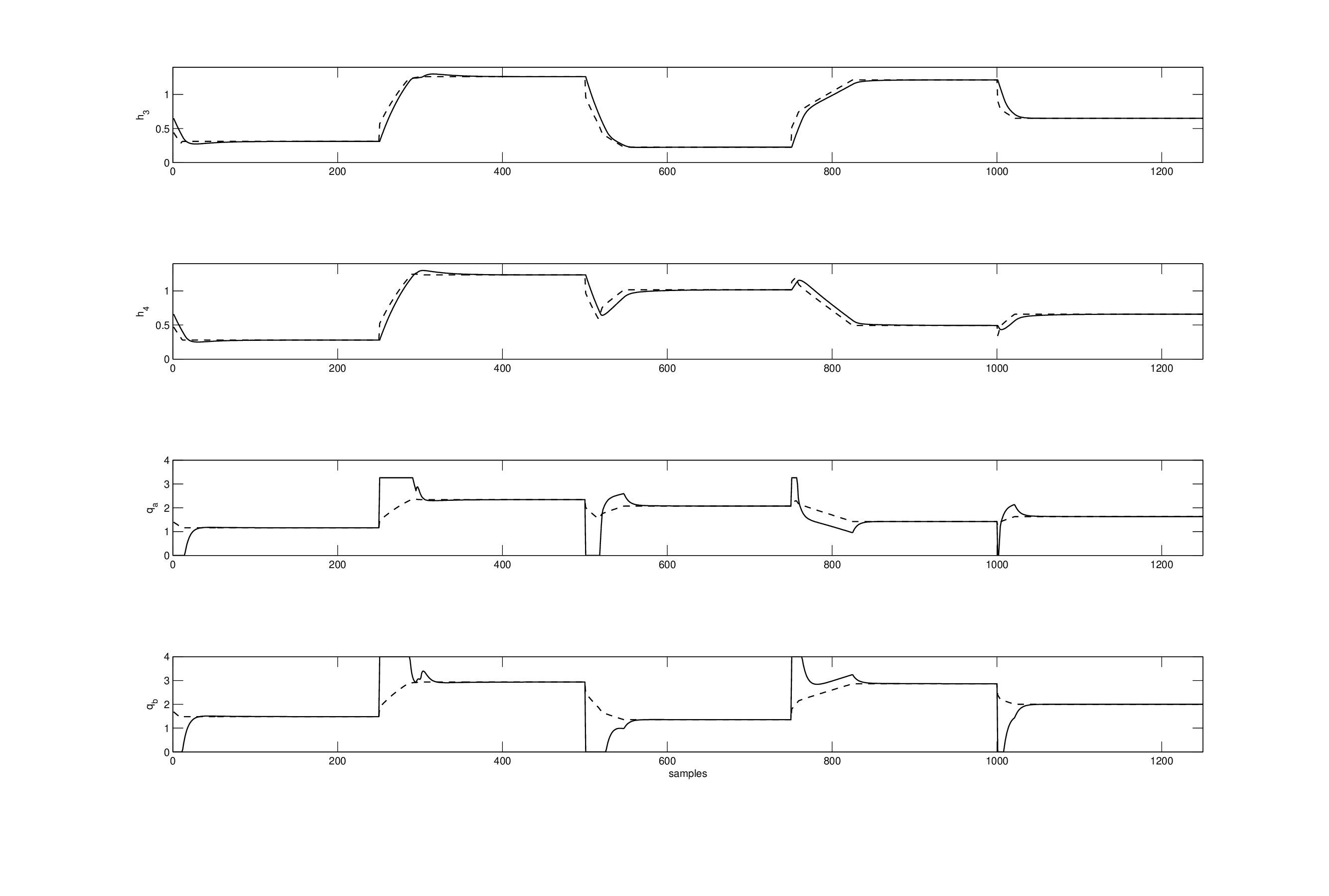}
		\caption{Application to the four tanks plant. Time evolution of the levels $h_3$ and $h_4$ and flows $q_a$ and $q_b$: artificial reference in dashed line, real evolution of the system in solid line.}\label{cha1:4T_h3h4qaqb}
	\end{figure}

	The projection of the maximal invariant set for tracking onto $y$,
	$\Omega_{y}$, the projection of the region of attraction onto $y$,
	$\setY_3$, the set of equilibrium levels $\setY_s$ and the
	ouput-space evolution of the levels $h_1$ and $h_2$ are shown in
	Figure \ref{cha1:4T_sets}. The time evolutions are shown in Figures
	\ref{cha1:4T_h1h2} and \ref{cha1:4T_h3h4qaqb}. The reference is
	depicted in dashed-dotted line, while the artificial reference and
	the real evolution of the system are depicted respectively in dashed
	and solid line. 
	
	Since the initial state is an admissible equilibrium 
	point the controller is feasible. The control law steers the system to the setpoints irrespective of the size of the steps. As it can be seen, when the desired setpoint is
	reachable, the closed-loop system converges to it  without any offset.
	When the reference changes to an unreachable setpoint, the
	controller drives the system to the closest equilibrium point, in the
	sense that the offset cost function is minimized.

	\subsection{Local Optimality}
	
\begin{figure}
		\begin{minipage}{0.58\textwidth}
			\begin{center}
				\centering
				\includegraphics[width=0.98\textwidth]{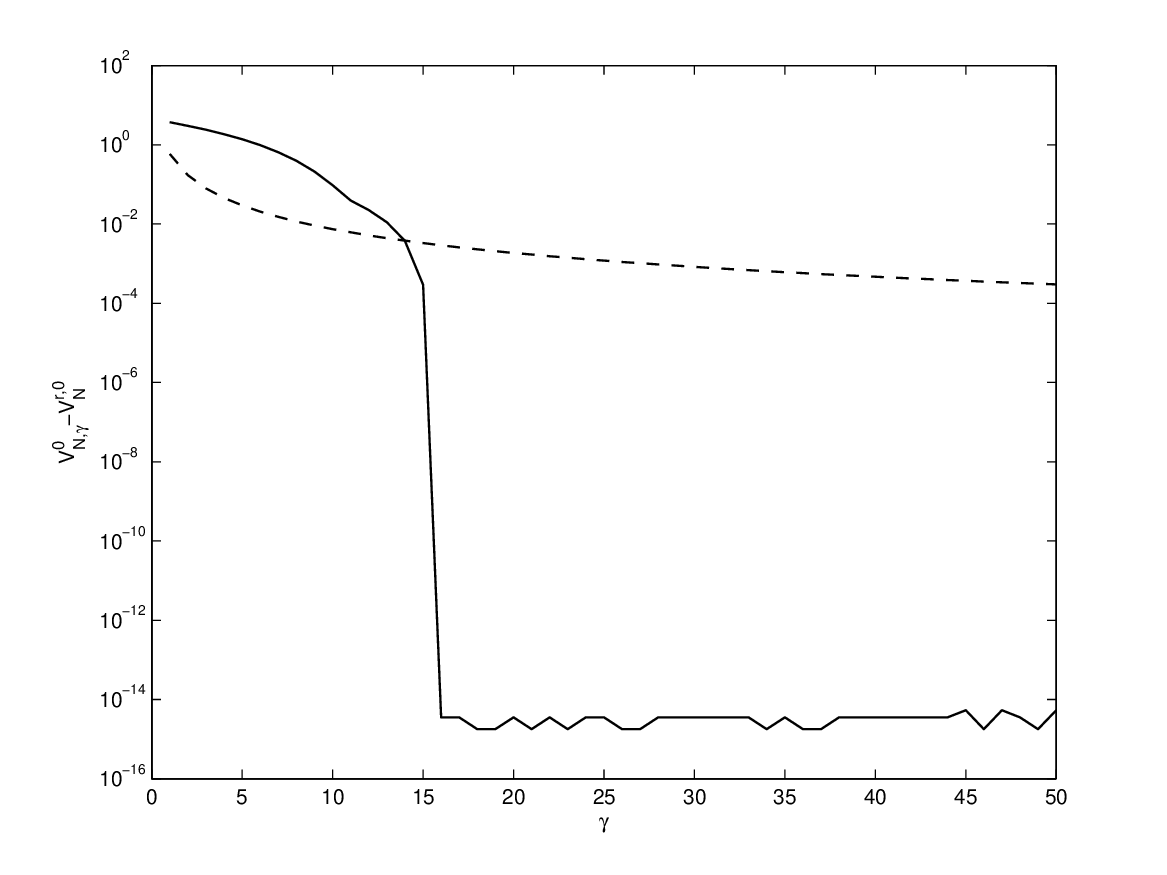}
			\end{center}
		\end{minipage}
		\begin{minipage}{0.42\textwidth}
			\begin{center}
				%\label{differentALFA}
				\begin{scriptsize}
					\begin{tabular}{|c|c|c|}
						\hline
						$\gamma$ & $V_{N,\gamma}N^0\!-\!V^{r,0}_N$\\
						\hline \hline
						$14$ & $ 0.0038$ \\
						$15$ & $ 2.93e\!-\!4$\\
						$15.3$ & $ 1.47e\!-\!5$\\
						$15.38$ & $ 1.47e\!-\!5$\\
						$15.386$ & $ 1.15e\!-\!9$\\
						$\gamma_{max}$ & $3.55e\!-\!15$\\
						$16$ & $3.55e-15$\\
						\hline
					\end{tabular}
				\end{scriptsize}
			\end{center}
		\end{minipage}
		\caption{Application to the four tanks plant. Difference between the optimal regulation cost and the optimal tracking
			cost versus $\gamma$.}\label{cha1:confronto_costi}
	\end{figure}

	To illustrate the property of local optimality, we compare  an MPC for tracking with linear offset cost function, with an MPC for tracking with quadratic offset cost function. To this aim, the
	difference of the optimal cost value provided by these two controllers, and the one of the MPC for regulation,
	$V_N^{r,0}$ has been compared. The quadratic offset
	cost function has been chosen as $V_O=\gamma^2\|y_a-y_{sp}\|^2_{T}$ with
	$T=100 I_4$. The linear offset cost function
	has been taken as a $1$-norm, $V_O=\gamma\|y_a-y_{sp}\|_{1}$.  The
	system has been considered to be steered to the point
	$y_{sp}=(h^0_{1},h^0_{2})$, with initial condition $y_0=(1.25,1.25)$. In
	Figure \ref{cha1:confronto_costi} the value $V_{N,\gamma}^{0}-V_N^{r,0}$
	versus $\gamma$ is plotted in solid line for the linear offset const function, and in dashed line for the quadratic one. As it can be
	seen, the dashed line tends to zero asymptotically while
	the solid line drops to (practically) zero when $\gamma=16$.\\
	%This result shows that the optimality gap can be made
	%arbitrarily small by means of a suitable penalization of the square
	%of a $2$-norm, and this value asymptotically converge to zero
	%\cite{AlvaradoPHD07}, while in case of a $1$-norm, the
	%difference between the optimal value of the MPCT cost
	%function and the standard MPC for regulation cost function becomes
	%zero.\\
	%Note how the value of $V_{N,\gamma}^{0}-V_N^{r,0}$ drops to practically zero
	%when $\gamma=16$. 
	This happens because the value of $\gamma$ becomes greater than the value
	of the Lagrange multiplier of the equality constraint of the
	regulation Problem \eqref{chaMPCT:prob_regulation}. In this test, the
	equality constraint of Problem \eqref{chaMPCT:prob_regulation} has been
	chosen as an $\infty$-norm, and hence, to obtain an exact penalty
	function, the offset cost function of Problem \eqref{chaMPCT:prob_locopt} has
	been chosen as a $1$-norm. To point out this fact, consider that,
	for this example, the maximum value of the Lagrange multipliers of
	the equality constraint Problem \eqref{chaMPCT:prob_regulation} is $\gamma_{max}=15.3868$. In the
	table of Figure \ref{cha1:confronto_costi}, the value of $V_{N,\gamma}^{0}-V_N^{r,0}$in case of different values
	of the parameter $\gamma$ is presented. Note how the value seriously
	decreases when $\gamma$ becomes equals to $\gamma_{max}$.
	
\begin{figure}
		\begin{minipage}{0.58\textwidth}
			\begin{center}
				\centering
				\includegraphics[width=0.98\textwidth]{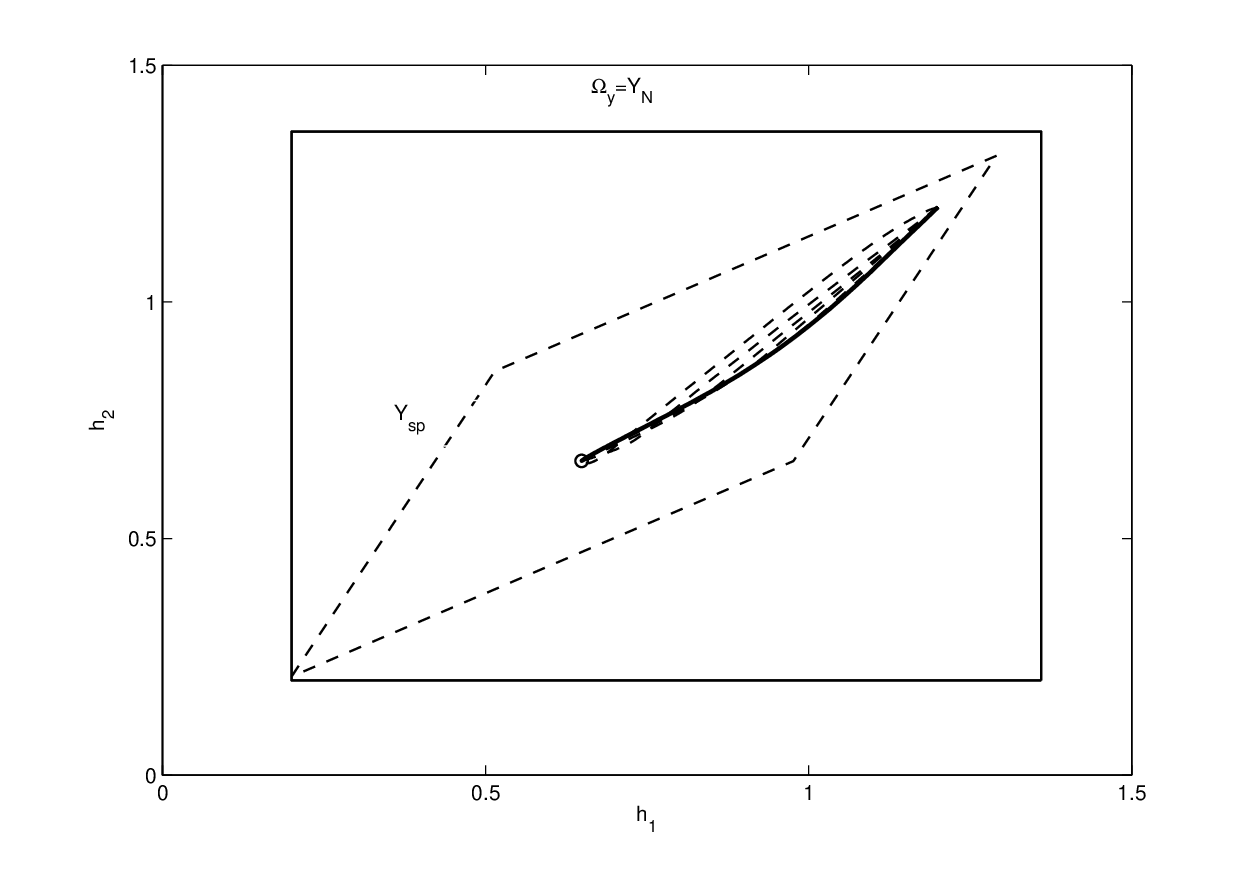}
			\end{center}
		\end{minipage}
		\begin{minipage}{0.42\textwidth}
			\begin{center}
				%\label{differentALFA}
				\begin{scriptsize}
					\begin{tabular}{|c|c|c|}
						\hline
						$\gamma$ & $V_{N,\gamma}^0$\\
						\hline \hline
						$2$ & $84.2693$ \\
						$4$ & $40.6147$\\
						$6$ & $22.4751$\\
						$8$ & $13.8412$\\
						$10$ & $10.0638$\\
						$12$ & $8.9572$\\
						$14$ & $8.6989$\\
						$15.3868$ & $8.6084$\\
						$18$ & $8.6084$\\
						$20$ & $8.6084$\\
						\hline
					\end{tabular}
				\end{scriptsize}
			\end{center}
		\end{minipage}
		\caption{Application to the four tanks plant. State-space optimal trajectories and optimal cost for $\gamma$
			varying.}\label{cha1:opt_traj}
	\end{figure}

	Last, the optimal trajectories from the point
	$y_0=(1.25,1.25)$ to the point $y=(h^0_{1},h^0_{2})$ have been
	calculated, for a value of $\gamma$ that varies in the set
	$$\gamma=\{2,4,6,8,10,12,14,\gamma_{max},18,20\}$$ In figure
	\ref{cha1:opt_traj} the state-space trajectories and the values of
	the optimal cost $V_{N,\gamma}^0$ for $\gamma$ increasing are shown. See how the value of the optimal
	cost decreases as the value of $\gamma$ increases. The optimal
	trajectory, in solid line, is the one for which
	$\gamma=\gamma_{max}$. Notice that value of the optimal cost
	decreases from $V_{N,2}^0=84.2693$ to $V_{N,\gamma_{max}}^0=8.6084$ when
	$\gamma$ reaches the value of $\gamma_{max}$.

	%%%%%%%%%%%%%%%%%%%%%%%%%%%%%%%%%%%%%%%%%%%%%%%%%%%%%%%%%%%%%%%%%
	\clearpage
	\section{Conclusions}
	In this report the MPC for tracking formulation has been
	presented.
	This formulation is based on four main ingredients:
	\begin{enumerate}
		\item[(i)] an artificial steady state and input, considered as decision
		variables;
		
		\item[(ii) ] a stage cost that penalizes the deviation of the predicted
		trajectory from the artificial steady conditions;
		
		\item[(iii) ] an extra cost, the \emph{offset cost} function, added to penalize the
		deviation of the artificial steady state from the target setpoint;
		
		\item[(iv) ] a relaxed terminal constraint.
	\end{enumerate}
	
	It has been shown that the proposed controller ensures recursive feasibility, and convergence to the desired setpoint (or to the best admissible equilibrium point), for the case of both terminal equality constraint and terminal inequality constraint. Asymptotic stability has been proved, providing a Lyapunov function.
	
	It has been also shown how to cast the MPCT problem as a Quadratic Programming problem, and how to calculate an explicit off-line control law, resorting to well known Multi Parametric Programming tools.
	
	Furthermore, it has been proved that, under some mild assumptions on the offset cost function, the MPCT guarantees the Local Optimality property, that is the optimal value of the MPCT cost function is equal to the one of an MPC for regulation, and that the region, where this property is ensured, is the MPC for regulation feasible set.

%%%%%%%%%%%%%%%%%%%%%%%%%%%%%%%%%%%%%%%%%%%%%%%%%%%%%%%%%%%%%%%%%%%%%%%%%%%%%%%

\section{Appendix: technical lemmas}

\begin{lem}\label{chaMPCT:Lem_Func_DefPos}
	Consider that Assumptions \ref{assumption1_OPT}-\ref{chaMPCT:assumptionQR_eqconst} hold. Let $x_{t}$ be the optimal steady state, such that function $V_O(y,y_{sp})$ is minimized. For all $x \in \setX_{N}$ and $x_a^0(x) \in \setX_{sp}$, define the function $e(x)=x-x_a^0(x)$. Then, there exists a $\mathcal{K}$-function $\alpha_e$ such that
	\beqn
	\|e(x)\|\geq \alpha_e(\|x-x_{t}\|)
	\eeqn
\end{lem}

\begin{pf}
	Notice that, due to convexity, $e(x)$ is a continuous function \cite{RawlingsLIB09}. Moreover, let us consider these two cases.
	\begin{enumerate}
		\item $\|e(x)\|=0$ iff $x=x_{t}$. In fact, (i) if $e(x)=0$, then $x=x_a^0(x)$, and from Lemma \ref{chaMPCT:lema_optimalidad} (or Lemma \ref{chaMPCT:lema_optimalidad_ineqconst} in case of terminal inequality constraint), this implies that $x_a^0(x)=x_{t}$; (ii) if $x=x_{t}$, then by optimality $x_a^0(x)=x_{t}$, and then $x=x_a^0(x)$. Then, $\|e(x)\|=0$.
		
		\item $\|e(x)\|>0$ for all $\|x-x_{t}\|>0$. In fact, for any $x \not = x_{t}$, $\|e(x)\|\not = 0$ and moreover $\|x-x_{t}\|>0$. Then, $\|e(x)\|>0$.
	\end{enumerate}
	Then, since $\setX_{N}$ is compact, in virtue of \cite[Ch. 5, Lemma 6, pag. 148]{VidyasagarLIB93}, there exists a $\mathcal{K}$-function $\alpha_e$ such that $\|e(x)\|\geq \alpha_e(\|x-x_{t}\|)$ on $\setX_{N}$.
\end{pf}

\begin{lem}\label{chaMPCT:lema_optimalidad}
	Consider that Assumptions \ref{assumption1_OPT}-\ref{chaMPCT:assumptionQR_eqconst} hold. Let the optimal solution to Problem \eqref{chaMPCT:optprob_eqconst}, at time $k$, be such that $x(k)=x_a^0(x(k))$, $u(k)=u_a^0(x(k))$, and $y(k)=y_a^0(x(k))$, and $x(k+1)=x_a^0(x(k))$. Let $(x_{t},u_{t},y_{t})$ be the optimal triplet, such that function $V_O(y,y_{sp})$ is minimized. Then $x(k)=x_{t}$, $u(k)=u_{t}$, and $y(k)=y_{t}$.
\end{lem}

\begin{pf}
	Consider that $(x_a^0(x(k)),u_a^0(x(k)),y_a^0(x(k)))$ is the optimal solution to \eqref{chaMPCT:optprob_eqconst} at time $k$. Then $$V_{N}^0(x(k))=V_O(y_a^0(x(k)),y_{sp})$$
	In what follows, the time dependence is removed for the sake of clarity.
	
	This Lemma will be proved by contradiction. Assume that the stationary point at time $k$ is not the optimal one, that is $(x_a^0(x),u_a^0(x))\neq (x_{t},u_{t})$. Then, by convexity, there exists a $\beta \in [0,1]$ such that
	$$(\tx_a,\tu_a)=\beta(x_a^0(x),u_a^0(x))+(1-\beta)(x_{t},u_{t})$$
	characterizes a stationary point and moreover
	\beqn \label{chaMPCT:eq_coste_offset_no_optimo} V_O(\tilde y_a,y_{sp})\leq V_O(y_{t},y_{sp})
	\eeqn
	That is, since the real system is not at the optimal point $(x_{t},u_{t})$, it is more convenient to move towards $(\tx_a,\tu_a)$, than to remain in $(x_a^0(x),u_a^0(x))$. Define as $\tilde \vu=\{\tu(0),\tu(i),...,\tu(N-1)\}$ a feasible sequence, to Problem \eqref{chaMPCT:optprob_eqconst}, that drives the system from $(x_a(x)^0,u_a(x)^0)$ to $(\tx_a,\tu_a)$. This sequence is such that, the $j$-th element is given by $\tu(j)=K_{db}(\tilde x(j)-\tilde x_a)+\tu_a$, and $\tilde x(j+1)= A\tilde x(j)+ B \tu(j)$, $\tilde x(0)=x_a^0(x)$. Then, the cost to drive the system from $(x_a^0(x),u_a^0(x))$ to $(\tx_a,\tu_a)$ is given by
	\beqnan
	V_{N}(x_a^0(x),y_{sp};\tilde \vu,\tx_a,\tu_a)\!\!&=&\!\!\sum\limits_{j=0}^{N-1} \!\|\tx(j)\!-\!\tx_a\|^2_Q\!+\!\|\tu(j)\!-\!\tu_a\|^2_R\!\!+V_O(\tilde y_a,y_{sp})\\
	\!\!&=&\!\!\sum\limits_{j=0}^{N-1}\!\overbrace{\|\tx(j)\!-\!\tx_a\|^2_Q\!+\!\|K_{db}(\tx(j)\!-\!\tx_a)\|^2_R}^{\|\tx(j)\!-\!\tx_a\|^2_{Q+K_{db}'RK_{db}}}\!\!+V_O(\tilde y_a,y_{sp})\\
	\!\!&=&\!\!\|x_a^0(x)\!-\!\tx_a\|^2_{\tilde P}+V_O(\tilde y_a,y_{sp})\\
	\!\!&=&\!\!(1-\beta)^2\|x_a^0(x)\!-\!x_{t}\|^2_{\tilde P}\!+V_O(\tilde y_a,y_{sp})
	\eeqnan

	Now define $W(\beta)=(1-\beta)^2\|x_a^0(x)\!-\!x_{t}\|^2_{\tilde P}\!+\!V_O(\tilde y_a,y_{sp})$ and notice that for $\beta=1$, $W(1)=V_O(y_a^0(x),y_{sp})$. Taking the partial of this function with respect to $\beta$, and evaluating it for $\beta=1$ we obtain:
	\beqnan \left.\frac{\partial
		W}{\partial\beta}\right|_{\beta=1} =g^{0'}(y_a^0(x),y_{sp})
	\eeqnan
	where $g^{0'} \in \partial V_O(y_a^0(x),y_{sp})$, defining $\partial
	V_O(y_a^0(x),y_{sp})$ as the subdifferential of $V_O(y_a^0(x),y_{sp})$. From convexity and from \eqref{chaMPCT:eq_coste_offset_no_optimo},
	\beqnan \left.\frac{\partial
		W}{\partial\beta}\right|_{\beta=1}\!\!&=&\!\!g^{0'}(y_a^0(x),y_{sp})\\
	\!\!&\geq&\!\! V_O(y_a^0(x),y_{sp}) \!-\!V_O(\tilde y_a,y_{sp})\! >\!0
	\eeqnan
	This means that there exists a value of $\beta \in [0,1)$ such
	that $V_{N}(x_a^0(x),y_{sp};\tilde \vu,\tx_a,\tu_a)$ is smaller than the value of the cost
	$V_{N}(x_a^0(x),y_{sp};\tilde \vu,\tx_a,\tu_a)$ for $\beta=1$, which is
	$V_O(y_a^0(x),y_{sp})$. This contradicts the optimality of the solution to Problem \eqref{chaMPCT:optprob_eqconst} at time $k$, and the assumption that $(x_a^0(x),u_a^0(x))$ is a fixed point, that is the optimal solution to Problem \eqref{chaMPCT:optprob_eqconst} at time $k+1$ is still $(x_a^0(x),u_a^0(x))$. Then it has to be that $(x_a^0(x),u_a^0(x))=(x_{t},u_{t})$. Moreover, from Definition \ref{chaMPCT:def_optimo}, we can state that this point is the one that minimizes the offset cost function $V_O(y,y_{sp})$. So the Lemma is proved. 
\end{pf}

\begin{lem}\label{chaMPCT:lema_optimalidad_ineqconst}
	Consider that Assumptions \ref{assumption1_OPT} and \ref{chaMPCT:assumption_STAB_ineqcons} hold. Let the optimal solution to Problem \eqref{chaMPCT:optprob_ineqconst}, at time $k$, be such that $x(k)=x_a^0(x(k))$, $u(k)=u_a^0(x(k))$, and $y(k)=y_a^0(x(k))$, and $x(k+1)=x_a^0(x(k))$. Let $(x_{t},u_{t},y_{t})$ be the optimal triplet, such that function $V_O(y,y_{sp})$ is minimized. Then $x(k)=x_{t}$, $u(k)=u_{t}$, and $y(k)=y_{t}$.
\end{lem}

\begin{pf}
	The proof to this Lemma follows same arguments as the proof to Lemma \ref{chaMPCT:lema_optimalidad}. Only notice that, in this case, the feasible sequence $\tilde \vu=\{\tu(0),\tu(i),...,\tu(N-1)\}$ that drive the closed-loop system from $(x_a^0(x),u_a^0(x))$ to $(\tx_a,\tu_a)$, is such that $\tu(j)=K(\tilde x(j)-\tilde x_a)+\tu_a$, and $\tilde x(j+1)= A\tilde x(j)+ B \tu(j)$, $\tilde x(0)=x_a^0(x)$. Then, the cost to drive the system from $(x_a^0(x),u_a^0(x))$ to $(\tx_a,\tu_a)$ is given by
	\beqnan
	V_{N}(x_a^0(x),y_{sp};\tilde \vu,\tx_a,\tu_a)\!\!&=&\!\!\sum\limits_{j=0}^{N-1} \!\|\tx(j)\!-\!\tx_a\|^2_Q\!+\!\|\tu(j)\!-\!\tu_a\|^2_R\!+\|\tx(N)\!-\!\tx_a\|^2_P\!\\
	\!\!&&\!\!+V_O(\tilde y_a,y_{sp})\\
	\!\!&=&\!\!\sum\limits_{j=0}^{N-1}\overbrace{ \!\|\tx(j)\!-\!\tx_a\|^2_Q\!+\!\|K(\tx(j)\!-\!\tx_a)\|^2_R}^{\|\tx(j)\!-\!\tx_a\|^2_{Q+K'RK}}\!+\|\tx(N)\!-\!\tx_a\|^2_P\!\\
	\!\!&&\!\!+V_O(\tilde y_a,y_{sp})\\
	\!\!&=&\!\!\|x_a^0(x)\!-\!\tx_a\|^2_{P}+V_O(\tilde y_a,y_{sp})\\
	\!\!&=&\!\!(1-\beta)^2\|x_a^0(x)\!-\!x_{t}\|^2_{P}\!+V_O(\tilde y_a,y_{sp})
	\eeqnan

	In this case also, it can be verified that the partial of $W(\beta)=(1-\beta)^2\|x_a^0(x)\!-\!x_{t}\|^2_{P}\!+\!V_O(\tilde y_a,y_{sp})$ with respect to $\beta$, is strictly positive for $\beta=1$. This means that there exists a value of $\beta \in [0,1)$ such
	that $V_{N}(x_a^0(x),y_{sp};\tilde \vu,\tx_a,\tu_a)$ is smaller than the value of the cost
	$V_{N}(x_a^0(x),y_{sp};\tilde \vu,\tx_a,\tu_a)$ for $\beta=1$, which is
	$V_O(y_a^0(x),y_{sp})$. This contradicts the optimality of the solution to Problem \eqref{chaMPCT:optprob_ineqconst} at time $k$, and the assumption that $(x_a^0(x),u_a^0(x))$ is a fixed point, that is the optimal solution to Problem \eqref{chaMPCT:optprob_ineqconst} at time $k+1$ is still $(x_a^0(x),u_a^0(x))$. Then it has to be that $(x_a^0(x),u_a^0(x))=(x_{t},u_{t})$.
\end{pf}

% % % % % % % % % % % % % % % % % % % % % % % % % % % % % % % % % % % % % % % % % % % % % % % % % % % % %

\clearpage
\bibliographystyle{elsarticle-num}
\bibliography{MPCT_int_rep}

\begin{thebibliography}{10}
\expandafter\ifx\csname url\endcsname\relax
  \def\url#1{\texttt{#1}}\fi
\expandafter\ifx\csname urlprefix\endcsname\relax\def\urlprefix{URL }\fi
\expandafter\ifx\csname href\endcsname\relax
  \def\href#1#2{#2} \def\path#1{#1}\fi

\bibitem{PannocchiaAICHEJ03}
G.~Pannocchia, J.~B. Rawlings, Disturbance models for offset-free
  model-predictive control, AIChE Journal 49 (2003) 426--437.

\bibitem{MaederAUT10}
U.~Maeder, M.~Morari, Offset-free reference tracking with model predictive
  control, Automatica 46~(9) (2010) 1469--1476.

\bibitem{BemporadTAC97}
A.~Bemporad, A.~Casavola, E.~Mosca, Nonlinear control of constrained linear
  systems via predictive reference management., IEEE Transactions on Automatic
  Control 42 (1997) 340--349.

\bibitem{ChisciIJC03}
L.~Chisci, G.~Zappa, Dual mode predictive tracking of piecewise constant
  references for constrained linear systems, Int. J. Control 76 (2003) 61--72.

\bibitem{LimonAUT08}
D.~Limon, I.~Alvarado, T.~Alamo, E.~F. Camacho, {MPC} for tracking of
  piece-wise constant references for constrained linear systems, Automatica
  44~(9) (2008) 2382--2387.

\bibitem{RossiterCTA96}
J.~A. Rossiter, B.~Kouvaritakis, J.~R. Gossner, Guaranteeing feasibility in
  constrained stable generalized predictive control., IEEE Proc. Control theory
  Appl. 143 (1996) 463--469.

\bibitem{RawlingsLIB09}
J.~B. Rawlings, D.~Q. Mayne, Model Predictive Control: Theory and Design, 1st
  Edition, Nob-Hill Publishing, 2009.

\bibitem{RaoAIChE99}
C.~V. Rao, J.~B. Rawlings, Steady states and constraints in model predictive
  control, AIChE Journal 45~(6) (1999) 1266--1278.

\bibitem{FerramoscaAUT09}
A.~Ferramosca, D.~Limon, I.~Alvarado, T.~Alamo, E.~F. Camacho, {MPC} for
  tracking with optimal closed-loop performance, Automatica 45~(8) (2009)
  1975--1978.

\bibitem{AlvaradoPHD07}
I.~Alvarado, Model predictive control for tracking constrained linear systems,
  Ph.D. thesis, Univ. de Sevilla. (2007).

\bibitem{FerramoscaPHD11}
A.~Ferramosca, Model predictive control for systems with changing setpoints,
  Ph.D. thesis, Univ. de Sevilla.,
  http://fondosdigitales.us.es/tesis/autores/1537/ (2011).

\bibitem{BoydLIB06}
S.~Boyd, L.~Vandenberghe, Convex Optimization, Cambridge University Press,
  2006.

\bibitem{Nocedal99numerical}
J.~Nocedal, S.~J. Wright, Numerical Optimization, Vol.~2, Springer New York,
  1999.

\bibitem{GilbertTAC91}
E.~G. Gilbert, K.~Tan, Linear systems with state and control constraints: The
  theory and application of maximal output admissible sets, IEEE Transactions
  on Automatic Control 36 (1991) 1008--1020.

\bibitem{BemporadAUT02}
A.~Bemporad, M.~Morari, V.~Dua, E.~Pistikopoulos, The explicit linear quadratic
  regulator for constrained systems, Automatica 38 (2002) 3--20.

\bibitem{GrimmAUT04}
G.~Grimm, M.~J. Messina, S.~Z. Tuna, A.~R. Teel, Examples when nonlinear model
  predictive control is nonrobust, Automatica 40 (2004) 1729--1738.

\bibitem{LimonISS09}
D.~Limon, T.~Alamo, D.~M. Raimondo, D.~M. de~la Pe{\~n}a, J.~M. Bravo,
  A.~Ferramosca, E.~F. Camacho, Input-to-state stability: an unifying framework
  for robust model predictive control, in: L.~Magni, D.~M. Raimondo,
  F.~Allg{\"o}wer (Eds.), Nonlinear Model Predictive Control - Towards New
  Challenging Applications, Springer, 2009, pp. 1--26.

\bibitem{MessinaAUT05}
M.~J. Messina, S.~E. Tuna, A.~R. Teel, Discrete-time certainty equivalence
  output feedback: allowing discontinuous control laws including those from
  model predictive control., Automatica 41 (2005) 617--628.

\bibitem{MayneAUT00}
D.~Q. Mayne, J.~B. Rawlings, C.~V. Rao, P.~O.~M. Scokaert, Constrained model
  predictive control: Stability and optimality, Automatica 36~(6) (2000)
  789--814.

\bibitem{HuTAC02}
B.~Hu, A.~Linnemann, Towards infinite-horizon optimality in nonlinear model
  predictive control, IEEE Transactions on Automatic Control 47~(4) (2002)
  679--682.

\bibitem{LuenbergerLIB84}
D.~E. Luenberger, Linear and Nonlinear Programming, Addison-Wesley, 1984.

\bibitem{FerramoscaIJSS11}
A.~Ferramosca, D.~Limon, I.~Alvarado, T.~Alamo, F.~Casta{\~n}o, E.~F. Camacho,
  Optimal {MPC} for tracking of constrained linear systems, Int. J. of Systems
  Science 42~(8) (2011) 1265--1276.

\bibitem{jones2006primal}
C.~N. Jones, J.~Maciejowski, Primal-dual enumeration for multiparametric linear
  programming, Mathematical Software-ICMS 2006 (2006) 248--259.

\bibitem{morari2008multiparametric}
M.~Morari, C.~N. Jones, M.~N. Zeilinger, M.~Baric, Multiparametric linear
  programming for control, in: Control Conference, 2008. CCC 2008. 27th
  Chinese, IEEE, 2008, pp. 2--4.

\bibitem{JohanssonTCST00}
K.~H. Johansson, The quadruple-tank process: A multivariable laboratory process
  with an adjustable zero, IEEE Transaction on Control Systems Technology 8
  (2000) 456--465.

\bibitem{VidyasagarLIB93}
M.~Vidyasagar, Nonlinear Systems Theory, 2nd Edition, Prentice-Hall, 1993.

\end{thebibliography}

\end{document}